\documentclass[twoside,11pt]{article}
\usepackage{amsfonts}
\usepackage{ulem}
\usepackage{graphicx}
\newtheorem{thm}{Theorem} 
\newtheorem{cor}[thm]{Corollary} 
\newtheorem{lemma}[thm]{Lemma} 
\newtheorem{prop}[thm]{Proposition}

\newtheorem{rem}[thm]{Remark}

\begin{document} 
 
 \def\reals{\hbox{\sl I\kern-.18em R \kern-.3em}} 
 \def\ch{{\cal H}}
 \def\cA{{\cal A}}
 \def\cB{{\cal B}}
 \def\cK{{\cal K}}
 \def\cC{{\cal C}}
 \def\cH{{\cal H}}
 \def\cE{{\cal E}}
 \def\cR{{\cal R}}
 \def\cS{{\cal S}}
 \def\cT{{\cal T}}
 \def\cU{{\cal U}}
 \def\cW{{\cal W}}
 \def\cX{{\cal X}}
 \def\cI{{\cal I}}
 \def\cO{{\cal O}}
 \def\axis{{\bf A}}
 \def\fib{{\bf H}}
 \def\ba{{\bf a}}
 \def\bb{{\bf b}}
 \def\bc{{\bf C}}
 \def\be{{\bf e}}
 \def\sa{{\bf sa\ }}
 \def\scc{{\bf scc\ }}
 \def\UU{{{{\mathbf{X}}}^{\e \cdot 2}}}
 \def\arcpres{{X^\eta}}
 \def\arcpresprime{{X^{\prime\eta}}}
 \def\sheararcpres{{X^\eta_{\cI}}}
 \def\tilesheararcpres{{X^{\kern-1.05emT \kern.46em\eta}_{\cI}}}
 \def\sheararcpresprime{{X^{\prime \eta}_{\cI}}}
 \def\sheararcpresprimeprime{{X^{\prime \prime \eta}_{\cI}}}
 \def\notchdisc{{\Delta^\eta_\v}}
 \def\notchdiscd{{\Delta^\eta_d}}
 \def\notchdisce{{\Delta^\eta_e}}
 \def\notchdiscf{{\Delta^\eta_f}}
 \def\tilenotchdisc{{\Delta^{\kern-.98emT \kern.4em\eta}_\v}}
 \def\tilenotchdiscd{{\Delta^{\kern-.98emT \kern.4em\eta}_d}}
 \def\tilenotchdisce{{\Delta^{\kern-.98emT \kern.4em\eta}_e}}
 \def\tilenotchdiscf{{\Delta^{\kern-.98emT \kern.4em\eta}_f}}
 \def\notchdiscprime{{\Delta^{\prime \eta}_\v}}
 \def\notchdiscprimed{{\Delta^{\prime \eta}_d}}
 \def\notchdiscprimee{{\Delta^{\prime \eta}_e}}
 \def\notchdiscprimef{{\Delta^{\prime \eta}_f}}
 \def\hn{{h^\eta}}
 \def\vn{{v^\eta}}
 \def\an{{\a^\eta}}
 \def\bu{{\dot \cup}}
 \def\ver{{\bf v}}
 \def\sin{{\bf s}}
 \def\d{{\delta}} 
 \def\ci{{\circ}} 
 \def\e{{\epsilon}} 
 \def\l{{\lambda}} 
 \def\L{{\Lambda}} 
 \def\m{{\mu}} 
 \def\n{{\nu}} 
 \def\o{{\omega}} 
 \def\s{{\sigma}} 
 \def\v{{\varepsilon}} 
 \def\a{{\alpha}} 
 \def\b{{\beta}} 
 \def\p{{\partial}} 
 \def\r{{\rho}} 
 \def\ra{{\rightarrow}} 
 \def\lra{{\longrightarrow}}
 \def\g{{\gamma}} 
 \def\D{{\Delta}} 
 \def\La{{\Leftarrow}} 
 \def\Ra{{\Rightarrow}} 
 \def\x{{\xi}} 
 \def\c{{\mathbb C}} 
 \def\z{{\mathbb Z}} 
 \def\2{{\mathbb Z_2}} 
 \def\q{{\mathbb Q}}
 \def\r{{\mathbb R}} 
 \def\t{{\tau}} 
 \def\u{{\Upsilon}}
 \def\U{{\Upsilon}}
 \def\th{{\theta}} 
 \def\la{{\leftarrow}}
 \def\ra{{\rightarrow}}
 \def\lla{{\longleftarrow}} 
 \def\da{{\downarrow}} 
 \def\ua{{\uparrow}} 
 \def\nwa{{\nwtarrow}} 
 \def\swa{{\swarrow}} 
 \def\nea{{\netarrow}} 
 \def\sea{{\searrow}} 
 \def\hla{{\hookleftarrow}} 
 \def\hra{{\hookrightarrow}}
 \def\rhu{{\rightharpoonup}}
 \def\ntran{{\longrightarrow^{\kern-1emN \kern.37em }}}
 \def\btran{{\longrightarrow^{\kern-1emB \kern.37em }}}
 \def\stran{{\longrightarrow^{\kern-1emS \kern.37em }}}
 \def\Ttran{{\longrightarrow^{\kern-1emT \kern.37em }}}
 \def\qed{{\hfill$\diamondsuit$}} 
 \def\pf{{\noindent{\bf Proof.\hspace{2mm}}}} 
 \def\ni{{\noindent}} 
 \def\sm{{{\mbox{\tiny M}}}} 
 \def\sc{{{\mbox{\tiny C}}}} 
 \def\ke{{\mbox{ker}(H_1(\p M;\2)\ra H_1(M;\2))}} 
 \def\et{{\mbox{\hspace{1.5mm}}}}
 \def\Xb{{X_{\overline{\beta}}}}
 \def\ob{{\overline{\beta}}}
 \def\bi{\begin{itemize}}
 \def\ei{\end{itemize}}

%
%

\pagestyle{myheadings}
\markboth{William W. Menasco}{Recognizing destabilization, exchange moves and flypes}
\title{Recognizing when closed braids admit
a destabilization, an exchange move or an elementary flype} 

 
\author{William W. Menasco
\\  
{\small University at Buffalo} 
 \\ {\small Buffalo, New York 14260-2900}
 \\ {\small {\em Dedicated to Joan Birman on the occasion of her advancement to
active retirement.} } 
}
\date{revised  June 24, 2010
}

\maketitle

\begin{abstract}
The Markov Theorem Without Stabilization (MTWS) established the existence
of a calculus of braid isotopies that can be used to move between closed braid
representatives of a given oriented link type without having to increase the
braid index by stabilization \cite{[BM4]}.  Although the calculus is extensive there are three
key isotopies that were identified and analyzed---destabilization, exchange moves
and elementary braid preserving flypes.  One of the critical open problems left in
the wake of the MTWS is the {\em recognition problem}---determining when a given
closed $n$-braid admits a specified move of the calculus.
In this note we give an algorithmic solution to the recognition problem for
these three key isotopies of the MTWS calculus.
The algorithm is ``directed'' by a complexity measure
that can be {\em monotonically simplified} by the application of {\em elementary moves}.
\end{abstract}
 
\section{Introduction.}
\label{section:introduction}
\subsection{Preliminaries.}
\label{subsection:preliminaries}
Given an oriented link $ X \subset S^3$ it is a classical result of Alexander \cite{[A]}
that $X$ can be represented as a closed $n$-braid.  For expository purposes
it is convenient to translate Alexander's result into the following setting.
Let $S^3 = \mathbb R^3 \cup \{ \infty \} $ and give $\mathbb R^3$ an open-book decomposition, i.e.
$\mathbb R^3 \setminus \{z-{\rm axis}\}$ is fibered by a collection of half-plane fibers
$\fib = \{ H_\theta | \theta \in [0,2\pi]/0 \sim 2\pi \ (= S^1) \} $ where
$\axis = \partial \bar{H}_\theta$ is the
$z$-axis.  Equivalently, we consider the cylindrical coordinates
$(r,\theta,z)$ on $\mathbb R^3$ and $H_\theta$ is the set of points in $\mathbb R^3$ having their
second coordinate the fixed $\theta \in S^1$.

An oriented link $X$ in $\mathbb R^3 ( \subset S^3 )$ is a {\em closed $n$-braid} if
$X \subset \mathbb R^3 \setminus \{z-{\rm axis}\}$ such that $X$ transversely intersects
each fiber of $\fib$ in $n$ points.  The {\em braid index of $X$} is
the cardinality $b(X) = |X \cap H_\theta| = n$ which is invariant
for all $H_\theta \in \fib$.

\begin{figure}[htpb]
\centerline{\includegraphics[scale=.4, bb=27 402 589 751]{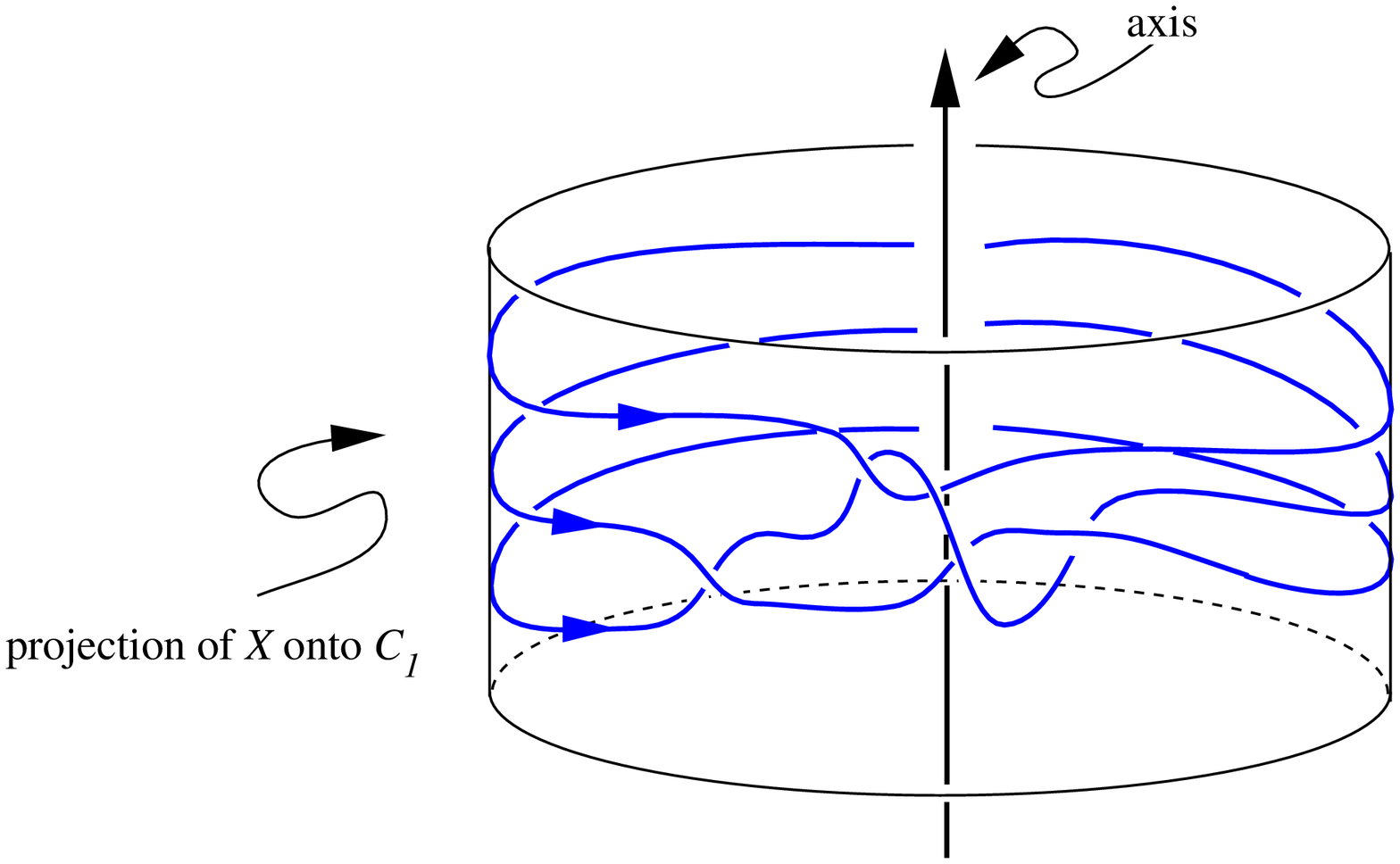}}
\caption{}
\label{figure:braid on cylinder}
\end{figure}

\begin{figure}[htpb]
\centerline{\includegraphics[scale=.5, bb=50 727 927 967]{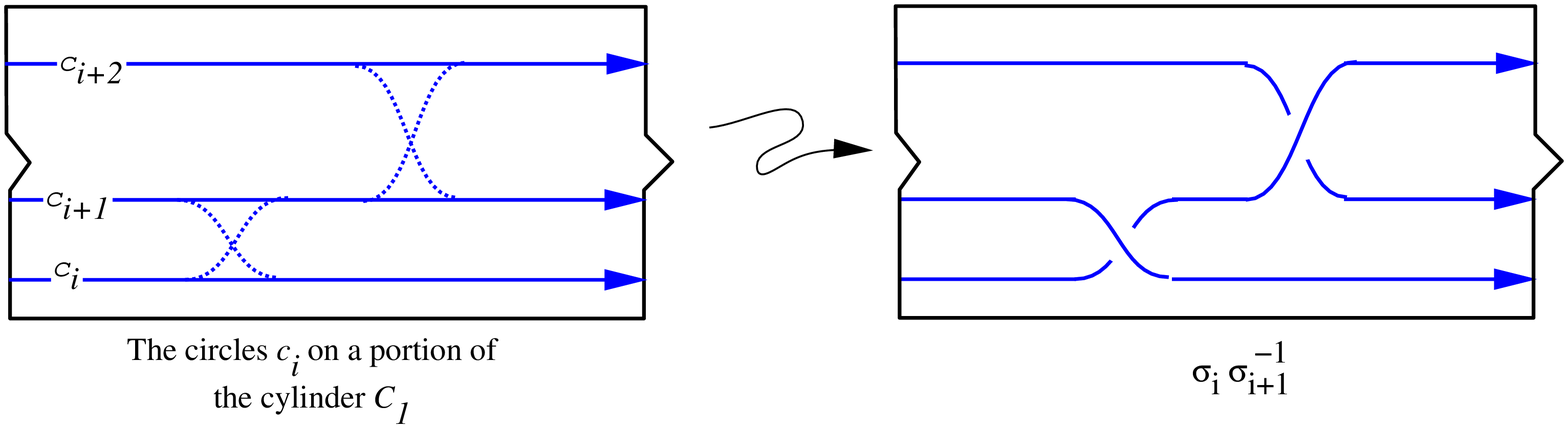}}
\caption{}
\label{figure:standard projection}
\end{figure}

Possibly after a small isotopy of $X$ in $\reals^3 \setminus \{z-{\rm axis}\}$, 
we can consider a regular projection $\pi: X \rightarrow C_1$ given by
$\pi:(r,\theta,z) \mapsto (1,\theta,z)$, where $C_1 = \{ (r,\theta,z) | r=1\}$.  The projection
$\pi(X) \subset C_1$ is isotopic to a {\em standard projection} that is constructed as follows.
For $n = b(X)$ we first consider the circles:
$c_i = \{ (r,\theta,z) | r=1 , z=\frac{i}{n} \}$, $1\leq i \leq n$.  Next, we alter the
projection of this trivial unlink of $n$ components to construct the projection
of $\pi(X)$ by having adjacent circles $c_i$ and $c_{i+1}$ cross via the addition of
a positive or negative crossing at the needed angle to produce a projection of $X$.
It is clear through an ambient isotopy of $\mathbb R^3$ which preserves the fibers of $\fib$
that we can always reposition $X$ in $\mathbb R^3 \setminus \{z-{\rm axis}\}$ 
so that $\pi(X)$ is a standard projection.  (See Figures \ref{figure:braid on cylinder}
and \ref{figure:standard projection}.)

From a standard projection $\pi(X)$ on $C_1$ we can read off a cyclic word in the classical
Artin generators $\sigma^{\pm 1}_i$, $1 \leq i \leq n-1 $.
Specifically, $\b (X)$ is the cyclic
word that comes from recording the angular occurrence of crossings where a positive (respectively
negative) crossing between the circles $c_i$ and $c_{i+1}$ contributes a $\sigma_i$ (respectively
$\sigma_i^{-1}$) element to $\b(X)$.

\begin{figure}[htpb]
\centerline{\includegraphics[scale=.45, bb=7 626 1139 1166]{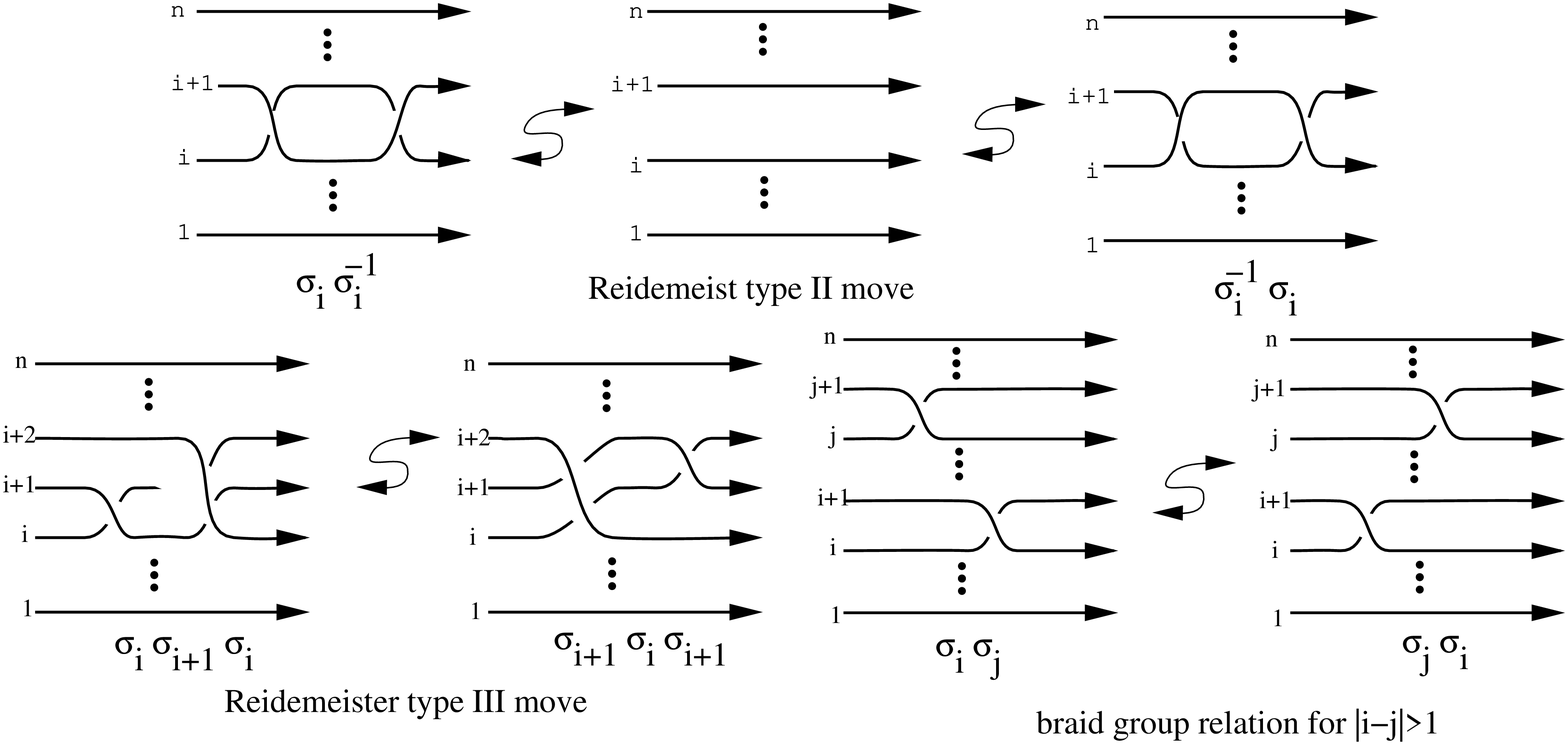}}
\caption{Generally, type III moves correspond to
$\sigma_i^\e \sigma_{i+1}^{\d} \sigma_i^{\d}$ interchanged with
$ \sigma_{i+1}^{\d} \sigma_i^{\d} \sigma_{i+1}^{\e}$ where $\e,\d \in \{-1 , 1\}$. Generally, the braid group
relation is $\sigma_i^{\e} \sigma_j^{\d}$ interchanged with
$\sigma_j^{\d} \sigma_i^{\e}$ where again $\e,\d \in \{-1 , 1\}$.}
\label{figure:braid isotopies}
\end{figure}

We can alter $\pi(X)$ (and, correspondingly, $\b(X)$) while preserving the $n$-braid structure
of $X$, plus its standard projection characteristic, by use of type-II and -III
Reidemeister moves and the braid group relation.  (See Figure \ref{figure:braid isotopies}.)

We consider the equivalence classes under the moves in Figure \ref{figure:braid isotopies}.
Specifically, $X$ and $X^\prime$ are {\em braid isotopic} if $\pi(X)$ can be
altered to produce $\pi(X^\prime)$ through a sequence of type-II and -III moves, and
braid group relations.  We will let $\cB_n(X)$ be notation for the equivalence class of $n$-braids
which are braid isotopic to $X$.

Next, let $\cW^t$ be all words generated by the set $\{\sigma_1^{\pm 1}, \cdots , \sigma_t^{\pm 1} \}$ and
$\cU^s$ be all words generated by the set $\{\sigma_s^{\pm 1}, \cdots , \sigma_{n-1}^{\pm 1} \}$.
Initially we allow for the possibilities of $\cW^t \cap \cU^s = \emptyset$ or
$\cW^t \cap \cU^s \not= \emptyset$, but exclude $\cW^t = \cU^s$.
We now have a sequence of definitions.

\begin{figure}[ht]
\centerline{\includegraphics[scale=.8, bb=0 0 488 366]{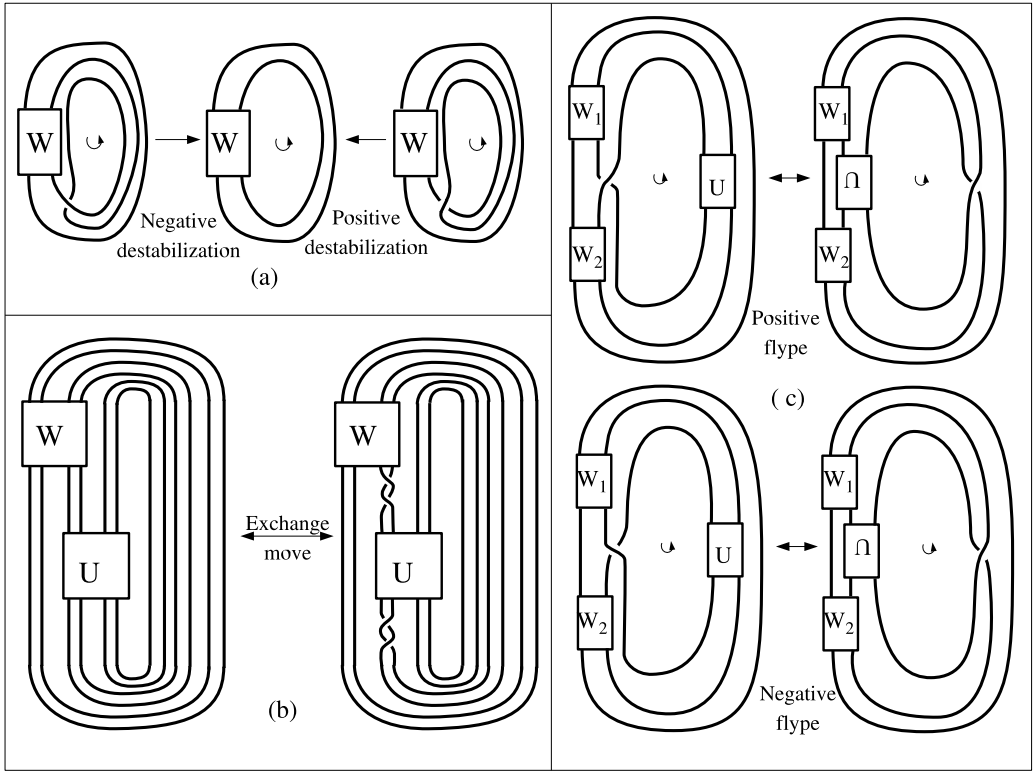}}
\caption{Destabilization, an exchange move and braid preserving flypes are
depicted respectively in illustrations (a), (b) and (c).  In particular, in (b) we have
illustrated an exchange move on a $6$-braid so as to emphasize the occurrence of the
full twists $\tau_{[3,4]}^{\ }$ \& $\tau_{[3,4]}^{-1}$.  The block labels correspond to
the associated braid syllable.  In particular, in (c) when $U = \sigma_{n-1}^{p}$ we have an
elementary flype.  The labels in the braid boxes correspond to the associated braid syllable.}
\label{figure:moves}
\end{figure}

An $n$-braid $X$ {\em admits a destabilization} if for $\pi(X)$, its associated braid
word $\b(X)$ is of the form $W \sigma_{n-1}^{\pm 1}$, where $W \in \cW^{n-2}$.
(See Figure \ref{figure:moves}(a).)
The $(n-1)$-braid $Y$ with $\b(Y) = W$
(and corresponding projection $\pi(Y)$) is obtained by a {\em destabilization} of $X$.

Consistent with the terminology in \cite{[M1]}, an $n$-braid $X$ {\em admits a double destabilization} for
if for $\pi(X)$, its associated braid word $\b(X)$ is of the form
$$W \UU = W \s^{\e}_{n-2} \s_{n-1}^{\e} \s_{n-3}^{\e} \s_{n-2}^{\e} ,$$
where $W\in \cW^{n-3}$ and $\e \in \{ -1 , +1\}$.  (Thus, $\UU$ is a positive ($\e = +1$) or negative
($\e = -1$) ``$2$-strand crossing''.)
%
%

$X$ {\em admits an exchange move} if $\b(X)$ is of the form
$WU$ where $W$ (respectively $U$) is a word of $\cW^t$
(respectively $\cU^s$) for some integers $s \leq t$.
The $n$-braid $Y$ that is associated with the cyclic word
$W \tau_{[s,t+1]^{\ }} U \tau_{[s,t+1]}^{-1}$, where
$\tau_{[s,t+1]}$ is a full (positive or negative) twist on the $s$ through $(t+1)$ strands,
is {\em exchange related} to $X$.  (See Figure \ref{figure:moves}(b).)

We further refine our concept of exchange move. Suppose $X$ admits an exchange move. Then $\b(X)$ is
of the form $WU$ where $W$ is in $\cW^t$ and $U$ is in $\cU^s$.
The reader should observe that if $W$ is a word in $\cW^t$ then it is also a word
in $\cW^{t^\prime}$ for any $ t < t^\prime \leq n-2$.  Similarly,
$U$ is also a word in $\cU^{s^\prime}$ for $2 \leq  s^\prime < s$.
Since the specification of $t$ and $s$ determine which strands are to be used in the full twist
$\tau_{[s,t+1]}$ it is useful to have a canonical way of viewing the words $W$ and $U$.
Thus, we say the product $WU = \b(X)$ corresponds to a {\em thin exchange move} if
for all choices of $s$ and $t$ with $s \leq t$ we have $t - s$ minimal.  In Figure \ref{figure:moves}(b)
we have $s = 3 = t$, so this illustration corresponds to a thin exchange move.
It is clear that if $X$ admits an
exchange move then it has a thin exchange move.

Continuing, $\b(X)$ {\em admits a (braid preserving) flype} if it is of the form
$W_1 U W_2 \sigma_{n-1}^{\pm 1}$ where $W_1 ,  W_2 \in \cW^{n-2}$ and $U \in \cU^s$
for some integers $s \leq n-1$.  When
$s = n-1$ we have $\b(X)$ of the form $W_1 \sigma_{n-1}^p W_2 \sigma_{n-1}^{\pm 1}$
where $p \in \mathbb Z - \{ 0 \}$ and $\b(X)$ {\em admits an elementary (braid preserving) flype}.
The $n$-braid $Y$ which has
$ \b(Y) = W_1 \sigma_{n-1}^{\pm 1} W_2 \sigma_{n-1}^{p}$ is {\em elementary flype related} to $X$.

We say that $\cB_n(X)$ {\em admits a destabilization, exchange move,
flype, or elementary flype},
respectively, if there exists a braid representative $X^\prime \in \cB_n(X)$ which
admits a destabilization, exchange move, flype, or elementary flype, respectively.
(As mentioned, if $\cB_n(X)$ admits an exchange move then there will be a $X^{\prime\prime} \in \cB_n(X)$ that will admit a thin exchange move.)

\begin{rem}
\label{remark:elementary exchange moves and flypes}
{\rm The reader should notice that when $s = t+1$ then the product $WU$ represents a
composite link.  When $s \geq t+2$ then the product $WU$ is a split link.
We will assume that $\b(X)$ is neither composite nor split.
The reader should further observe that when $\b(X)$ is of the form $WU$ with
$W \in \cW^t$, $U \in \cU^s$ and $t-s < n-3$ then $\b(X)$ admits infinitely many
distinct exchange moves.  To take an illustrative construction, suppose we have a $n$-braid
$X$ with $\b(X) = WU$ where with $W \in \cW^t$ and $ U \in \cU^s$ and $ s > 2$.
Let $\a$ be any braid word in $\cW^{s-1}$. Then for the cyclic word $WU$
we have
$$ WU = W \a \a^{-1}U = W \a U \a^{-1} = 
[\a^{-1} W \a] U = W^\prime U^\prime $$ where still $W^\prime \in \cW^t$ but
$U^\prime \in \cU^s$.
Now viewing $W^\prime \in \cW^t$ and $U^\prime \in \cU^{s-1}$ we perform
an exchange move on $W^\prime U^\prime$ to
obtain $W^\prime \tau_{[s-1,t+1]}^{\ } U^\prime \tau_{[s-1,t+1]}^{-1}$.
For a judicious choose of $W$, $U$ and $\a$ we can insure that
$W^\prime \tau U^\prime \tau^{-1}$ is not braid isotopic to $W \tau U \tau^{-1}$.  Thus,
by varying $\a$---for example, taking powers of a fixed $\a$---we can produce infinitely many exchange moves.
From this example it becomes clear that the exchange moves that will always ``respect braid
isotopy'' are the thin exchange moves.  That is, when performing a thin
exchange move we are requiring that the full twist $\tau$ be on the least number
of strands.  In our results this will have the implication that we will be able
to recognize when two braids $X$ and $Y$ are related by an thin exchange move.
Finally, we observe that this phenomenon is not an issue for elementary flypes.
As with exchange moves, it is true that we have a similar sequence of cyclic equalities
$$  W_1 \sigma_{n-1}^p W_2 \sigma_{n-1}^{\pm 1}
= W_1 \a \a^{-1} \sigma_{n-1}^p W_2 \sigma_{n-1}^{\pm 1} =
[W_1 \a] \sigma_{n-1}^p [\a^{-1} W_2] \sigma_{n-1}^{\pm 1} =
W_1^\prime \sigma_{n-1}^p W_2^\prime \sigma_{n-1}^{\pm 1}$$
for any word $\a$ using generators in  $\cW^{n-3}$.  But, after the elementary flype 
we also have
$$ W_1^\prime \sigma_{n-1}^{\pm 1} W_2^\prime \sigma_{n-1}^p =
[W_1 \a] \sigma_{n-1}^{\pm 1} [\a^{-1} W_2] \sigma_{n-1}^p = 
W_1 \a \a^{-1} \sigma_{n-1}^{\pm 1} W_2 \sigma_{n-1}^p =
W_1 \sigma_{n-1}^{\pm 1} W_2 \sigma_{n-1}^p.$$
Thus, although there are infinitely
many words admitting the same elementary flype, up to braid isotopy all of these
elementary flypes relate to the same two closed $n$-braids. \qed}
\end{rem}

\begin{rem}
\label{remark: equivalency class for EM and F}
{\rm We also observe that having $\b(X) = WU$ with
$W \in \cW^t$, $U \in \cU^s$ and $s > 2$ is not a unique format illustrating
even a thin exchange move.  Specifically, suppose $W = W^\prime \alpha$ with
$W^\prime \in \cW^t$ and $\alpha \in \cW^{n-s-1}$.  Then $\alpha$ and $U$ commute
in the braid group, and the closed braids having cyclic words
$[W \alpha] U$ and $[\alpha W] U$ are braid isotopic.
Moreover, the closed braids having cyclic words
$[W \alpha] \tau_{[s , t+1]}^{\ } U \tau_{[s, t+1]}^{-1}$ and
$[\alpha W] \tau_{[s , t+1]}^{\ } U \tau_{[s, t+1]}^{-1}$ are braid isotopic.
Thus, it is convenient to equate $[W \alpha] U$ and $[\alpha W] U$
as admitting the same exchange move.
\qed}
\end{rem}

\begin{rem}
\label{remark: double destabilization}
{\rm We can apply the discussion in the previous remarks to
an $n$-braid $\b(X) = W \UU$ where $\UU$ is a $2$-strand crossing so as to
have $X$ admitting a double destabilization.
Notice that $X= W\UU$ also admits an exchange move and is exchange equivalent to
$X^\prime = W \tau_{[(n - 2),(n -1)]}^{p} \UU \tau_{[(n-2),(n - 1)]}^{-p}$
for $p \in \mathbb{Z}$.
Now the reader should notice two features.  First, since $W\UU$ has $2$ parallel strands
we can slide the $\tau_{[(n - 2),(n -1)]}^{p}$-twist through these stands to cancel the
$\tau_{[(n-2),(n - 1)]}^{-p}$-twists.  Thus, $X$ and $X^\prime$ are braid isotopic.
Second, by destabilizing a single strand of the $2$-strands of $\UU$ we have a means by which
to produce infinitely distinct conjugacy classes.
This phenomena of $X$ destabilizing to infinitely
many possible distinct conjugacy classes was investigated by
A. V. Malyutin \cite{[M1], [M2]}.  (The author wishes to thanks
the referee for alerting him to this body of work.)
\qed}
\end{rem}

A long standing problem (Problem 1.84 in \cite{[K]}) is determining when an $n$-braid
equivalence class $\cB_n(X)$ contains a braid $X$ that admits either a destabilization,
an exchange move, or an elementary flype.  Our main result (Theorem \ref{theorem:main result})
states that there is a simple algorithmic method for making these determinations.
To understand this algorithm we need to consider our representation of $X$ in $\fib$ anew using
`rectangular diagrams'.

\noindent
{\bf Notational conventions--}
%
Our discussion will also require extensive use of arcs joined together at common endpoints
to form edge-paths that are homeomorphic to either a closed interval or circle.
We thus introduce notation for {\em ordered union} and {\em cyclic ordered union}.
So, for arcs $\{a_1 , \cdots , a_l \}$ with $a_{i}$ sharing a single endpoint with $a_{i+1}$,
$1 \leq i < l$,
we have the ordered union $a_1 \bu a_2 \bu \cdots \bu a_l$ being the edge-path homeomorphic to
a closed interval that is obtained by adjoining $a_i$ to $a_{i+1}$ 
at their common endpoint, $1 \leq i < l$.  If
$a_1$ and $a_l$ also share a common endpoint then the cyclic ordered union
$a_1 \bu \cdots \bu a_l \bu$ is homeomorphic to a circle.

Finally, in our use of
angle parameters $\theta$ we have identify
$S^1$ with the quotient space $[0,2\pi]/ 0 \sim 2\pi$.  Since $S^1$ is given the orientation
that corresponds to winding positively around the $z -{\rm axis}$, the points on $S^1$ are
cyclically ordered and it makes sense to talk about the oriented angle interval $[\theta_1 , \theta_2] \subset S^1$
that starts at $\theta_1$ and ends at $\theta_2$.  The arc length $|[\theta_1 , \theta_2]|$ will
necessarily be less than $2 \pi$.

\subsection{Rectangular diagrams and main results.}
\label{subsection:rectangular diagrams}

A {\em horizontal arc}, $h \subset C_1$,
is any arc having parametrization $\{ (1,t,z_0 ) | \ t \in
[\theta_1 , \theta_2] \}$.  The {\em horizontal position} of $h$ is the fixed constant $z_0$.
The {\em angular support} of $h$ is the angle interval
$[\theta_1, \theta_2] \subset S^1$.  Horizontal arcs inherit a
natural orientation from the forward direction of the $\theta$ coordinate.
A {\em vertical arc}, $v \subset H_{\theta_0^{\ }}$,
is any arc having parametrization $\{ (r(t), \theta_0 , z(t) ) | \ 
0 \leq t \leq 1, \ r(0)=r(1)=1; \  {\rm and} \  
r(t) > 1 , \  \frac{dz}{dt} > 0 {\rm \ for} \ t \in (0,1) \}$, where
$r(t)$ and $z(t)$ are real-valued functions that are continuous on $[0,1]$ and differentiable on
$(0,1)$.  The {\em angular position} of $v$ is $\theta_0$.
The {\em vertical support} of $v$ is the interval $[z(0),z(1)]$.  (We remark that the
parametrization of the the vertical arcs will not be used in assigning
orientation to the vertical arcs.)

\begin{figure}[htpb]
\centerline{\includegraphics[scale=.4, bb=19 408 991 1021]{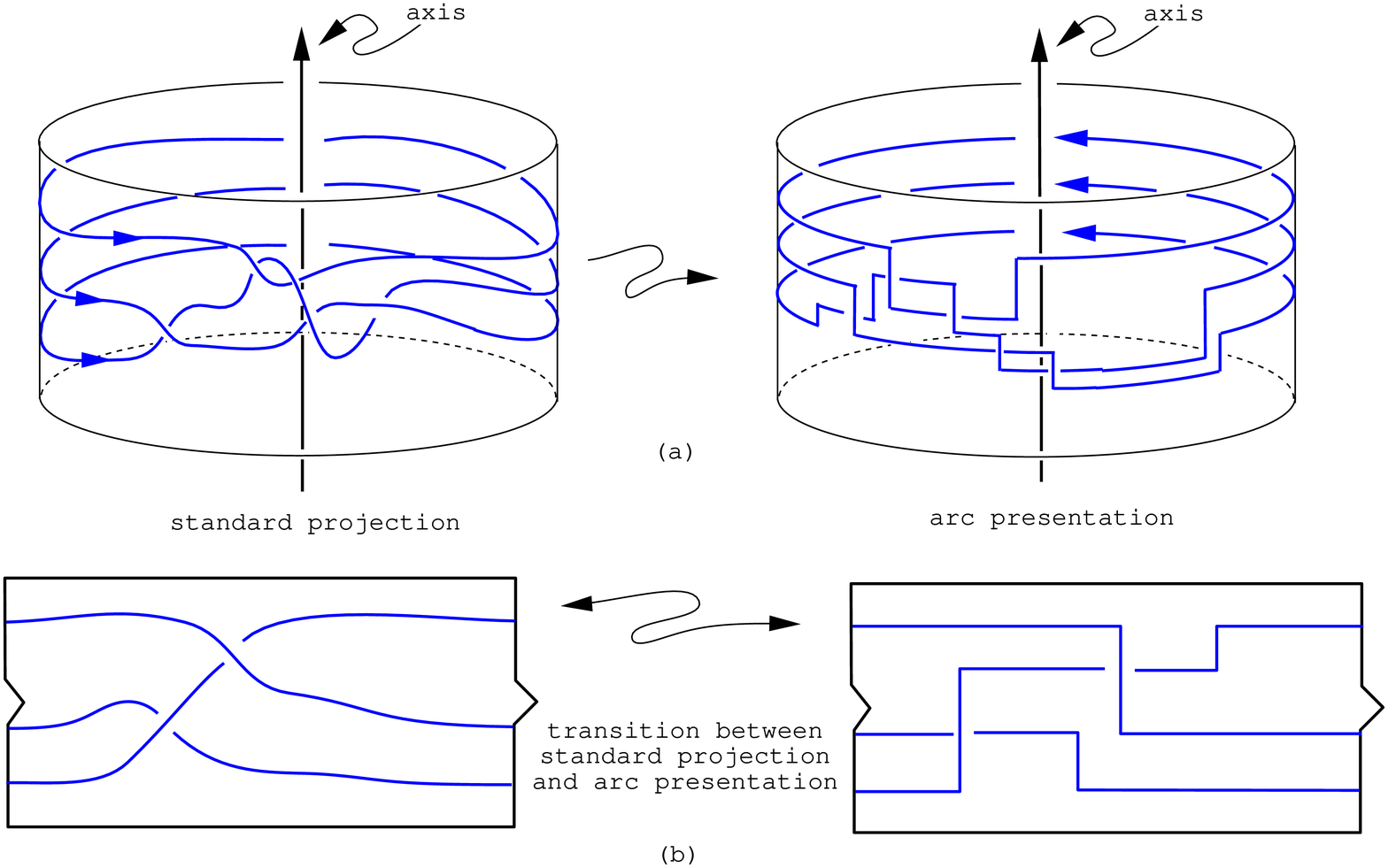}}
\caption{}
\label{figure:arc presentation}
\end{figure}

Let $\cX$ be an oriented link type in $S^3$.
$\arcpres \in \cX$ is an {\em arc presentation}
if each component $Y$ of $\arcpres$ is a cyclic union of arcs
$h_1^Y \bu v_1^Y \bu \cdots \bu h_k^Y \bu v_k^Y \bu$
with $k$ necessarily varying in value between components and:
\bi
\item[1.] each $h_i^Y$, $1 \leq i \leq k$;
is an oriented {\em horizontal arc} having orientation agreeing with the $Y$,
\item[2.] each $v_i^Y$, $1 \leq i \leq k$, is a {\em vertical arc} having orientation agreeing
with $Y$;
\item[3.] $h_i^Y \cap v_j^Y \subset \partial h_i^Y \cap \partial v_j^Y$, and
for $1 \leq i \leq k$ this intersection is a single point when
$j \ ({\rm mod} k) = \{i, i-1 \}$, otherwise it is empty;
\item[4.] the horizontal position of each horizontal arc is distinct over all components of $\arcpres$;
\item[5.] the angular position of each vertical arc is distinct over all components of $\arcpres$.
\ei

For a given arc presentation $\arcpres $
there is a cyclic order to the horizontal positions of the $h_i {\rm 's}$, as determined
by their occurrence on the $\axis$, and a cyclic
order to the angular position of the $v_j {\rm 's}$, as determined by their occurrence in $\fib$.
It is clear that given two arc presentations
with identical cyclic order for horizontal positions and angular positions there is an ambient isotopy
between the two presentations that corresponds to re-scaling of the horizontal positions and
angular positions along with the vertical and angular support of the arcs in the presentations.
Thus, we will think of two arc presentations as being equivalent if the cyclic orderings of their horizontal positions and vertical positions are equivalent. 

We define the {\em complexity} of the arc
presentation $\arcpres$, $\cC(\arcpres)$, as being the number of vertical arcs.

Given a closed $n$-braid $X$ and a corresponding standard projection
$\pi(X)$ we can easily produce a (not necessarily unique)
arc presentation $\arcpres$ as illustrated
by the transition in Figure \ref{figure:arc presentation}.  Clearly, there is also the
transition from an arc presentation $\arcpres$ to an $n$-braid $X$ with a standard
projection $\pi(X)$.  We will use the notation $X \ntran \arcpres$ or $\arcpres \btran X$ it indicate
these two presentation transitions.

\begin{figure}[htpb]
\centerline{\includegraphics[scale=.5, bb=0 0 845 366]{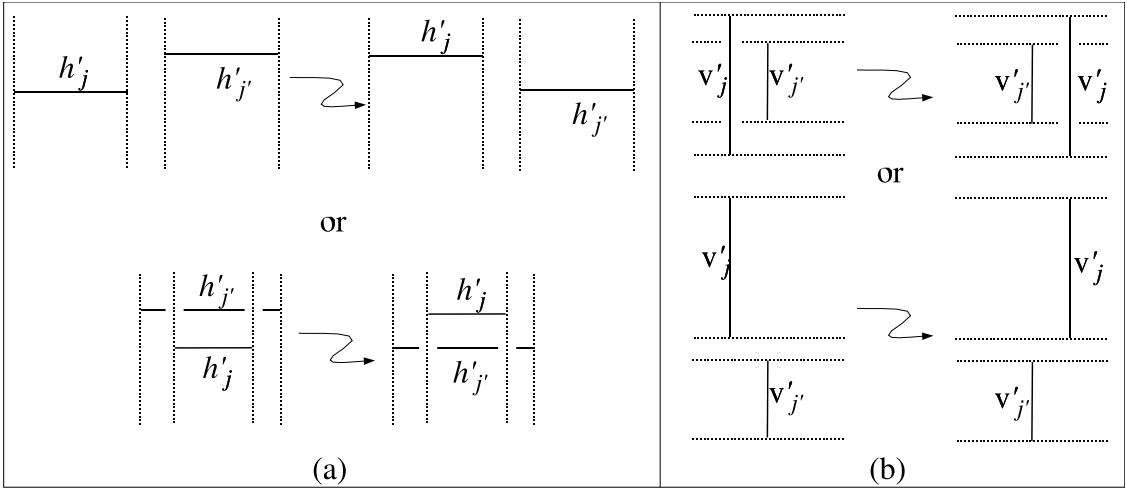}}
\caption{Illustration (a) corresponds to horizontal exchange moves and (b) corresponds
to vertical exchange moves.  In each of the illustrated sequences dotted arcs are
used to indicate that the adjoining vertical arcs (for (a)) and horizontal arcs
(for (b)) have two possible ways of attaching themselves to the labeled solid
arc.}
\label{figure:elementary moves}
\end{figure}

We next define {\em elementary moves} on an arc presentation $\arcpres$.
To setup these moves we let 
$$ \arcpres = \bigsqcup_{Y \in X} \left( h_1^Y \bu v_1^Y \bu \cdots \bu h_k^Y \bu v_{k^Y}^Y \bu \right) $$
with: $z_i^Y$ the horizontal position for $h_i^Y$;
$[\theta_{i_1^{\ }}^Y,\theta_{i_2^{\ }}^Y] \subset S^1$ the angular
support for $h_i^Y$; $\theta_i^Y$ the angular position for $v_i^Y$; and $[z_{i_1^{\ }}^Y , z_{i_2^{\ }}^Y]$ the
vertical support for $v_i^Y$ (where $1\leq i \leq k^Y$ in all statements).  It is convenient at times to
drop the uses of the $Y$ superscript.  The reader should also be alert to the use of the subscript index variable to
indicate needed associations, e.g. $h_i$ has angular support
$[\theta_{i_1^{\ }},\theta_{i_2^{\ }}]$
whereas $h_{j}$ has angular support $[\theta_{j_1^{\ }},\theta_{j_2^{\ }}]$.)

\noindent \underline{Horizontal exchange move}---The move comes in two
flavors.  The first allows us to take the horizontal arc of $\arcpres$ that is
of maximal (respectively, minimal) horizontal position and, without altering its
angular support, reposition it to be of minimal (respectively, maximal) horizontal
position.  The attached vertical arcs are adjusted in a corresponding manner.
The second flavor considers two distinct horizontal arcs,
$h_i$ and $h_{j}$, 
of $\arcpres$ that are, first, {\em consecutive} and, second, {\em nested}.  Namely, first,
their associated horizontal positions,
$z_i$ \& $z_j$, are consecutive
in the ordering of the horizontal positions.  And second, either
$ [\theta_{i_1^{\ }},\theta_{i_2^{\ }}] \subset [\theta_{j_1^{\ }},\theta_{j_2^{\ }}] $,
or $ [\theta_{j_1^{\ }},\theta_{j_2^{\ }}] \subset [\theta_{i_1^{\ }},\theta_{i_2^{\ }}] $,
or $ [\theta_{i_1^{\ }},\theta_{i_2^{\ }}] \cap [\theta_{j_1^{\ }},\theta_{j_2^{\ }}] = \emptyset $,
i.e. nested.
Then we can locally alter $\arcpres$ 
by replacing the edgepaths $ v_{i-1} \bu h_i \bu v_i $ and $ v_{j - 1} \bu h_{j} \bu v_{j} $
with, respectively, $v_{i-1}^\prime \bu h_i^\prime \bu v_i^\prime$
and $v_{j -1}^\prime \bu h_{j}^\prime \bu v_{j}^\prime$,
where for the corresponding horizontal positions we have
$ z_{i}^\prime = z_{j}$ and $z_{j}^\prime = z_i$, and the vertical support for
$v_{i-1}^\prime$, $v_i^\prime$, $v_{j -1}^\prime$ and $v_{j}^\prime$
are adjusted in a corresponding manner.  The boundary endpoints of these edgepaths
are fixed under this alteration and all other vertical and horizontal arcs of $\arcpres$ are fixed.
(See Figure \ref{figure:elementary moves}(a).)

\noindent \underline{Vertical exchange move}---Let $v_i$ and $v_{j}$ be two distinct vertical
arcs of $\arcpres$ that are, again, {\em consecutive} and {\em nested}.
Namely, the angular positions $\theta_i$ and $\theta_{j}$ are consecutive
in the cyclic ordering of the vertical arcs of $\arcpres$.  And, either
$ [z_{i_1^{\ }},z_{i_2^{\ }}] \subset [z_{j_1^{\ }},z_{j_2^{\ }}] $,
or $ [z_{j_1^{\ }},z_{j_2^{\ }}] \subset [z_{i_1^{\ }},z_{i_2^{\ }}] $,
or $ [z_{i_1^{\ }},z_{i_2^{\ }}] \cap [z_{j_1^{\ }},z_{j_2^{\ }}] 
= \emptyset $, i.e. nested.
Then we can locally alter $\arcpres$ 
by replacing the edgepaths $ h_{i-1} \bu v_i \bu h_i$ \& $h_{j-1} \bu v_j \bu h_j $
with, respectively, $ h_{i-1}^\prime \bu v_i^\prime \bu h_i^\prime$
\& $h_{j-1}^\prime \bu v_j^\prime \bu h_j^\prime $
where for the corresponding angular positions we have
$ \theta_{i}^\prime = \theta_{j}$
and $\theta_{j}^\prime = \theta_i$, and the angular support for
$h_{i-1}^\prime$, $h_i^\prime$, $h_{j -1}^\prime$ and $h_{j}^\prime$
are adjusted in a corresponding manner.
Again, the boundary endpoints of these edgepaths
are fixed under this alteration and all other vertical and horizontal arcs of $\arcpres$ are fixed.
(See Figure \ref{figure:elementary moves}(b).)

\noindent \underline{Horizontal simplification}---Let $h_i$ and $h_{i+1}$ be two horizontal
arcs that are consecutive (as previously defined in the horizontal exchange move), and
are adjacent to a common vertical arc $v_i$ so that
$h_i \bu v_i \bu h_{i+1} \bu v_{i+1}$ is an edgepath on a component of $\arcpres$.
Then we can locally alter $\arcpres$ by replacing $h_i \bu v_i \bu h_{i+1} \bu v_{i+1}$
with an edgepath $h_i^\prime \bu v_{i+1}^\prime$
where: the horizontal position of $h_i^\prime$ is $z_i$;
the angular support of $h_i^\prime$ is $[\theta_{i_1^{\ }},\theta_{i_2^{\ }}] \cup
[\theta_{{(i+1)}_1^{\ }},\theta_{{(i+1)}_2^{\ }}]$; and, the angular position of $v_{i+1}^\prime$ is
$\theta_{i+1}$.  (Notice that $\theta_{i_2^{\ }} = \theta_{{(i+1)}_1^{\ }}$.)
Since our vertical support notation does not have a correspondence to
the orientation of vertical arcs, the vertical support of $v_{i+1}^\prime$
must be specified with some care: it is in fact the closure on $\axis$ of the interval
$\{(z_{i_1^{\ }} , z_{i_2^{\ }}) \cup (z_{{(i+1)}_1^{\ }} , z_{{(i+1)}_2^{\ }})\}
\setminus \{(z_{i_1^{\ }} , z_{i_2^{\ }}) \cap (z_{{(i+1)}_1^{\ }} , z_{{(i+1)}_2^{\ }})\}$.
Again, the boundary endpoints of these edgepaths
are fixed under this alteration and all other vertical and horizontal arcs of $\arcpres$ are fixed.

\noindent \underline{Vertical simplification}---Let $v_i$ and $v_{i+1}$ be two vertical
arcs that are consecutive (as previously defined in the vertical exchange move),
and are adjacent to a common horizontal arc $h_i$ so that
$v_i \bu h_i \bu v_{i+1} \bu h_{i+1}$ is an edgepath on a component of $\arcpres$.
Then we can locally alter $\arcpres$ by replacing $v_i \bu h_i \bu v_{i+1} \bu h_{i+1}$
with an edgepath $v_i^\prime \bu h_{i+1}^\prime$
where: the angular position of $v_i^\prime$ is $\theta_i$;
the vertical support of $v_i^\prime$ is
the closure on $\axis$ of the interval
$\{(z_{i_1^{\ }} , z_{i_2^{\ }}) \cup (z_{{(i+1)}_1^{\ }} , z_{{(i+1)}_2^{\ }})\}
\setminus \{(z_{i_1^{\ }} , z_{i_2^{\ }}) \cap (z_{{(i+1)}_1^{\ }} , z_{{(i+1)}_2^{\ }})\}$;
the horizontal position of $h_{i+1}^\prime$ is
$z_{i+1}$; and the angular support of $h_{i+1}^\prime$
is $[\theta_{i_1^{\ }} , \theta_{i_2^{\ }}] \cup [\theta_{{(i+1)}_1^{\ }} , \theta_{{(i+1)}_2^{\ }}]$.
As before, the boundary endpoints of these edgepaths
are fixed under this alteration and all other vertical and horizontal arcs of $\arcpres$ are fixed.

Given an arc presentation $\arcpres$ we notice that for any sequence of elementary moves
applied to $\arcpres$, the complexity measure $\cC(\arcpres)$ is non-increasing.  That is,
any sequence of elementary moves which includes the uses of either horizontal or vertical
simplification will be {\em monotonically simplified}.

One would hope that for a closed $n$-braid $X$ which admits, respectively, a destabilization,
exchange move, or elementary flype, there exists a sequence of elementary moves
to the arc presentation $\arcpres$ (coming from $ X \ntran \arcpres$) such that for the
resulting arc presentation $X^{\prime\eta}$, the closed $n$-braid $X^\prime$ coming
from $X^{\prime\eta} \btran X^\prime$ admits, respectively, a destabilization, exchange move
or elementary flype (as seen from the standard projection $\pi(X^\prime)$).  Unfortunately,
this is too good to be true.  In order to produce an $X^\prime$ that admits the assumed
isotopy it may be necessary to increase the number of arcs in the arc presentation.
At first glance this seems to disturb our ability to
monotonically simplify.  However, it is possible to control the manner in which
we introduce additional arcs in the arc presentation to monotonically
simplify using an altered complexity measure.  

To accomplish this controlled addition of arcs we introduce the notion of
"shearing intervals".  For a given arc presentation
$\arcpres$ let
$\cI \subset S^1$ be a union of disjoint closed angle intervals of the form
$[\vartheta^k - \e^k , \vartheta^k + \e^k], 0\leq k \leq l$ such that for $ \vartheta \in \cI$ we have that
$H_{\vartheta} \in \fib$ contains no vertical arc of $\arcpres$.  The value of $l$ will be
$1$, $2$, or $3$ when discussing, respectively, destabilization, exchange move, or elementary flype.
Then $\sheararcpres$ is an arc presentation $\arcpres$ along with a specification of where
the {\em shearing intervals} of $\cI$ are to be initially positions in $S^1$.
We require that this initial positioning of $\cI$ be away from the
vertical arcs of $\arcpres$.  As such notice that the positioning of a component of $\cI$ is
characterized by which vertical arcs it lies between.  Thus, once we specify how many
components $\cI$ should contain, up to angular rescaling there are only
finitely many possible initial $\sheararcpres$ for a given $\arcpres$.

We define the complexity measure
$\cC(\sheararcpres)$ to be the number of vertical arcs in the angle interval(s) $S^1 \setminus \cI$.
The reader should notice that for an initial $\sheararcpres$ we have $\cC(\sheararcpres) = \cC(\arcpres)$.

We notice that
our previous elementary moves---horizontal or vertical exchange moves and simplification---can
also be applied to $\sheararcpres$ in the angular intervals of $S^1 \setminus \cI$
so as to monotonically simplify $\sheararcpres$ with respect to its complexity.
Specifically, if any of the edgepaths used in the description of our elementary moves
are totally contained in a component of $S^1 \setminus \cI$ then applying the move
to $\arcpres$ can be seen as applying to move to $\sheararcpres$.
So for horizontal or vertical exchange moves
$\cC(\sheararcpres)$ is unchanged.  Similarly, for the simplification moves
the complexity measure of the resulting
$\sheararcpres$ decreases.



\begin{figure}[htpb]
\centerline{\includegraphics[scale=.6, bb=0 0 707 174]{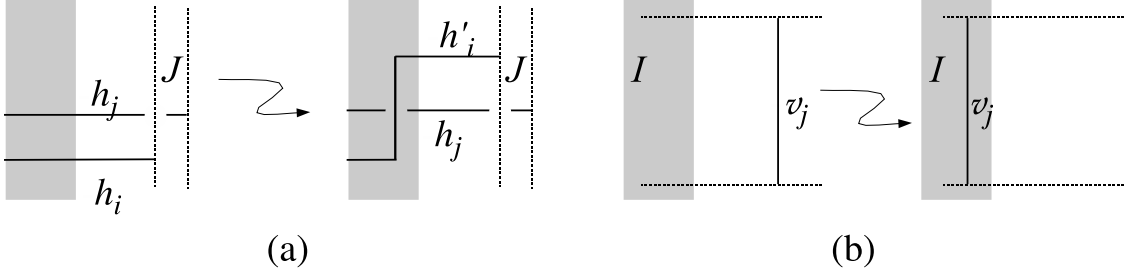}}
\caption{Figure (a) illustrates a shear-horizontal exchange and (b) illustrates
a shear-vertical simplification. Again, in each of the illustrated sequences dotted arcs are
used to indicate that the adjoining vertical arcs (for (a)) and horizontal arcs
(for (b)) have two possible ways of attaching themselves to the labeled solid
arc.}
\label{figure:shear exchange}
\end{figure}

We now add two new elementary moves that utilizes the intervals of $\cI$.

\noindent \underline{Shear horizontal exchange move}---We refer to Figure
\ref{figure:shear exchange}(a).  Let $h_i$ and $h_{j}$ be two
horizontal arcs of $\sheararcpres$ that are {\em consecutive} and {\em nested}
with respect to $\cI$.  That is, first, for some angle interval $J \subset S^1 \setminus \cI$
we have that: a) $h_i$ and $h_j$ intersect $J$; and, b) over 
all of the horizontal arcs of $\sheararcpres$ which intersect $J$,
$h_i$ and $h_{j}$ have consecutive horizontal positions in the ordering
along $\axis$.  Second, the angular support of $h_i \cap J$ is contained
inside the angular support of $h_{j} \cap J$.  Then we can interchange the
horizontal position of these to arcs.  This is achieved by a {\em horizontal shear} inside 
an interval of $\cI$---the
introduction of a vertical (contained in $\cI$) and a horizontal arc as shown in Figure \ref{figure:shear exchange}(a)---
to the portion of $h_i$ and $h_{j}$ that is contained in $\cI$.  Thus, when we
consider the resulting arc presentations an original horizontal exchange move
is realizable.  We abuse notation still referring to the resulted arc presentation as
$\sheararcpres$.  Notice that $\cC(\sheararcpres)$ remains constant.

\noindent \underline{Shear vertical simplification}---We refer to Figure
\ref{figure:shear exchange}(b).  Let $v_j$ be a
vertical arc of $\sheararcpres$ that is {\em consecutive}
with respect to $\cI$.  That is, for an interval $I \in \cI$
there is no vertical arc whose angular position is between $v_j$ and $I$.  We can
then push $v_j$ into $I$.  Again, we abuse notation by referring to the resulted
arc presentation as $\sheararcpres$.  Notice that $\cC(\sheararcpres)$ is decreased
by a count of one.

From now on we refer to horizontal or vertical exchange moves and simplification along
with shear horizontal exchange moves and shear vertical simplification as our
collection of elementary moves on arc presentations with shearing intervals.
Since our notation, $\sheararcpres$, is for an arc presentation along
with unchanging shearing intervals under elementary moves, it is
very convenient to abuse notation and use $\sheararcpres$ when making
set operation statements involving the underlying arc presentation.  For example,
by $\sheararcpres \cap C_1$ we will mean $\arcpres \cap C_1$.

To connect the dots, 
for a given $\arcpres$ and fixing the number of angle intervals in $\cI$,
up to re-scaling, there are only finitely many initial $\sheararcpres$,
i.e. only finitely many {\em combinatorial distinct} initial $\sheararcpres$.
The complexity measure $\cC(\sheararcpres)$ is just a count on the number of
vertical arcs in $S^1 \setminus \cI$ and that elementary moves on a $\sheararcpres$
never increases the number of vertical arcs in $S^1 \setminus \cI$.
So starting with $\cC(\arcpres)$ vertical arcs in the intervals of
$S^1 \setminus \cI$, after any sequence of elementary moves there is
a finite number of possible combinatorial distinct rectangular (non-closed) braid presentations in the intervals of $S^1 \setminus \cI$.  (Again, ``combinatorial distinct'' refers to
equivalency up to rescaling of angle and height positions of arcs.)
Thus, starting with an initial $\sheararcpres$ (having no vertical arcs in $\cI$)
there are only finitely many possible rectangular diagrams occurring in the intervals
$S^1 \setminus \cI$ after any sequence of elementary moves to $\sheararcpres$
up to our combinatorial equivalence.

Due to the finite number of rectangular diagrams occurring in the
$S^1 \setminus \cI$ intervals, for any sequence of elementary moves
on $\sheararcpres$ that produces only
combinatorial distinct diagrams in $S^1 \setminus \cI$,
a bounded number of additional horizontal and vertical arcs
inside the intervals of $\cI$ through the application of shear horizontal exchange moves
and shear vertical simplifications will be introduced.
However, it is possible to have arbitrarily long sequences of elementary moves on $\sheararcpres$
containing the occurrence of same combinatorial
distinct rectangular diagrams in the intervals $S^1 \setminus \cI$ arbitrarily many times.
For example, starting with
a fixed $\sheararcpres$ one could produce a finite cyclic sequence of
rectangular diagrams---starting and ending at the same fixed $\sheararcpres$---all
having the same complexity measure by applying a
sequence of horizontal, vertical and shear horizontal exchange moves.  If we
repeat such a sequence any number of times we
can create the canceling $\a$-braiding phenomena of Remarks
\ref{remark:elementary exchange moves and flypes} \&
\ref{remark: equivalency class for EM and F}, or the canceling $\tau$-braiding phenomena
of Remark \ref{remark: double destabilization}.
We concluded that although what can occur in the $\cI$ intervals may be infinite, what can occur
in $S^1 \setminus \cI$ is finite, and recognizing when a closed braid
admits one of our isotopies will dependent on interrupting the diagrams
in $S^1 \setminus \cI$.

We are now in a position to state our main results. 

\begin{thm}
\label{theorem:main result}
Let $X$ be a closed $n$-braid such that $\cB_n(X)$ admits, respectively, a destabilization,
exchange move or elementary flype.  Consider any arc presentation coming from the
a presentation transition $X \ntran \arcpres$.  Then there exists
a set of intervals $\cI$ and a sequence of arc presentations
$$ \sheararcpres = X^0 \ra X^1 \ra \cdots \ra X^l = X^{\prime\eta}_{\cI} $$
such that:
\bi
\item[1.] If $\cB_n(X)$ admits, respectively, a destabilization, exchange move, flype
or elementary flype then $\cI$ has, respectively, one, two or three intervals. 
\item[2.] $X^{i+1}$ is obtained from $X^i$ via one of the elementary moves.  All of these moves
are with respect to the intervals of $\cI$.
\item[3.] The closed $n$-braid obtained from the presentation transition
$X^{\prime\eta}_{\cI} \btran X^\prime$ admits, respectively, a destabilization, exchange move,
flype or elementary flype (as seen from the standard projection $\pi(X^\prime)$).
\ei
\end{thm}

We observe that $\cC(X^{i+1}) \leq \cC(X^i)$ for $0 \leq i \leq l$ for
all applications of our elementary moves.  In particular,
if $X^i \ra X^{i+1}$ corresponds to a horizontal, vertical or shear horizontal exchange move
then $\cC(X^i) = \cC(X^{i+1})$.  If it corresponds to a horizontal, vertical or shear
vertical simplification then $\cC(X^{i+1}) < \cC(X^i)$, i.e. monotonically simplified.

As previously remarked, we will restrict ourselves to only
three possible choices for $\cI$---it has one, two or three intervals---and we recall the positioning
of these intervals is characterized by which vertical arcs of a initially given $\arcpres$ they lie between.
Thus, there are only finitely many possible initial $\sheararcpres$.  Also, our previous remarks
gives us that there are only finitely many resulting $\sheararcpres$ after elementary moves.
The production of such a finite set is easily seen as algorithmic.
Therefore, Theorem \ref{theorem:main result}
implies the following corollary.

\begin{cor}
\label{corrollary:main result}
There exists an algorithm for deciding whether a closed $n$-braid is braid
isotopic to one that admits either a destabilization, exchange move, or elementary flype.
\end{cor}

The construction of our algorithmic solutions
comes from utilizing the braid foliation machinery that was first
developed in \cite{[BF],[BM1],[BM2],[BM3]} and further
refined in the beautiful work of I.A. Dynnikov \cite{[D]}.

In \cite{[M1]}, exploiting the interpretation of $n$-braids as elements
of the mapping class group of the $n$-punctured disc and Nielsen-Thurston's theory,
an alternate algorithm is established for determining when a closed braid admits
a destabilization.  This algorithm is based upon an analysis of
the action of the mapping class group on
the geodesics of the $n$-punctured disc endowed with a fixed hyperbolic metric.

In a strict sense, Theorem \ref{theorem:main result} is an existence result---it tells
us if a given braid admits a given move.  Based upon Remarks
\ref{remark:elementary exchange moves and flypes} \&
\ref{remark: double destabilization} we know
that determining whether two fixed braids are related by a particular move is problematic.
But, by paying close attention to the machinery in the proof of Theorem \ref{theorem:main result}
we can make such a determination in the cases of a double destabilization,
thin exchange move, and elementary flype.  Thus, we have the following theorem.

\begin{thm}
\label{Theorem:destably equivalent}
There exists an algorithm for deciding whether closed braid $Y$ is related to $X$ by
a thin exchange move, an elementary flype.  And, there exists an algorithm for
deciding whether $X$ admits a double destabilization.
\end{thm}

\noindent
{\footnotesize ACKNOWLEDGMENTS---The author wishes to thank the
referee's considerable input and energy.
Specifically, the author was alerted to the connections to the work
of A.V. Malyutin by the referee; and, Theorem \ref{Theorem:destably equivalent} is
essentially the referee's formulation.
The author also acknowledges NSF partial support through grant \#DMS 0306062.}

\section{The cylinder machinery.}
\label{section:cylinder machinery}
\subsection{Destabilizing, exchange and flyping discs.}
\label{subsection:d,e,f discs}
Our first objective is the give a geometric characterization for recognizing
when a closed $n$-braid is braid isotopic to one that admits either a destabilization,
exchange move or elementary flype.  All geometric characterizations will depend on
the existence of a specified embedded disc.  Our characterizations will, in fact, occur
in pairs: one for the braid presentation and one for the arc presentation.

All of our geometric characterizing discs, $\Delta_\v$ will be {\em above the braid}.  That is,
$\Delta_\v = D_{+1} \cup N$ where:
\bi
\item[1.] $D_{+1}$ is the disc in a plane $z = z_{max}$ having
$0 \leq r \leq 1 , \theta \in S^1$ and $z_{max}$ being a constant greater the
horizontal positions of all the horizontal arcs of $\arcpres$. (Thus,
$D_{+1} \cap \arcpres = \emptyset$.)
\item[2.] $\Delta_\e$ is oriented so that $D_{+1}$ necessarily intersects $\axis$
geometrically and algebraically $+1$ at a {\em vertex} point $\ver_{max}$.
\item[3.] $N$ is an annulus having $r \geq 1$ for all its points.
\ei

\noindent {\bf Destabilizing disc}---(\underline{Braid presentation})
Let $X$ be a closed $n$-braid which admits
a destabilization, i.e. the corresponding braid word $\b (X) = W \sigma_{n-1}^{\pm 1}$
with $W \in \cW^{n-2}$.
Then there exists a {\em destabilizing disc} $\Delta_d$ having the following properties.
\bi
\item[D-a.] $\partial \Delta_d = \a_{h} \bu \a_{v}^\partial \bu$ where we have the
{\em horizontal boundary} $\a_{h} \subset X$ and the {\em $\partial$-vertical arc}
$\a_{v}^\partial \subset H_{\theta_{\ }^\partial}$ for some $H_{\theta_{\ }^\partial} \in \fib$
\item[D-b.] $\Delta_d \cap X = \a_{h}$.
\item[D-c.] $\Delta_d$ transversely intersects $\axis$ at a single vertex point $\ver_{max}$.
\item[D-d.] $a_\theta = \Delta_d \cap H_\theta$ is a single arc having an endpoint on $\a_{h}$ and $\ver_{max}$ as the other endpoint, for $\theta \in S^1 - \{\theta_{\ }^\partial\}$.  If $\theta = \theta_{\ }^\partial$
then $\a_{v}^\partial \subset a_{\theta_{\ }^\partial}^{\ } (= \Delta_d \cap H_{\theta^\partial})$.  We require all
leaves to homeomorphic to the unit interval of $\mathbb R$.
To summarize, the braid fibration induces a {\em radial foliation} on $\Delta_d$.
\ei
If $X^\prime$ is braid isotopic to $X$ then we can extend the braid isotopy which takes
$X$ to $X^\prime$ to an ambient isotopy of $S^3 \setminus \axis$.  This ambient isotopy takes
$\Delta_d$ to a destabilizing disc for $X^\prime$, i.e. properties D-a through D-d are
still satisfied.  Thus, every $n$-braid representative of $\cB_n (X)$ will have a
destabilizing disc.

Further analysis of the leaves containing the $\partial$-vertical arcs is useful.
Specifically, let $a_{{\theta_{\ }^\partial}}^{\ }$ be a leaf in the
radial foliation of $\Delta_d$ that
contains an $\partial$-vertical arc $\a_{v}^\partial$.  If there is an angle interval
$[\theta^\prime , \theta^\partial] (\subset S^1)$ such that
for any $\theta \in [\theta^\prime, \theta_{\ }^\partial]$ pushing the leaf $a_\theta$ forward in
the radial foliation to $a_{\theta_{\ }^\partial}^{\ }$
corresponds to a homeomorphism between $a_\theta$ and
$a_{\theta_{\ }^\partial}^{\ }$ then we say that $\a_{v}^\partial$ has a {\em front edge}.
Similarly, if there is an angle interval
$[\theta_{\ }^\partial, \theta^\prime] (\subset S^1)$ such that
for any $\theta \in [\theta_{\ }^\partial, \theta^\prime]$ pushing the leaf $a_\theta$ backwards in
the radial foliation
to $a_{\theta_{\ }^\partial}^{\ }$
corresponds to a homeomorphism between $a_\theta$ and
$a_{\theta_{\ }^\partial}^{\ }$ then we say that $\a_{v}^\partial$ has a {\em back edge}.
Notice that by the definition of $\Delta_d$, $\a_{v}^\partial$ has either a front edge or a back
edge, but not both.

Without loss of generality we make the convenient assumption that if
$\Delta_d$ is a positive (respectively, negative) destabilizing disc then
the $\partial$-vertical arc is a front (respectively, back) edge.  Notice that
an edge assignment is stable under braid isotopy.

\noindent
(\underline{Arc presentation})  Again with $\b(X) = W \sigma_{n-1}^{\pm 1}$,
we consider the arc presentation coming from
the transition $X \ntran \arcpres$.
Then we also have a disc, which is call an {\em obvious destabilizing
disc}, that has the following properties.
\bi
\item[DA-a.] $\partial \Delta_d = h_1 \bu v_1 \bu h_2 \bu v_2^\partial \bu$ where
$v_1$ is a vertical arc of $\arcpres$, $h_1$ \& $h_2$ are contained in horizontal
arcs of $\arcpres$, and $v_2^\partial \subset H_{\theta_{\ }^\partial}$ for $H_{\theta_{\ }^\partial} \in \fib$
is the {\em $\partial$-vertical arc} of $\Delta_d$.
\item[DA-b.] $\Delta_d \cap \arcpres = h_1 \bu v_1 \bu h_2$.
\item[DA-c.] $\Delta_d$ transversely intersects $\axis$ at a single vertex point $\ver_{max}$, the
{\em vertex} of the foliation on $\Delta_d$.
\item[DA-d.] $a_\theta = \Delta_d \cap H_\theta$ is a single arc having
an endpoint $\ver_{max}$.  Moreover,
when $H_\theta$ does not contain $v_1$ or $v_2^\partial$, $a_\theta$ has an endpoint on either
$h_1$ or $h_2$.  When $H_\theta$ does contain $v_1$ (respectively $v_2^\partial$),
$v_1 \subset a_\theta$ (respectively $v_2^\partial \subset a_\theta$).
In particular, the braid fibration induces a {\em radial foliation} on $\Delta_d$.
\ei

\begin{figure}[htpb]
\centerline{\includegraphics[scale=1, bb=0 0 210 95]{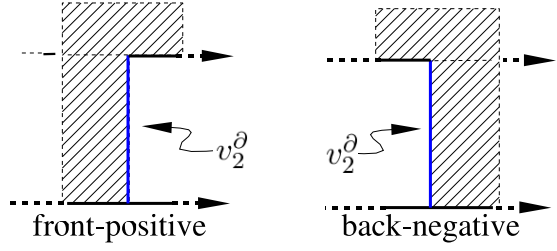}}
\caption{The illustration depicts the braid above the disc.
The viewpoint is of an observer at a point with $r > 1$.
The left configuration has the horizontal positions of the horizontal arcs
increasing, whereas the right configuration has the horizontal
positions of the horizontal arcs decreasing.}
\label{figure:boundary of disc}
\end{figure}

Recall that $\Delta_d = D_{+1} \cup N$ where the points of $N$ have coordinates
with $r \geq 1$.  Near $\partial \Delta_d$ we can use $N$ to determine the
parity of the associated destabilization.  The possibilities are easily
listed.  As we traverse in the positive direction an edgepath neighborhood in
$\partial \Delta_d$ of $v_2^\partial$, the horizontal position of the horizontal
arcs adjacent to $v_2^\partial$
is either increasing or decreasing.  Then it is easily checked that our parity scheme is as follows.
In $N$ near $v_2^\partial$ if we have $\partial \Delta_d$ with increasing horizontal position then the parity of the destabilization is positive.
For decreasing horizontal position
the associated destabilization is negative.
Refer to Figure \ref{figure:boundary of disc}.
We remind the reader that our choice of edge assignment has the $\partial$-vertical arc
as a front (respectively, back) edge
for a positive (respectively, negative) destabilizing disc.

We will also see that this parity scheme and edge characterization can be used for
the exchange move and flyping discs.

\noindent {\bf Exchange move disc}---(\underline{Braid presentation})
Let $X$ be a closed $n$-braid which admits
an exchange move, i.e. the corresponding braid word $\b (X) = W U$
where $W \in \cW^t$, $U \in \cU^s$ and $1< s \leq t < n-1$.
Then there exists a {\em exchange disc} $\Delta_e$ above the braid having the following properties.
\bi
\item[E-a.] $\partial \Delta_e = \a_{h_1^{\ }}^{\ } \bu \a_{v_1^{\ }}^\partial \bu \a_{h_2^{\ }}^{\ } \bu \a_{v_2^{\ }}^\partial \bu$
where we have the {\em horizontal boundary} $\a_{h_1^{\ }}^{\ }, \a_{h_2^{\ }}^{\ } \subset X$; and
the {\em $\partial$-vertical arcs}
$\a_{v_i^{\ }}^\partial \subset H_{\theta_i^{^\partial}}^{\ }$ where $H_{\theta_i^{^\partial}}^{\ } \in \fib$, $i = 1,2$.  We require that there be both a front and a back edge $\partial$-vertical arc.
\item[E-b.] $\Delta_e \cap X = \a_{h_1^{\ }}^{\ } \cup \a_{h_2^{\ }}^{\ }$.
\item[E-c.] $\Delta_e$ transversely intersects $\axis$ at a single vertex point $\ver_{max}$.
\item[E-d.] $\{\a_{h_1^{\ }}^{\ } \cup \a_{h_2^{\ }}^{\ } \} \cap H_\theta \not= \emptyset$ for $H_\theta \in \fib$
and $\theta \not\in \{ \theta^\partial_1 , \theta^\partial_2 \}$.
\item[E-e.] $a_\theta = \Delta_e \cap H_\theta$ is a single arc
having: an endpoint on $\ver_{max}$;
an endpoint on $\a_{h_1^{\ }}^{\ } \cup \a_{h_2^{\ }}^{\ }$ when $\theta \not\in \{ \theta^\partial_1 , \theta^\partial_2 \}$;
and contains $\a_{v_i^{\ }}^\partial$ when $\theta = \theta^\partial_i$, $i=1,2$.
In particular, the braid fibration induces a {\em radial foliation} on $\Delta_e$.
\ei
If $X^\prime$ is braid isotopic to $X$ then we can extend the braid isotopy which takes
$X$ to $X^\prime$ to an ambient isotopy of $S^3 \setminus \axis$.  This ambient isotopy takes
$\Delta_e$ to a exchange disc for $X^\prime$, i.e. properties E-a through E-e are
still satisfied.  Thus, every $n$-braid representative of $\cB_n (X)$ will have an
exchange disc.  Again, the edge assignment of $\partial$-vertical arcs is stable
under braid isotopy.

\noindent
(\underline{Arc presentation})  Again with $\b(X) = WU$,
we consider the arc presentation coming from
the transition $X \ntran \arcpres$.
Then we also have a disc, which is call an {\em obvious exchange
disc} above the braid, that has the following properties.
\bi
\item[EA-a.] $\partial \Delta_e = h_1 \bu v_1^\partial \bu h_2 \bu v_2^\partial \bu$ where
$h_1$ and $h_2$ are subarcs of two different horizontal arcs
of $\arcpres$; and $v_i^\partial \subset H_{\theta_i^{^\partial}}$ for $H_{\theta_i^{^\partial}} \in \fib$,
$i=1,2$, are the $\partial$-vertical arcs of $\partial \Delta_e$.
We require that there be both a front and a back edge $\partial$-vertical arc.
\item[EA-b.] $\Delta_e \cap \arcpres = h_1 \cup h_2$.
\item[EA-c.] $\Delta_e$ transversely intersects $\axis$ at a single vertex point $\ver_{max}$.
\item[EA-d.] $\{ h_1 \cup h_2 \} \cap H_\theta \not= \emptyset$ for $H_\theta \in \fib$
and $\theta \not\in \{ \theta^\partial_1 , \theta^\partial_2 \}$.
\item[EA-e.] $a_\theta = \Delta_e \cap H_\theta$ is a single arc having an endpoint $\ver_{max}$.  Moreover,
when $H_\theta$ does not contain $v_i^\partial$, $a_\theta$ has an endpoint on either
$h_1$ or $h_2$.  When $H_\theta$ does contain $v_1^\partial$ (respectively $v_2^\partial$),
$v_1^\partial \subset a_\theta$ (respectively $v_2^\partial \subset a_\theta$).
In particular, the braid fibration induces a {\em radial foliation} on $\Delta_d$.
\ei

The reader should observe that when we consider $N \subset \Delta_e$ of
$\partial \Delta_e$ near $v_1^\partial$ and $v_2^\partial$,
we must necessarily have one horizontal boundary arc
being positive and one being negative in the sense of Figure \ref{figure:boundary of disc}.

\noindent {\bf Flyping disc}---(\underline{Braid presentation})
Let $X$ be a closed $n$-braid which admits
an elementary flype, i.e.
$\b (X) = W_1 \sigma_{n-1}^p W_2 \sigma_{n-1}^{\pm 1}$.  The embedded disc
we will use to illustrate the presence of a flype can conceptually be seen
as an amalgamation of a destabilizing disc and an exchange disc, since the flype
involves a 'flyping crossing' and 'flyping block'.
Specifically,  there exists a {\em flyping disc} $\Delta_f$ having the following properties.
\bi
\item[F-a.] $\partial \Delta_f = \a_{h_1^{\ }} \bu \a_{v_1^{\ }}^\partial \bu \a_{h_2^{\ }} \bu \a_{v_2^{\ }}^\partial
\bu \a_{h_3^{\ }} \bu \a_{v_3^{\ }}^\partial \bu$
where we have the {\em horizontal boundary arcs} $\a_{h_i^{\ }} \subset X$; and
the {\em $\partial$-vertical arcs}
$\a_{v_i^{\ }}^\partial \subset H_{\theta_i^{^\partial}}$ where $H_{\theta_i^{^\partial}} \in \fib$, $i = 1,2,3$.  When $\Delta_f$ corresponds to a positive (respectively, negative) flype
we have edge assignments as follows: $\a_{v_1^{\ }}^\partial$ front,
$\a_{v_2^{\ }}^\partial$ back, $\a_{v_3^{\ }}^\partial$ front (respectively,
$\a_{v_1^{\ }}^\partial$ back,
$\a_{v_2^{\ }}^\partial$ front, $\a_{v_3^{\ }}^\partial$ back).
\item[F-b.] $\Delta_f \cap X = \a_{h_1^{\ }} \cup \a_{h_2^{\ }} \cup \a_{h_3^{\ }}$.
\item[F-c.] $\Delta_f$ transversely intersects $\axis$ at a single vertex point $\ver_{max}$.
\item[F-d.] $\{\a_{h_1^{\ }} \cup \a_{h_2^{\ }} \cup \a_{h_3^{\ }}\} \cap H_\theta \not= \emptyset$
for $H_\theta \in \fib$ and $\theta \not\in \{ \theta^\partial_1 , \theta^\partial_2 , \theta^\partial_3 \}$.
\item[F-e.] $a_\theta = \Delta_f \cap H_\theta$ is a single arc
having: an endpoint on $\ver_{max}$; an
endpoint on $\a_{h_1^{\ }} \cup \a_{h_2^{\ }} \cup \a_{h_3^{\ }}$ when
$\theta \not\in \{ \theta^\partial_1 , \theta^\partial_2 , \theta^\partial_3 \}$;
and contains $\a_{v_i^{\ }}^\partial$ when angle is one of the angles $\theta^\partial_i$, $i=1,2,3$.
In particular, the braid fibration induces a {\em radial foliation} on $\Delta_f$.
\ei
If $X^\prime$ is braid isotopic to $X$ then we can extend the braid isotopy which takes
$X$ to $X^\prime$ to an ambient isotopy of $S^3 \setminus \axis$.  This ambient isotopy takes
$\Delta_f$ to a flyping disc for $X^\prime$, i.e. properties F-a through F-e are
still satisfied.  Thus, every $n$-braid representative of $\cB_n (X)$ will have an
flyping disc.

\begin{figure}[htpb]
\centerline{\includegraphics[scale=.9, bb=0 0 458 137]{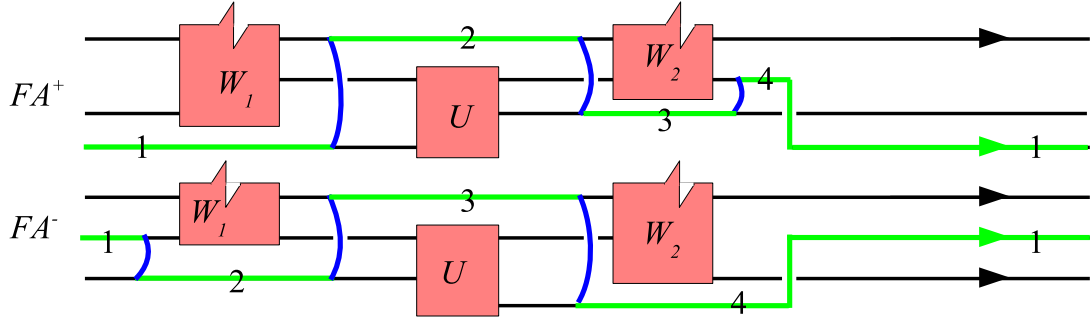}}
\caption{The figure illustrates the boundaries of the obvious flyping discs
when the vertical arcs and boundary arcs are above.  $FA^+$ is the positive flype
and $FA^-$ is the negative flype.  To reduce clutter
we simply label the horizontal arcs $h_i$ (in green) with their subscripts.
The three $\partial$-vertical arcs $v^\partial_i, \ i=,1,2,3,$ are in blue and
the single vertical arc, $v_4$ is in green.}
\label{figure:types of flyping disc}
\end{figure}

\noindent
(\underline{Arc presentation})  Again with $\b(X) = W_1 U W_2 \sigma_{n-1}^{\pm 1}$,
we consider the arc presentation coming from
the transition $X \ntran \arcpres$.
Then we also have a disc, which is call an {\em obvious flyping
disc}.
\bi
\item[FA-a.] $\partial \Delta_f = h_1 \bu v_1^\partial \bu h_2 \bu v_2^\partial
\bu h_3 \bu v_3^\partial \bu h_4 \bu v_4 \bu$ where
$h_1, h_2, h_3, h_4$ are subarcs of differing horizontal arcs of
$\arcpres$; $v_i^\partial \subset H_{\theta^\partial_i}$ for $H_{\theta^\partial_i} \in \fib$,
$i=1,2,3$ are {\em $\partial$-vertical arcs} with $r \geq 1$; and $v_4$ is a vertical arc of $\arcpres$.
When $\Delta_f$ corresponds to a positive (respectively, negative) flype
we have edge assignments as follows: $v_1^\partial$ front,
$v_2^\partial$ back, $v_3^\partial$ front (respectively,
$v_1^\partial$ back,
$v_2^\partial$ front, $v_3^\partial$ back).
\item[FA-b.] $\Delta_f \cap X = h_2 \cup h_3 \cup (h_4 \bu v_4 \bu h_1)$.
\item[FA-c.] $v_1^\partial , v_2^\partial \subset \Delta_f$
are above the braid with opposite parity in the sense of Figure \ref{figure:boundary of disc}.
\item[FA-d.] $v_3^\partial$ is above and positive (respectively, negative)
in which case $X$ admits a positive (respectively, negative) flype.
\item[FA-e.] $\Delta_f$ transversely intersects $\axis$ at a single vertex point $\ver_{max}$.
\item[FA-f.] The braid fibration induces a radial foliation on $\Delta_f$.
That is, for all $H_\theta \in \fib$, $\Delta_f \cap H_\theta$ is a single arc having
$\ver_{max}$ as one endpoint; when $H_\theta$ does not contain a vertical arc of $\arcpres$
or a $\partial$-vertical arc of $\partial \Delta_f$ then the intersection arc also has an
endpoint on a horizontal portion of $\partial \Delta_f$; and, when $H_\theta$ does contain
a vertical arc or $\partial$-vertical arc then that arc is in the single arc of
$\Delta_f \cap H_\theta$.
\ei

Returning our focus on establishing Theorem \ref{theorem:main result}, we now
have the following proposition.

\begin{prop}
\label{proposition:obvious disc}
Let $X$ be a closed $n$-braid such that $\cB_n(X)$ admits, respectively, a destabilization,
exchange move or elementary flype.  Consider any arc presentation coming from
a presentation transition $X \ntran \arcpres$.  Then there exists
a set of intervals $\cI$ and a sequence of arc presentations
$$ \sheararcpres = X^0 \ra X^1 \ra \cdots \ra X^l = X^{\prime\eta}_{\cI} $$
such that:
\bi
\item[1.] If $\cB_n(X)$ admits, respectively, a destabilization, exchange move
or elementary flype then $\cI$ has, respectively, one, two or three intervals. 
\item[2.] $X^{i+1}$ is obtained from $X^i$ via one of the elementary moves---horizontal
exchange move, vertical exchange move, horizontal simplification, vertical simplification,
shear horizontal exchange move, and shear vertical simplification.
\item[3.] $\cC(X^{i+1}) \leq \cC(X^i)$ for $0 \leq i \leq l$.  In particular,
if $X^i \ra X^{i+1}$ corresponds to a horizontal, vertical or shear horizontal exchange move
then $\cC(X^i) = \cC(X^{i+1})$.  If it corresponds to a horizontal, vertical or shear
vertical simplification then $\cC(X^{i+1}) < \cC(X^i)$. Thus, our sequence will
monotonically simplify.
\item[4.] If $\cB_n(X)$ admits, respectively, a destabilization, exchange move
or elementary flype then here exists an obvious, respectively,
destabilizing, exchange, or flyping disc $\Delta_\v , \ \v \in \{d,e,f\}$ for
$X^{\prime\eta}_{\cI}$ such that the angular positions at which the
$\partial$-vertical arcs of $\Delta_\v$ occur are in the set of $H_\theta$
containing $\partial \cI$.
\ei
\end{prop}

The proof of Proposition \ref{proposition:obvious disc} requires understanding ``notch discs''.

\subsection{Notch discs.}
\label{subsection:notch discs}

Let $X$ be a closed $n$-braid presentation that
is braid isotopic to a braid that admits either a destabilization, exchange move or
elementary flype.  Let $\Delta_\v$, $\v \in \{d,e,f\}$, be an appropriate disc illustrating
the isotopy.
We take a transition from a braid to an arc presentation, $X \ntran \arcpres$,
and consider the impact of this transition on
$\Delta_\v$.
As before, we denote the leaves of the radial foliation of $\Delta_\v$ by
$a_\theta (= \Delta_\v \cap H_\theta)$; and $a_\theta$ contains a $\partial$-vertical
arc when $\theta$ is $\theta^\partial_i$ for an appropriate $i = \{1,2,3\}$.  
Since $X$ can be positioned to be arbitrarily close to $\arcpres$
(where closeness is measured by
the standard metric for the $(r,\theta,z)$ coordinates), by an ambient isotopy of
$S^3$ that preserves the boundary and foliation properties of $\Delta_\v$, we can
assume that $\arcpres$ intersects $\Delta_\v$ such that we have the following:
\bi
\item[i.] For each leaf $a_\theta$, where $\theta \in S^1$ is not an angular position
of a $\partial$-vertical arc,
we have that $a_\theta$ intersects either the interior of a single horizontal arc or
the interior of a single vertical arc of $\arcpres$.
\item[ii.] For $a_\theta$ where $\theta \in S^1$ is not an angular position
of a $\partial$-vertical arc,
if $a_\theta$ intersects the interior of a horizontal arc of $h \subset \arcpres$ then
$a_\theta \cap h$ is a single point.
\item[iii.] For $a_\theta$ where $\theta \in S^1$ is not an angular position
of a $\partial$-vertical arc,
if $a_\theta$ intersects the interior of a vertical arc of $v \subset \arcpres$ then
$v \subset a_\theta$.
\item[iv.] For $a_\theta$ containing $\partial$-vertical arc $v^\partial$,
there are two horizontal arcs $h^\prime ,h^{\prime\prime} \subset \arcpres$
such that
\bi
\item[a.] $ a_{\theta} \cap \arcpres = a_{\theta} \cap \{ h^\prime \cup h^{\prime\prime} \}$
\item[b.] $ v^\partial $ has its endpoints on $h^\prime$ and $h^{\prime\prime}$.
\ei
\ei

\begin{figure}[htpb]
\centerline{\includegraphics[scale=.65, bb=0 0 771 378]{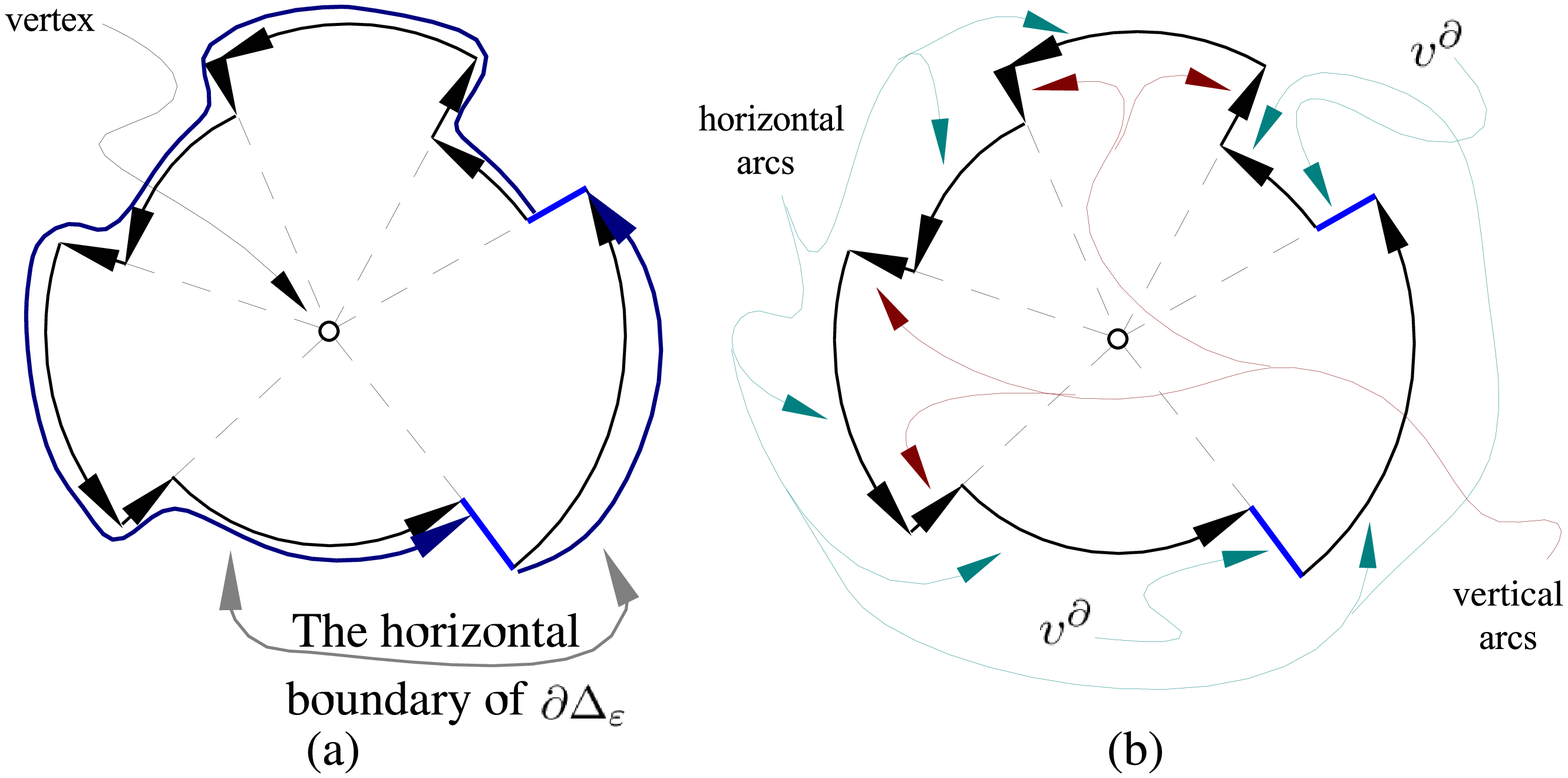}}
\caption{In illustrate (a) we have drawn $\notchdisc \subset \Delta_\v$.  In
(b) we indicate what the structure of $\partial\notchdisc$.}
\label{figure:notch}
\end{figure}

With these conditions holding we can {\em notch} $\Delta_\v$ to
produce $\notchdisc$.  That is, $\notchdisc \subset \Delta_\v$
is the sub-disc whose boundary is obtain by projecting $X \cap \Delta_\v$ along
the leaves $a_\theta$ to $\arcpres \cap \Delta_\v$.  (See Figure \ref{figure:notch}.)
The boundary of $\notchdisc$ is then a union of three types of arcs: {\em horizontal arcs},
$h^\eta_j$, that can be either arcs or sub-arcs of the horizontal arcs of $\arcpres$;
{\em vertical arcs}, $v^\eta_j$,
that are in fact arcs coming from the vertical arcs of $\arcpres$;
and {\em $\partial$-vertical arcs}, $v^\partial_i$.
We then have $\partial \notchdisc$ being a cyclic ordered union of arcs
alternating between horizontal and vertical arc types, i.e.
$$ \partial \notchdisc = h^\eta_1 \bu v^\eta_1 \bu \cdots \bu h^\eta_{j_{\ }^1} \bu
v^\partial_1 \bu h^\eta_{j_{\ }^1 +1} \bu \cdots \bu h^\eta_{j_{\ }^p} \bu
v^\partial_p \bu, $$
where $p$ is $1$ for $\v = d$, $2$ for $\v = e$ and $3$ for $\v = f$.

Extending our notation for the transition between braid presentations and arc presentations,
we will use $(X , \Delta_\v) \ntran (\arcpres , \notchdisc)$ and
$(\arcpres , \notchdisc) \btran (X, \Delta_\v ) $ for indicating the transition
between presentation-disc pairs.  Clearly, given a pair $(\arcpres , \notchdisc )$
when elementary moves are applied to $\arcpres$ there can be an alteration to
the positioning of $\notchdisc$ and/or the arc decomposition of $\partial \notchdisc$.

\begin{lemma}[First simplification of $(\arcpres , \notchdisc)$]
\label{lemma:first simplification of notch disc}
Let $X$ be a closed $n$-braid such that $\cB_n(X)$ admits, respectively, a destabilization,
exchange move or elementary flype and let $\Delta_\v$ be, respectively, a
destabilizing disc, exchange move disc, or elementary flyping disc, i.e.
$\v \in \{ d , e , f \}$.  Then there exists an alternate disc $\Delta_\v^\prime$
such that for $(X , \Delta_\v^\prime) \ntran (\arcpres , \notchdiscprime)$ we have
\bi
\item[a.] $\arcpres$ is unchanged.
\item[b.] $ \partial \notchdiscprime = h^\eta_1 \bu \a^\eta_{\theta_1^{\ }} \bu
h^\eta_{2} \bu v^\eta_2 \bu \cdots \bu h^\eta_{l} \bu v^\eta_l \bu $ when $\v = d$.
\item[c.] $ \partial \notchdiscprime = h^\eta_1 \bu \a^\eta_{\theta_1^{\ }} \bu h^\eta_{2}
\bu v^\eta_2 \bu \cdots \bu  h^\eta_{l} \bu \a^\eta_{\theta_2^{\ }} \bu$ when $\v = e$.
\item[d.] $ \partial \notchdiscprime = h^\eta_1 \bu \a^\eta_{\theta_1^{\ }} \bu h^\eta_{2}
\bu \a^\eta_{\theta_2^{\ }} \bu h^\eta_{3} \bu v^\eta_3 \bu \cdots \bu
h^\eta_{4} \bu v^\eta_4 \bu \cdots \bu h^\eta_{l} \bu \a^\eta_{\theta_3^{\ }} \bu$ when
$\v = f$.
\ei
\end{lemma}

\pf
The statements b. through d. are achieved by performing an isotopy of the arcs
${\a^\eta_{\theta_i^{\ }}} {\rm 's}$.  In particular, statement b. is true by construction.
To achieve statement c. while maintaining the truth of statement a. we start with
the $\partial$-vertical arcs of $\notchdisc$, $\a^\eta_{\theta_1^{\ }}, \a^\eta_{\theta_2^{\ }}$.
If these two arcs have endpoints on a common horizontal arc of $\arcpres$ then we are done.
If not then we push $\a^\eta_{\theta_2^{\ }}$ backward (or forward) through the disc fibers of
$\fib$.  This push will naturally isotop $\a^\eta_{\theta_2^{\ }}$ in the disc fibers.  We stop
our push when $\a^\eta_{\theta_2^{\ }}$ has an endpoint
on a horizontal arc that $\a^\eta_{\theta_1^{\ }}$
also has an endpoint on.  This corresponds to an ambient isotopy of the graph
$\arcpres \cup \a^\eta_{\theta_1^{\ }} \cup \a^\eta_{\theta_2^{\ }}$
in $\reals^3 \setminus \axis$.
There is still a disc whose boundary is the union of two subarc and the resulting
two $\partial$-vertical arcs.  This new disc is our $\notchdiscprime$.  It is easy to see
that $\notchdiscprime$ is in fact a notch disc.

Similarly, for achieving statement d. while maintaining the truth of statement a.
we first push $\a^\eta_{\theta_2^{\ }}$ backward until it has an endpoint on a horizontal
arc that also contains an endpoint of $\a^\eta_{\theta_1^{\ }}$.  Not we push
$\a^\eta_{\theta_3^{\ }}$ backward until it has an endpoint on a horizontal arc that
contains an endpoint of $\a^\eta_{\theta_2^{\ }}$.  Since both pushes are ambient
isotopies of $\reals^3 \setminus \axis$ we again have an new notch disc.
\qed

We will refer to the portion
$h^\eta_{2} \cup v^\eta_2 \cup \cdots \cup h^\eta_{l} \cup v^\eta_l
\subset \partial \Delta^{\prime \eta}_d$ as the {\em middle boundary}
of $\Delta^{\prime \eta}_d$.
The {\em middle boundary} of $\partial \Delta^{\prime \eta}_e$
(respectively $\partial \Delta^{\prime \eta}_f$) is
$ h^\eta_{2} \bu v^\eta_2 \bu \cdots \bu  h^\eta_{l}$
(respectively, $h^\eta_{3} \bu v^\eta_3 \bu \cdots \bu
h^\eta_{4} \bu v^\eta_4 \bu \cdots \bu h^\eta_{l}$).
Proposition \ref{proposition:obvious disc} and, thus, Theorem \ref{theorem:main result}
will be established when the middle boundary of $\notchdisc$ has been simplified
so that it contains only a single horizontal boundary arc.
Our edge assignment assumptions give us that middle boundary edgepath for
$\partial \Delta^{\prime \eta}_e$ start on back edge $\partial$-vertical arc
and ends on a front edge $\partial$-vertical arc.
Whereas, for a positive (respectively, negative) flyping disc
$\partial \Delta^{\prime \eta}_e$ the middle boundary edgepath starts and ends on
a front edge (respectively, back edge) $\partial$-vertical arc.

Recalling our decomposition $\Delta_\v = D_{+1} \cup N$,
the reader should notice that the notching transition and the argument of Lemma
\ref{lemma:first simplification of notch disc} alters only the annulus $N$ and leaves
$D_{+1}$ untouched.  In particular, the plane $z = z_{max}$ containing $D_{+1}$
still has $z_{max}$ being greater than the horizontal position of all horizontal
arcs of $\arcpres$.  Abusing notation we will have the decomposition
$\notchdisc = D_{+1} \cup N$.

For the remainder of our discussion we will assume that our notch disc
satisfies the conclusion
of Lemma \ref{lemma:first simplification of notch disc}.

\subsection{The intersection of $C_1$ \& $\notchdisc$.}
\label{subsection:intersection of c1 and notchdisc}

In this subsection we start with a given initial pair $(\arcpres , \notchdisc)$ and
analyze the intersection $C_1 \cap \notchdisc$.
Our overall strategy is to simplify $C_1 \cap \notchdisc$ until $\notchdisc$
is an obvious disc illustrating either a destabilization, exchange move or elementary
flype.

We consider the intersection $C_1 \cap \notchdisc$.  Notice that for each
horizontal arc $h^\eta_j \subset \partial \notchdisc$ we necessarily have
$h^\eta_j \subset C_1$; for each vertical arc $v^\eta_j \subset \partial \notchdisc$
we have $v^\eta_j \cap C_1 = \partial v^\eta_j$; and we can assume that the interior of each vertical
boundary arc $v^\partial$ transversally intersects $C_1$ at finitely many points.
We can assume that $C_1$ and $int(\notchdisc)$ intersect transversely.  Thus,
$\overline{C_1 \cap int( \notchdisc )}$ is a union of simple arcs (\sa) and simple closed curves (\scc).

From our decomposition $\notchdisc = D_{+1} \cup N$ we know that we have a
distinguished \scc $\partial D_{+1} \subset C_1 \cap \notchdisc$ which we will
denote by $c_{max}$.  All other \scc and \sa of $C_1 \cap \notchdisc$ are intersections
in $C_1 \cap N$.

This next lemma allows us to get some initial
control over the behavior of $C_1 \cap \notchdisc$ without altering $\arcpres$.

\begin{lemma}[Second simplification of $(\arcpres , \notchdisc)$]
\label{lemma:control of arc intersections}
Let $(X , \Delta_\v)$ be a braid presentation/disc pair where $\v \in \{d,e,f\}$,
and consider an arc presentation/disc pair coming from the transition
$(X , \Delta_\v) \ntran (\arcpres , \notchdisc )$. 
Then we can replace
the pair $(\arcpres , \notchdisc)$ with $(\arcpres , \notchdiscprime)$
such that no \sa has an endpoint on $\arcpres \cap \partial \notchdiscprime$.
In particular, all \sa of $\notchdiscprime \cap C_1$ have their endpoints on
the $\partial$-vertical arcs of $\partial \notchdiscprime$.
\end{lemma}

\pf
By an isotopy of a collar neighborhood of $\partial \notchdisc$ in $N (\subset \notchdisc)$
we can assume that there is a neighborhood ${\bf n} \subset N \setminus c_{max}$ which has the structure $(\arcpres \cap \partial \notchdiscprime) \times I$ such that
${\bf n} \cap C_1 = \arcpres \cap \partial \notchdiscprime $.  After this isotopy
the only place where any \sa can have its endpoints is on the $v^\partial$
$\partial$-vertical arcs.
\qed

After Lemma \ref{lemma:control of arc intersections} the reader should notice that we can assume
$ \overline{C_1 \cap int( \notchdisc )} = C_1 \cap (\notchdisc \setminus \arcpres)$.

Next, we set $\cT_0 = \{(r,\theta,z) | r<1 \} $ and $ \cT_\infty = \{(r,\theta,z) | r>1 \}$.
Let $R \subset \notchdisc \setminus (C_1 \cap \notchdisc)$ be any component.
If $R \subset \cT_0$ (respectively $R \subset \cT_\infty$)
then we assign $R$ a ``$0$'' (respectively ``$\infty$'') label, i.e. $R^0$
(respectively $R^\infty$).

\begin{lemma}[Initial position of $C_1 \cap \notchdisc$-part 1.]
\label{lemma:discs that sa and scc bound}
Let $(X , \Delta_\v)$ be a braid presentation/disc pair where $\v \in \{d,e,f\}$,
and consider an arc presentation/disc pair coming from the transition
$(X , \Delta_\v) \ntran (\arcpres , \notchdisc )$. 
Then we can replace
the pair $(\arcpres , \notchdisc)$ with $(\arcpres , \notchdiscprime)$
such that the following hold:
\bi
\item[a.] Every \scc of $C_1 \cap \notchdiscprime$ bounds a subdisc of $\notchdiscprime$
whose associated label is $0$.
\item[b.] For every \scc of $c \subset C_1 \cap \notchdiscprime$
either $c$ bounds a subdisc of $C_1$ which contains a single horizontal
arc of $\arcpres$, or $c = c_{max}$ and bounds $D_{+1}$.
\item[c.] For every \sa of $C_1 \cap \notchdiscprime$ having both endpoints on the same
$\partial$-vertical arc
$v^\partial \subset H_{\theta^\partial}^{\ }$, it is outer-most in $\notchdiscprime$
and splits off a subdisc of $\notchdiscprime $ whose associated label is $0$.
\item[d.] For every \sa $\g \subset C_1 \cap \notchdiscprime$ that has both endpoints
on the same $\partial$-vertical arc $v^\partial$, there exists a sub-arc
$\g^\prime \subset H_{\theta^\partial}^{\ } \cap C_1$ with $\partial \g = \partial \g^\prime$
such that the bounded disc components of $C_1 \setminus [ \g \cup \g^\prime ] (\subset C_1)$
intersects exactly one horizontal arc of $\arcpres$.  (Note: $\g \bu \g^\prime \bu$ is not
necessarily a \scc, thus there may be more than one disc component of $C_1 \setminus (\g \cup \g^\prime)$
with only one intersecting $\arcpres$.)
\ei
\end{lemma}

\pf
Our argument for all four statements involves understanding the behavior of
$C_1 \cap \notchdisc$ in the radial foliation
of $\notchdisc$.  To start, we assume that 
all but finitely many points in the components of
$C_1 \cap \notchdisc$ are transverse to the radial foliation
of $\notchdisc$.  Moreover, we assume that the points tangency are generic
(local max or min) and each leaf in the radial foliation has at most one
point of tangency.

We first deal have a {\em simple situation}.  Let $c \subset C_1 \cap \notchdisc$
be a \scc such that: $c$ is innermost on both $C_1$ and $\notchdisc$;
$c$ bounds a disc $R^{\infty} \subset \notchdisc$; and, $c$ is tangent
at exactly two points to the radial foliation of $\notchdisc$.
Then $R^{\infty}$ splits off a $3$-ball in $\cT_{\infty}$ whose interior has
empty intersection with $\arcpres$.  We can
isotop $R^{\infty}$ through this $3$-ball so as to eliminate $c$ and reducing
$|C_1 \cap \notchdisc|$.

Next, let $p \in C_1 \cap \notchdisc$ be a point of tangency.
Suppose there exists a closed subarc in the leaf of
the radial foliation
$\g \subset \notchdisc \cap H_\theta$ such that: $p \in \g$; $\partial \g \subset C_1 \cap \notchdisc$
(we allow for $p \in \partial \g$); and,
$(\g \setminus p) \cap \cT_{\infty} = int(\g) \cap \cT_{\infty} ( \not= \emptyset )$.
Then we say $p$ is a {\em extraneous tangency}.
(We observe that in our simple situation of $c$ being a \scc having exactly two
points of tangency, neither tangency was extraneous.)

Let $E$ be the number of extraneous tangency of $C_1 \cap \notchdisc$.
Our immediate goal is to show how we can reduce $E$ to zero.
To do this we first observe that there are two types of tangent points
in $C_1 \cap \notchdisc \subset \notchdisc$. For a tangent point
$p \in C_1 \cap \notchdisc$ let $\g$, again, be a small enough subarc of the leaf
in the radial foliation of $\notchdisc$ such that: $p \in \g$ and
$(\g \setminus p) \cap C_1 = \emptyset$.  Then either
[type-1] $\g \cap \cT_{0} \not= \emptyset$ or [type-2] $\g \cap \cT_{\infty} \not= \emptyset$.

It is the type-1 where $\g \cap \cT_{0} \not= \emptyset$ that creates the
possibility of an extraneous tangency.  One can see this as follows.
Let $c \subset C_1 \cap \notchdisc$ be the \scc or \sa component with
$p \in c$.  Without loss of generality suppose that as we push $\g$
forward in the radial foliation of $\notchdisc$ it becomes a ``secant'' in $\cT_{\infty}$
near the tangency $p$, intersecting $c$ twice.  (Pushing $\g$ in the backward
direction would move $\g$ off of $c$.)  Now thinking of $\g$ as this newly formed
secant it will have its two endpoints on $c$ and it will intersect $\cT_{\infty}$.
If we continue to push our secant $\g$ forward (maintaining the feature that
its endpoints are sliding along $c$) we stop when one of three events
occurs: 1) $\g$ encounters another tangent point $p^\prime \subset C_1 \cap \notchdisc$;  2) $\g$
encounters a $\partial$-vertical arc; or 3) $\g$ shrinks down and becomes a tangent
point.  If event-3 occurs then $c$ was a \scc meeting the assumptions of
our simple situation.  For the moment we allow for the occurrence of event-2.
If event-1 occurs then the new tangent point $\g$ encounters is
extraneous.  Moreover, the configuration in event-1 satisfies the
following features.
\bi
\item[a.] $B \subset \cT_{\infty}$ is an open $3$-ball with
$\overline{B} \cap \arcpres = \emptyset$.
\item[b.] $\partial \overline{B}$ is $\overline{R_1 \cup R_2 \cup R_3}$.
\item[c.] $R_1 \subset \notchdisc \cap \cT_{\infty}$ is a subdisc with
$\partial \overline{R_1} = \g \bu \a \bu$ where:
\bi
\item[i.] $\g \subset \notchdisc \cap H_{\theta_0}^{\ }$ is an event-1 secant containing the extraneous tangent point $p^\prime$;
\item[ii.] $\a \subset c$ with $p \in \a$;
\ei
\item[d.] $R_2 \subset C_1$ is a subdisc
with $\partial \overline{R_2} = \a \bu \a^\prime \bu$ where
$\a$ is from c-ii and $\a^\prime \subset C_1$ is a subarc of $C_1 \cap H_{\theta_0}^{\ }$.
\item[e.] $R_3 \subset H_{\theta_0}^{\ } \cap \cT_{\infty} $ is a subregion
with $\partial \overline{R_3} = \g \bu \a^\prime \bu$.  (Depending on whether
$p$ is in $\partial \g$ or $int(\g)$, $R_3$ is either a subdisc or the wedge of
two subdiscs with $p^\prime$ as the wedge point.)
\ei  
We are now in a position we see how $E$ can be reduced to zero.
We consider the disc $R_3 \subset H_{\theta_0}^{\ }$.  Since
$\g \subset \overline{R_3}$, see that $p^\prime \subset \overline{R_3}$.  Assuming
that $R_3 \cap \notchdisc = \emptyset$, in a product neighborhood of $R_3$
we can push $\g$ into $\cT_0$ dragging $\notchdisc$ along to eliminate $p^\prime$
as an extraneous tangent point.
(The key feature of this isotopy is that the radial foliation of $\notchdisc$
is unchanged.)
For the new $\notchdiscprime$,
$C_1 \cap \notchdiscprime$ will have a new simple situation \scc.
If $p^\prime$ had been a boundary endpoint of $\g$ then
the number of tangent points in $C_1 \cap \notchdiscprime$ is the same.
If $p^\prime$ had been in $int(\g)$ then $C_1 \cap \notchdiscprime$ has two additional
tangent points in the radial foliation.  However, $E$ has gone down by one.

For $R_3 \cap \notchdisc \not= \emptyset$ we recall that $p^\prime$
is the only tangent point in $H_{\theta_0}^{\ }$.  So $R_3 \cap \notchdisc$ is a collection
of \sa.  Starting with an outermost such \sa, in a product neighborhood of
$R_3$ we push these \sa's into $\cT_0$, dragging $\notchdisc$ along.  It is
readily seen that no new extraneous tangent points are introduced.
Iterating this procedure we can assume $E=0$.  Moreover,
after eliminating simple situation \scc we can assume
statement-a is true.
%

To establish statement b. let $c \subset C_1 \cap \notchdisc$.  Assume that $c$ bounds
a sub-disc $R \subset \notchdisc$.  Now that we have established the validity of
statement a. we can assume that $int(R) \subset \cT_0$.
If $R \cap \axis \not= \emptyset$ then $R$ is in fact $D_{+1}$ and we have
$c = c_{max}$.  Moreover, $c_{max}$ is unique in having the feature that it bounds
a disc which intersects $\axis$

Now assume $R \cap \axis = \emptyset$.  Thus, $c$ must bound on $C_1$.
Let $\Delta_c \subset C_1$ having $\partial \Delta_c = c$.
Moreover, $c$ contains exactly two tangent points for otherwise $c$ would
contain an extraneous tangent point.  This allows us to assume
that we can isotop $R$, first, arbitrarily close to $C_1$, then
onto a subdisc of $C_1$ by radially pushing $R$ out of $\cT_0$
along $\theta$-rays.  The obstruction to pushing $R$ totally out of $\cT_0$ will be
horizontal arcs of $\arcpres$ that are contained in $\Delta_c$.  Thus, this radially
isotopy of $R$ will result in a new \scc for each horizontal arc of
$\Delta_c \cap \arcpres$.

The arguments establishing statements c. \& d. are similar.  For statement c. let
$\a \subset C_1 \cap \notchdisc$ be a \sa that splits off a region $R \subset \notchdisc$.
Assume that $\partial \a $ is contained in the $\partial$-vertical arc
$v^\partial \subset H_{\theta^\partial_{\ }}^{\ }$ where
$H_{\theta^\partial_{\ }}^{\ } \in \fib$.
The first two possibilities to consider is either $R \cap \cT_{\infty} = \emptyset$
or $R \cap \cT_{\infty} \not= \emptyset$.

Suppose $R \cap \cT_{\infty} \not= \emptyset$ and assume that $E=0$ and all
simple situation \scc's have been eliminated.  Then $\a$ contains either one or three
tangent points with the radial foliation of $\notchdisc$.  The
reader should readily see that the arc $v^\partial \cap R (\subset H_{\theta^\partial_{\ }}^{\ })$
splits off a subdisc of $\cT_{\infty} \cap H_{\theta^\partial_{\ }}^{\ }$ that does in
intersect $\arcpres$.  Thus, we can push this arc into $\cT_0$ reducing
$|C_1 \cap v^\partial|$.  We eliminate any extraneous tangencies or simple situation \scc that are produced.
Ultimately, this will result in pushing the $R$-portion of $\notchdisc$ into
$\cT_0$ thus reducing $|C_1 \cap \notchdisc|$.

Now suppose $R \cap \cT_{0} = \emptyset$ and assume that $E=0$.
It is easily seen that $\a$ contains either [case-1] one or [case-2] three points of tangencies.
For case-1 the one point of tangency must be a type-2.
For case-2 two of the tangencies must be type-2 and one a type-1.

Considering case-1, $\a$ has exactly one point of tangency with the radial
foliation of $\notchdisc$.  To be descriptive, the angular position of $\g ( = R \cap v^\partial)$ is an
endpoint of the interval that is the angular support of $\a$; and,
the angular position of the type-2 point of tangency is
the other endpoint of the angular support of $\a$.
The vertical support of $\a$ and $\g$ are equal.

Recall that
$\g \subset v^\partial \subset H_{\theta^\partial_{\ }}^{\ }$.
Let $\g^\prime \subset C_1 \cap H_{\theta^\partial_{\ }}^{\ }$ such that
$\g^\prime \bu \a \bu$ bounds a subdisc $R^\prime \subset C_1$.
In $\cT_0$ we can isotop $R$ close enough to $C_1$ such that
by radially pushing $R$ out of $\cT_0$
along $\theta$-rays $R$ is bijectively mapped onto $R^\prime$.
The obstruction to pushing $R$ totally out of $\cT_0$ will be
horizontal arcs of $\arcpres$ that are contained in $R^\prime$.  Thus, this radially
isotopy of $R$ will result in a new \sa's for each horizontal arc of
$R^\prime \cap \arcpres$.
Thus, for case-1 \sa's we obtain statements c. and d.

Statements c. and d. for case-2 \sa is achieve in almost similar fashion.
Notice that with $E=0$, our arc $\a$ will have exactly three points
of tangency.  Being descriptive, we can
traverse $\a$ so that we have, in order points of tangency,
$p_1 , p_2 , p_3 \subset \a$ where $p_1$ is a type-1 tangency that leads to
an event-2 secant-push.
The angular position of $\g$ \& $\g^\prime$ is in the interior
of the interval that is the angular support of $\a$.  And, the vertical
support of $\a$ properly contains the vertical support of $\g$ \& $\g^\prime$.
(For measures of vertical and angular support it is best to visualize
these support projected onto $C_1$.)
We can move the angular position of $p_1$ and $p_2$ arbitrarily close to that
of $\g$ \& $\g^\prime$ thus making the vertical support of $\a$ arbitrarily close
to that of $\g$ \& $\g^\prime$.  (Due to the fact that we are maintaining
the constant nature of the edge assignment of $v^\partial$, the geometry
of $R$ is creating a ``fold'' near $\g$.)

The closed curve $\g^\prime \bu \a \bu$  is either simple or has two points of
self intersection.  Either way $C_1 \setminus [ \g^\prime \bu \a \bu ]$ has bounded disc components.
In $\cT_0$ we can isotop $R$ close enough to $C_1$
by radially pushing $R$ out of $\cT_0$
along $\theta$-rays.
(The image of $R$ onto $C_1$ under this radial push will always be more than
a disc region due to the previously mentioned fold.) 
The obstruction to pushing $R$ and $\g$ totally out of $\cT_0$ will be
horizontal arcs of $\arcpres$ that are contained in the bounded disc components
of $C_1 \setminus [ \g^\prime \bu \a \bu ]$.  Thus, this radially
isotopy of $R$ will result in a new \sa's for each horizontal arc of
$R^\prime \cap \arcpres$.
Thus, for case-1 \sa's we obtain statements c. and d.
\qed

\begin{figure}[htb]
\centerline{\includegraphics[scale=.75, bb=0 0 471 208]{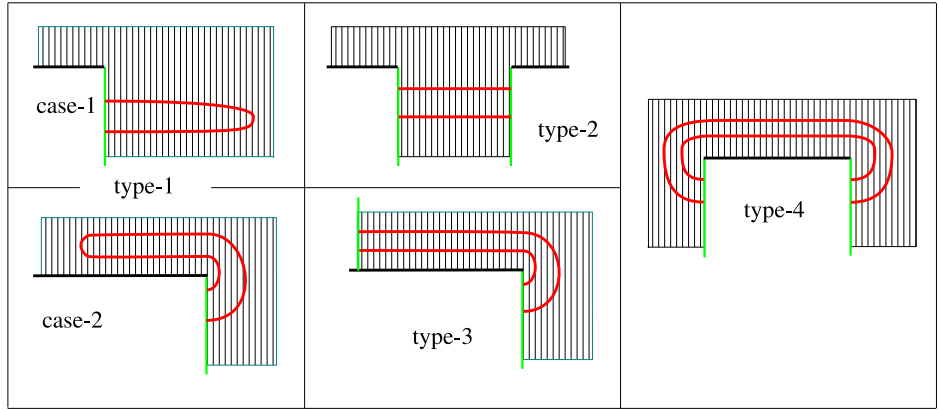}}
\caption{The vertical cross-hatching illustrates the local leaves of the radial foliation.
The $\g$ \sa's are the red arcs.}
\label{figure:sa arc types}
\end{figure}

Statements c \& d of Lemma \ref{lemma:discs that sa and scc bound} deal with
just one type of \sa which we will be interested in.  In particular, there are
four types of \sa's in our argument.  All will occur in angular intervals that is less
than $2 \pi$.  (See Figure \ref{figure:sa arc types}.)

\bi
\item[type-1] $\g \subset C_1 \cap \notchdiscprime$ is a \sa
having both endpoints on the same $\partial$-vertical arc
$v^\partial \subset \partial \notchdisc$.  From the proof of
Lemma \ref{lemma:discs that sa and scc bound} we know we have case-1 (one point of tangency) \& case-2 
(three points of tangency) for $\g$.
\item[type-2] $\g \subset C_1 \cap \notchdiscprime$ is a \sa
having its endpoints on different $\partial$-vertical arc
$v^\partial, v^{\prime \partial} \subset \partial \notchdisc$.
Moreover, $\g$ is transverse to the leaves of the radial foliation of $\notchdisc$.
Thus, the edge assignments of $v^\partial$ and $v^{\prime \partial}$
must be different.
\item[type-3] $\g \subset C_1 \cap \notchdiscprime$ is a \sa
having its endpoints on different $\partial$-vertical arc
$v^\partial, v^{\prime \partial} \subset \partial \notchdisc$.
Moreover, $\g$ has one point where it is not
transverse to the radial foliation of $\notchdisc$.
Thus, the edge assignments of $v^\partial$ and $v^{\prime \partial}$
must be the same.
\item[type-4] $\g \subset C_1 \cap \notchdiscprime$ is a \sa
having its endpoints on different $\partial$-vertical arc
$v^\partial, v^{\prime \partial} \subset \partial \notchdisc$.
Moreover, $\g$ has two points where it is not
transverse to the radial foliation of $\notchdisc$ with $\g$ transversely intersecting
any leaf at most twice.  
Thus, the edge assignments of $v^\partial$ and $v^{\prime \partial}$
must be different.
\ei

\begin{lemma}[Initial position of $C_1 \cap \notchdisc$-part 2.]
\label{lemma:discs with label 0}
Let $(X , \Delta_\v)$ be a braid presentation/disc pair where $\v \in \{d,e,f\}$,
and consider an arc presentation/disc pair coming from the transition
$(X , \Delta_\v) \ntran (\arcpres , \notchdisc )$. 
We can replace
the pair $(\arcpres , \notchdisc)$ with $(\arcpres , \notchdiscprime)$
such that for every component of
$\delta \subset \notchdiscprime \setminus (C_1 \cap \notchdiscprime)$ with label $0$,
if $\delta \cap \axis = \emptyset$ and $\delta \cap \arcpres = \emptyset$ then
$\delta$ is a disc of the one of the following type:
\bi
\item[a.] {\rm A whole disc}--A subdisc whose boundary in $C_1 \cap \notchdiscprime$ is a \scc.
\item[b.] {\rm A half disc}--A subdisc whose boundary is the cyclic ordered
union of two arcs $\g_1 \bu \g_2 \bu$
where $\g_1 \subset C_1 \cap \notchdisc$ is a type-1 \sa; and,
$\g_2 \subset v^\partial$ for some
some $\partial$-vertical arc $v^\partial$.
\item[c.] {\rm a rectangle}--A subdisc those boundary is the cyclic ordered union of
four arcs $\g_1 \bu \g_2 \bu \g_3 \bu \g_4 \bu$
where: $\g_1 , \g_3 \subset C_1 \cap \notchdisc$;
$\g_1$ \& $\g_3$ are \sa's of the same type-j in $\notchdisc$, $j \in \{ 2,3,4 \}$;
$\g_2 \subset v^\partial_1$ and
$\g_4 \subset v^\partial_2$ where $v^\partial_1 (\subset H_{\theta_1^{\ }}^{\ })$
\& $v^\partial_2 (\subset H_{\theta_2^{\ }}^{\ })$ are two
different $\partial$-vertical arcs.  Moreover, there exists
a unique horizontal arc $h \subset \arcpres$ whose angular support contains the
angular interval over which $\delta$ occurs; and, whose horizontal position is contained
in the horizontal interval over which $\delta$ occurs.
\ei
\end{lemma}

\pf
Consider a component $ \delta^0 \subset \notchdisc \cap \cT_0$
where $\delta^0 \cap \axis = \emptyset$.  If $|\partial \overline{\delta^0} \cap C_1| = 1$
then either $\partial \overline{\delta^0} \cap C_1$ is a \scc or a \sa.  If it is a \scc
then by statements a \& b of Lemma \ref{lemma:discs that sa and scc bound} we have our whole-disc-statement a.
If it is a \sa then this \sa must be a case-1 or -2 type-1 \sa and
by statements c \& d of Lemma \ref{lemma:discs that sa and scc bound}
we have our half-disc-statement b.

More generally, if $|\partial \overline{\delta^0} \cap C_1| > 1$ then, again,
by statements a \& c of Lemma \ref{lemma:discs that sa and scc bound}, $\delta^0$ is a disc
planar region of $\notchdisc$.  Moreover, since we
can assume there are no extraneous tangencies, we conclude that
$\partial \overline{\delta^0} \cap C_1$ is a collection of type-2, -3 \& -4 \sa's.
Also, since there are no extraneous tangencies we can place $\delta^0$ arbitrarily close to $C_1$
so that any radial ray, $\{(r,0,0) | r \geq 0 \}$, intersects $\delta^0$ at most once.
Thus, we can consider the attempt to radially push $\delta^0$ out of
$\cT_0$.  The obstructions to pushing all of $\delta^0$ out of $\cT_0$
are any horizontal arcs of $\arcpres$ that intersect the radial projection $\pi (\overline{\delta^0}) \subset C_1$.
We list the possibilities.  If $h \subset \arcpres \cap \pi (\overline{\delta^0})$
is a horizontal arc in the interior of
$\pi (\overline{\delta^0})$ then $h$ obstructs a whole disc in the radial push of $\delta^0$ out of $\cT_0$.
If $h \subset \arcpres$ is a portion of a horizontal arc that intersects $\partial \pi (\overline{\delta^0})$
once then $h$ obstructs a half disc in the radial push of $\delta^0$ out of $\cT_0$.
If $h \subset \arcpres$ is a portion of a horizontal arc that intersects $\partial \pi (\overline{\delta^0})$
twice then $h$ obstructs a rectangle in the radial push of $\delta^0$ out of $\cT_0$.
(See Figure \ref{figure:sa arc types} for illustrations of the three types of rectangles.)
\qed

In Lemmas \ref{lemma:control of arc intersections},
\ref{lemma:discs that sa and scc bound} and \ref{lemma:discs with label 0}
for the pair $(\arcpres, \notchdisc)$ we did
not alter $\arcpres$, only $\notchdisc$.  Thus, our complexity measure
$\cC ( \arcpres)$ remained constant.
The resulting pair $(\arcpres , \notchdiscprime)$ coming from the application
of Lemmas \ref{lemma:discs that sa and scc bound} and \ref{lemma:discs with label 0}
will be referred to as an
{\em initial position} for the braid presentation/disc pair.

\begin{rem}
\label{remark:on rectangle discs}
{\rm In \S\ref{subsection:the tiling machinery on notchdisc} we will develop the
machinery for eliminating the whole discs and half discs of Lemma \ref{lemma:discs with label 0}.
The reader should notice that rectangle disc of statement c. of Lemma
\ref{lemma:discs with label 0} only occur for $\Delta_e$ and $\Delta_f$ discs.
It is possible to replace $(X , \Delta_\v)$ and a corresponding
$(X , \Delta_\v) \ntran (\arcpres , \notchdisc )$, $ \v \in \{e,f\}$, with a new
$(\arcpres , \notchdiscprime)$
such that no component of
$\notchdiscprime \cap \cT_0$ is a rectangle subdisc.  The replacement
will still characterize the occurrence of an
exchange move or flype.  However, in line with our discussion in
Remark \ref{remark:elementary exchange moves and flypes}, this may not be the same exchange move or flype as characterized
by the original $\notchdisc$.  Thus, in order to establish Theorem
\ref{Theorem:destably equivalent} we will need to maintain the integrity of our $\Delta_e$
and $\Delta_f$ discs.

Finally, consider $\g_1$ \& $\g_3$ of statement c. 
having endpoints on $\partial$-vertical arcs $v_1^\partial$ and $v_2^\partial$, where
$\theta_1$ and $\theta_2$ are their angular positions, respectively.
Then the angular interval over which $\g_1$ \& $\g_3$
occur contains, say, $[\theta_1 , \theta_2]$.
(The other possibility is $[\theta_2 , \theta_1]$.)  
By a slight isotopy we can assume that the
$\g_1 ,\g_3$ have constant $z$-coordinates over
$[\theta_1 , \theta_2]$ and we will refer to these constant coordinates as the
horizontal positions of $\g_1$ \& $\g_3$.  We will refer to them
as {\em horizontal boundary arcs} and $\g_2, \g_4$ as {\em $\partial$-vertical arcs} of the rectangle.  Since $[\theta_1, \theta_2]$ will be arbitrarily close to the
angular interval over which $\g_1$ \& $\g_3$ occur (for type-2 it is equal) we will refer to
it as the angular support of these \sa's and the associated rectangle $\delta$.
Similarly, the vertical interval over which $\delta$ occurs is arbitrarily close to the
interval between the horizontal positions of $\g_1$ and $\g_3$.
We will refer to this interval as the vertical support of $\delta$.
We will let $\cR \subset \notchdisc$ be the set of all
rectangle subdiscs.}
\qed
\end{rem}

\subsection{The tiling machinery on $\notchdisc$.}
\label{subsection:the tiling machinery on notchdisc}
We now make the transition from an arc presentation to a shear presentation.
To do this we need to specify the intervals for $\cI$.
Each $\partial$-vertical arc
of $\Delta_\v$ will contribute one component to $\cI$.
So let $v^\partial_k \subset \partial \Delta_\v$ 
be a $\partial$-vertical arc having angular position $\vartheta_k$.
Then $\cI$ contains an angle interval
$[\vartheta_k - \e_k, \vartheta_k + \e_k ]$.
($\e_k$ is small enough so that the components of $\cI$
contains no vertical arcs and are pairwise disjoint; and,
$1 \leq k \leq p$ where $p$ is $1$, $2$ or $3$ depending upon whether
$\v$ is $d$, $e$ or $f$, respectively.) 
With this assignment for $\cI$ in place we introduce new transition notation
$(X , \Delta_\v ) \stran (\sheararcpres , \notchdisc)$ for going from
a braid/disc pair to a shear-arc-presentation/notch-disc pair. 

We now adapt the classical foliation/tiling machinery of surfaces in braid structures
that has been extensively developed and exploited in \cite{[BM1], [BM2], [BM3], [BM4], [BF]}.  This
machinery has built into it the ability to recognize when and where exchange moves
are admitted.  It is important to notice that we will be
creating a tiling only on the components of $\notchdisc \setminus int(\cR)$.
To start our adaptation we have the following definitions.

Given a pair $(\sheararcpres , \notchdisc)$ in initial position (as specified
by the conclusions of Lemma \ref{lemma:discs that sa and scc bound} \&
\ref{lemma:discs with label 0}),
let $ h \subset \sheararcpres$ be
a horizontal arc (respectively, portion of a horizontal arc)
that corresponds to a horizontal arc (respectively, portion of a horizontal arc)
associated with a whole disc as described in statement b.
(respectively, half disc as described in statement c.) of
Lemma \ref{lemma:discs that sa and scc bound}.

Let $[\theta^h_1 , \theta^h_2]$ be the angular support of $h$ and $z^h$ be the vertical
position of $h$.  The {\em coning disc of $h$} is the disc
$\delta^h = \{ (r,\theta,z^h) | 0 \leq r \leq 1 , \theta^h_1 \leq \theta \leq \theta^h_2 \}$.
Let $\ver^h = \delta^h \cap \axis$ is the vertex for $\delta^h$.
The reader should notice that a coning disc inherits a radial foliation from its
intersection with the disc fibers of $\fib$.

\begin{figure}[htb]
\centerline{\includegraphics[scale=.43, bb=0 0 510 608]{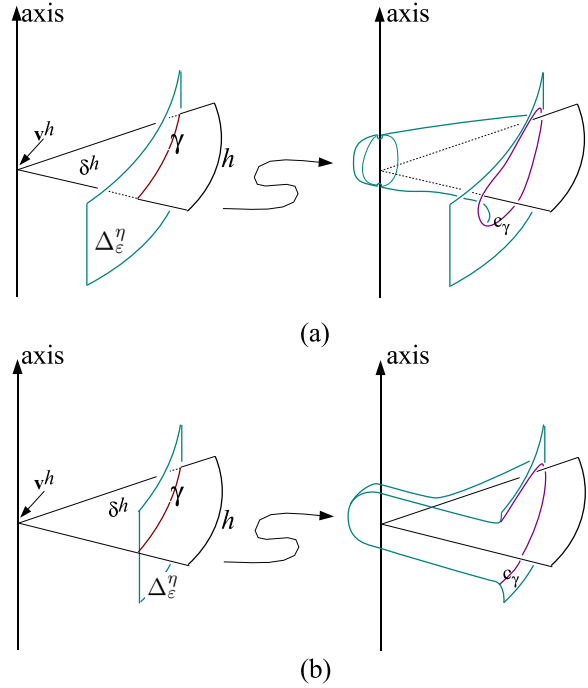}}
\caption{The corresponding alteration in the foliation in (a) is
depicted in Figure \ref{figure:tiles}(a).  Similarly, there is a correspondence between
illustration (b) and Figure \ref{figure:tiles}(b).  The curve $c_\gamma$ is used to
designate the grouping of the vertices and singularities introduced by the tiling
isotopy that is associated with the intersection arc $\gamma$.}
\label{figure:tiling isotopy}
\end{figure}

Now for any coning disc $\delta^h$ we consider the intersection set $\notchdisc \cap \delta^h$.
By the four Lemmas used to define initial position we know that this intersection set
will be a union of arcs that are transverse to disc fibers of $\fib$ and, thus, the
leaves of the radial foliation of $\delta^h$.  So viewed in $\delta^h$, any arc of intersection
with $\notchdisc$ will be seen as parallel to $ h \subset \partial \delta^h$.  In particular,
by a slight isotopy of the whole discs and half discs of
$\notchdisc \cap \cT_0$, we can assume that the cylindrical coordinates of any
intersection arc in the set $\notchdisc \cap \delta^h$ has constant $r$-coordinate along
with constant $z$-coordinate being $z^h$.

\begin{figure}[htb]
\centerline{\includegraphics[scale=.65, bb=0 0 713 333]{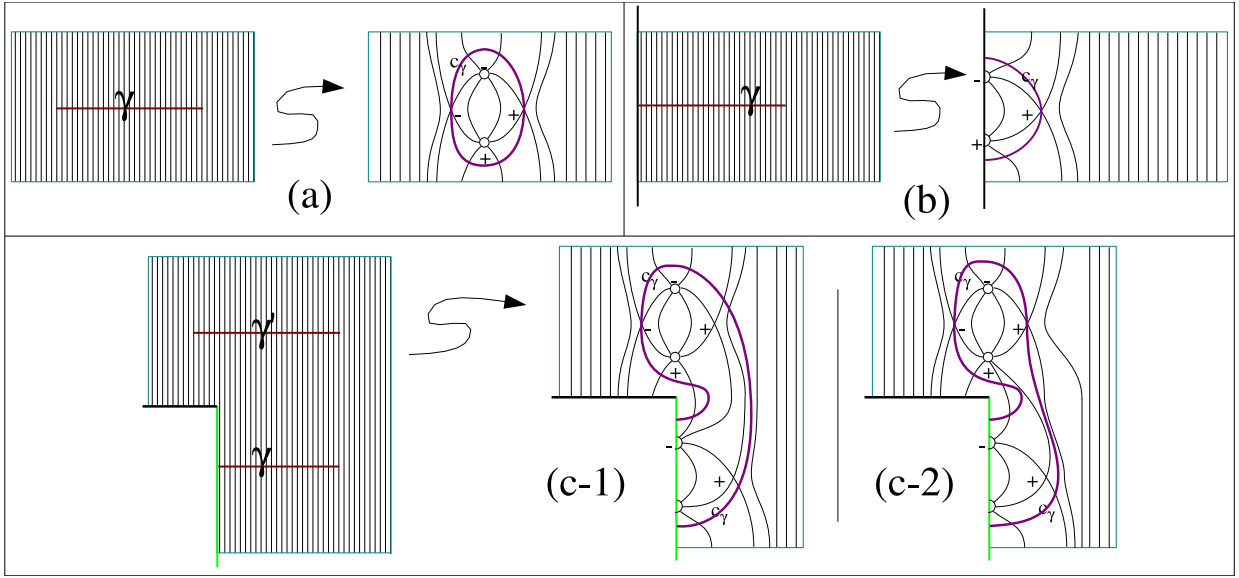}}
\caption{The parity signs indict one possibility for the parity assignment of
singularities and vertices.  In illustration (a) there is always the occurrence
of a positive/negative pair of vertices and singularities.  In illustration (b)
there is always the occurrence of a positive/negative pair of vertices.  The lone
singularity can be of either parity.  The curves
$c_\gamma(\subset C_1 \cap \tilenotchdisc$) is useful in
grouping of the vertices and singularities introduced by the tiling isotopy associated
with the intersection arc $\gamma$.}
\label{figure:tiles}
\end{figure}

Let $\gamma \subset \notchdisc \cap \delta^h$ be the intersection arc in $\delta^h$ having
the smallest $r$-coordinate which we call $r_\gamma$.  We can then isotop $\gamma$ through
$\delta^h$ and past the axis $\axis$ by letting $r_\gamma$ go to zero and past $\axis$.
Extending this isotopy of $\gamma$ to $\notchdisc$ we
produce a disc $\tilenotchdisc$ that has a {\em tiled foliation}.

As illustrated
in Figure \ref{figure:tiling isotopy}-right, there are two possibilities:
$c_\g$ is a \scc or $c_\g$ is a \sa.  However, since there are case-1 and case-2 \sa's, the type-1 \sa
situation
further bifurcates into two slightly differing alterations to the foliation.  In Figure \ref{figure:tiling isotopy}
we illustrate the {\em tiling isotopy} for a coning disc $\delta^h$ associated with
a whole disc (Figure \ref{figure:tiling isotopy}(a)) and a half disc split off
by a case-1 type-1 \sa (Figure \ref{figure:tiling isotopy}(b)).

As mentioned before, these tiling isotopies will occur away from the
rectangle discs (statement c. Lemma \ref{lemma:discs with label 0}) of our initial position. 
In Figure \ref{figure:tiling isotopy}(a) we have the case where $h$ is the
horizontal arc associated with a whole disc.  Here the tiling isotopy
introduces two intersection points or {\em vertices} of $\notchdisc$ with $\axis$; and
two points of tangency with disc fibers of $\fib$ or {\em singular points}.
Turning to Figure \ref{figure:tiles}(a) we see how the tiling isotopy alters the radial foliation in
a neighborhood of such an intersection arc $\gamma$.

The situation for half discs is similar but a little more complex since
the associated \sa splitting off the half disc in $\notchdisc$
can be either a case-1 or case-2, type-1.  Taking the case-1 first, in Figure
\ref{figure:tiling isotopy}(b)-left we illustrate a portion, $h$, of a horizontal
arc that is associated with said half disc.
The previously described tiling isotopy will isotopy one of the $\partial$-vertical arcs of
$v^\partial$ of $\sheararcpres$.  This will introduce two points of intersection
of $v^\partial$ with $\axis$ and one tangency point with a disc fiber of $\fib$
which we again call a {\em singular point}.  Figure \ref{figure:tiling isotopy}(b)-right
has this corresponding alteration to the foliation of $\notchdisc$.  For convenience
we will assume in all situations that the singularities introduces in tiling isotopy
does not occur in a disc fiber of $\fib$ that contains a vertical arc or $\partial$-vertical
arc.  Again, turning to Figure \ref{figure:tiles}(b) we see how the tiling isotopy alters the radial foliation in
a neighborhood of such an intersection arc $\gamma$.

When half disc $h$ is associated with is split off by a case-2 type-1 \sa then the coning disc $\delta^h$
will intersect $\notchdisc$ twice due to the ``fold'' in $\notchdisc$ near the $\partial$-vertical
arc.  Figure \ref{figure:tiles}(c)-left illustrates the local intersecting
arcs in the radial foliation of $\notchdisc$.  This situation, again bifurcates into two possibilities
since either $\g^\prime$ (the arc adjacent to a $\partial$-vertical arc)
is closest to $\axis$ or $\g$ (the arc not adjacent to a $\partial$-vertical arc) is closest to $\axis$.
Pushing these intersecting arcs through
$\delta^h$ yields to a change in foliation that corresponds to either
Figure \ref{figure:tiles}(c-1) (when $\g$ is nearest to $\axis$) or Figure \ref{figure:tiles}(c-2)
(when $\g^\prime$ is nearest to $\axis$).

In general, for a {\underline {single choice}} of a coning disc $\delta^h$ we will use the notation
$ (\sheararcpres,\notchdisc) \Ttran (\sheararcpres,\tilenotchdisc)$ to indicate
the tiling isotopy between $ (\sheararcpres,\notchdisc) $ and
$ (\sheararcpres,\tilenotchdisc) $.  (Our arc presentation and shearing intervals remain unchanged.)
The \scc/\sa illustrated in Figure \ref{figure:tiling isotopy}
corresponds to $c_\gamma \subset C_1 \cap \tilenotchdisc$. Specifically,
when $c_\gamma$ is a circle it will encircle the two vertices and pass through the
two singularities.  When $c_\gamma$ is a case-1 type-1 \sa it will split-off the two vertices that
are on the $\partial$-vertical arc and will pass through the single singularity.
When $c_\gamma$ is a case-2 type-1 \sa the disc it splits off will have four
vertices---two on a $\partial$-vertical arc and two in the interior of $\tilenotchdisc$---and
three singularities.
In Figure \ref{figure:tiles} we see how the tiling isotopy alters the radial foliation in
a neighborhood of the intersection arc $\gamma$.

\begin{figure}[htb]
\centerline{\includegraphics[scale=.9, bb=0 0 386 128]{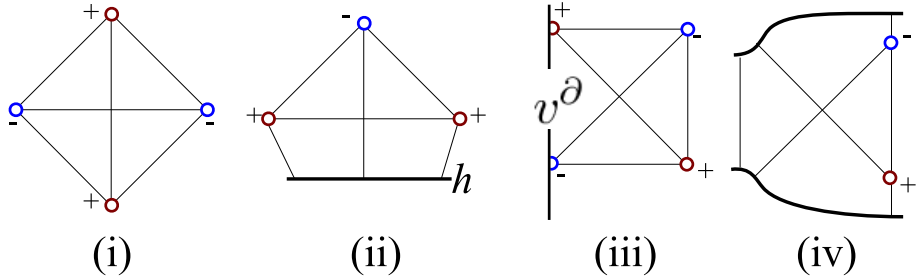}}
\caption{There are three possible means by which a singularity can be formed.  In
(i) two $b$-arcs (leaves in the foliation that have both endpoints on vertices)
can come together to form a singularity.  In (ii) the neighborhood of the singularity
has both $a$-arcs (leaves having an endpoint on a vertex and an endpoint on the boundary
of $\tilenotchdisc$) and $b$-arcs.  In (iii) the singularity has its vertices endpoints
on a $\partial$-vertical arc.  In (iv) the graphs in $bs$-singularities are illustrated.}
\label{figure:graphs}
\end{figure}

To review a little more thoroughly the tiling machinery, we consider the orientations
of the axis $\axis$, our disc $\tilenotchdisc$ and the disc fibers of $\fib$ which is
consistent with the orientation of $\axis$.  The orientation of $\tilenotchdisc$
is consistent with the orientation of the horizontal/vertical arcs on its boundary.
We can associate to each
vertex and singularity a parity as follows.  A vertex $\ver \subset \tilenotchdisc \cap \axis$
is {\em positive} (or $+$) if $\ver$ is a positive intersection.
Otherwise, $\ver$ is {\em negative} (or $-$).
A singular point $\sin \subset \tilenotchdisc$ is {\em positive} (or $+$)
if the orientation of the tangent plane
to $\notchdisc$ at $\sin$ agrees with the orientation of the disc fiber of $\fib$ that
contains $\sin$.  Otherwise, $\sin$ is {\em negative} (or $-$).
Thus, parity labeling of the vertices
in Figure \ref{figure:tiles} are the only possible assignments, whereas the parity labeling
of the singularities can either be as indicated or reversed.

Keeping with the literature that has developed around tiled foliations,
on the components of $\tilenotchdisc \setminus int(\cR)$ non-singular leaves must be
arcs that have their endpoints being either vertices of the foliation, or points in the set
$(\sheararcpres \cap C_1) \cup (\partial \cR \cap C_1)$, i.e. horizontal arcs.  
Due to reasons of orientation it is easily established that
there are three types of generic leaves in our tiled foliation: $a$-arcs which
have one endpoint being a vertex and one endpoint in the set
$(\sheararcpres \cap C_1) \cup (\partial \cR \cap C_1)$;
$b$-arcs which have both endpoints being
(differing) vertices; and, $s$-arcs which have endpoints on differing horizontal arcs
of $(\sheararcpres \cap C_1) \cup (\partial \cR \cap C_1)$.

Again, for reasons of orientation it is easily established that there are
only three types of singularities.  (See Figure \ref{figure:graphs}.)
An $ab$-singularity or $ab$-tile is formed by an $a$-arc and $b$-arc
coming together as illustrated in (ii) of
Figure \ref{figure:graphs}.  A $bb$-singularity
or $bb$-tile is formed by two $b$-arcs coming together as illustrated in
(i) \& (iii) of Figure \ref{figure:graphs}.  Finally, a $bs$-singular or $bs$-tile
is formed by a $b$-arc and $s$-arc coming together.  Notice that we could also refer to
this singularity as an $aa$-singular or $aa$-tile since the two $a$-arcs, a $b$-arc and
an $s$-arc make up the four sides of the tile.   Illustration (iv) depicts
this third singularity.

%
We note that one subdisc component of
$\tilenotchdisc \setminus int(\cR)$ will contain $D_{+1}$.
All the remaining components, for which $D_0$ will be convenient notation,
have $0$ algebraic intersection with $\axis$.

We now have the following tiling operations.

\noindent
{\bf $\cO_1$--Simplification moves on the middle boundary.}
We start with a simply alteration to horizontal arcs.

\begin{figure}[htb]
\centerline{\includegraphics[scale=.7, bb=0 0 316 151]{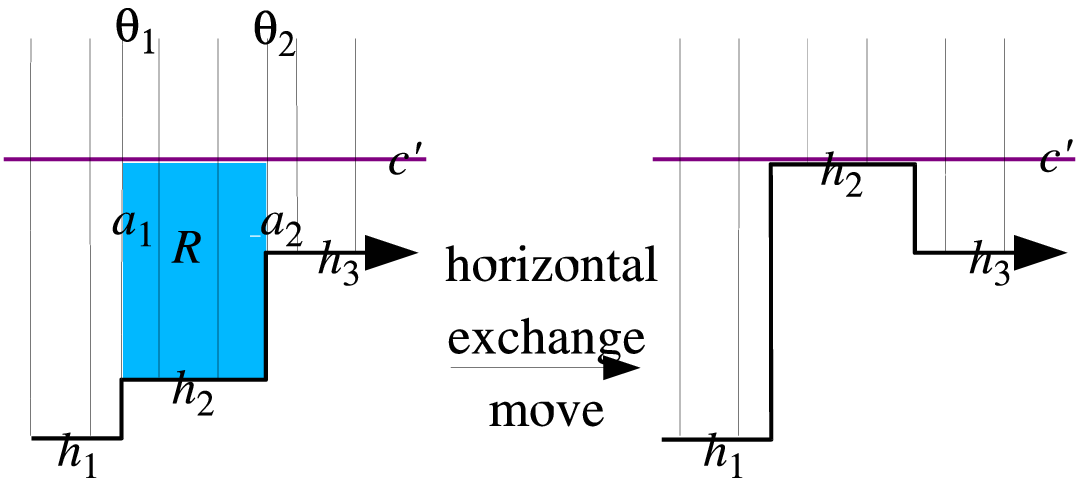}}
\caption{The sequence illustrates
a horizontal exchange move that alters $\tilenotchdisc$ near $c^\prime$.}
\label{figure:boundary exchange moves}
\end{figure}

\noindent
\underline{$\cO_{1.1}$--Horizontal arc near a horizontal boundary
arc of $\cR$ or $c_{max}$.}
We refer to Figure \ref{figure:boundary exchange moves} for this configuration
and appeal to its labeling.
Let $ R \subset \notchdisc$ be a rectangular region that satisfies the following:
i) $int(R) \cap C_1 = \emptyset$; ii) $\partial R = a_1 \bu c^\prime \bu a_2 \bu h_2 \bu$
where; iii) $a_1$ and $a_2$ are subarcs of leaves containing vertical arcs,
$v_1$ and $v_2$, respectively, of $\sheararcpres$;
iv) $c^\prime$ is an subarc of either a horizontal boundary arc of a rectangle or
$c_{max}$; and, v) $h_2 \subset \sheararcpres$ is a
horizontal arc.
Figure \ref{figure:boundary exchange moves} implies an isotopy of $h_2$ through $R$ so as to
place the horizontal position of the resulting horizontal arc nearest that of $c^\prime$.
Such an isotopy will leave all other features of $\notchdisc \cap C_1$ unaltered.  We wish to show
that this isotopy can be achieved through a sequence of our elementary moves.

The horizontal position of $v_1$ and $v_2$ leads to four cases since each vertical
arc can have vertical support that is either above or below the horizontal
position of $h_2$.  The case illustrated in Figure \ref{figure:boundary exchange moves}
has the vertical support of $v_1$ below $h_2$ whereas $v_2$ has vertical support
above the horizontal position of $h_2$. 

In Figure \ref{figure:boundary exchange moves} the angle
at which $a_1$ occurs is $\theta_1$ and the angle at which $a_2$ occurs is $\theta_2$
with $[\theta_1 ,\theta_2]$ being the angular support of $R$.
Having $c^\prime$ being a subarc of either a horizontal boundary arc of a rectangle of $\cR$ or
of $c_{max}$ allows us to speak of the horizontal position of $c^\prime$
and the vertical support of $R$ (which is between the horizontal
positions of $c^\prime$ and $h_2$).

Now we observe that the existence of $R \subset \notchdisc$ forces the vertical
support of every vertical arc of $\sheararcpres$
having angular position in $(\theta_1 ,\theta_2)$
to be either totally below, totally above, contained in,
or properly containing the vertical support of $R$.
(If $c^\prime$ is in $c_{max}$ then there is no `totally above'.)
Thus, for any two consecutive vertical arcs
$v , v^\prime \subset \sheararcpres$ having angular position in $(\theta_1 ,\theta_2)$,
if $v$ has vertical support contained in the vertical support of $R$ but
$v^\prime$ does not then $ v$ \& $ v^\prime$ are nested and we can move $v^\prime$
(forward/backward) past $v$.  
Thus, the only obstruction to moving vertical arcs having angular position
in $(\theta_1 , \theta_2)$ past $\theta_1$ (in the backward direction) or
past $\theta_2$ (in the forward direction) are the vertical arcs $v_1$ and $v_2$.
This is where our previously mention four cases come into play.
If $v_1$ (respectively, $v_2$) has vertical support below
(respectively, above) $h_1$ (as illustrated in
\ref{figure:boundary exchange moves}) then we can move all vertical arcs that
are above (respectively, below) $h_1$ with angular position in $(\theta_1 ,\theta_2)$
backwards (respectively, forward) and past $v_1$ (respectively, $v_2$).
The remaining three cases are easily listed and we leave them to the reader.
  
The conclusions are that after some number of vertical
exchange moves we can assume that over the angular interval
$(\theta_1 ,\theta_2)$ either: 1) there are not vertical arcs of $\sheararcpres$;
2) all vertical arcs of $\sheararcpres$ have vertical support below $h_2$;
or 3) all vertical arcs of $\sheararcpres$ have vertical support contained in the
vertical support of $R$.

With these conclusions in place we now see how the isotopy depicted in
Figure \ref{figure:boundary exchange moves} is realized through a
sequence of horizontal exchange moves.  If we have conclusion 1 or 3
we can through a sequence of second flavor horizontal exchange moves
move $h_2$ to be below and consecutive with the horizontal position of $c^\prime$.
If we have conclusion 3 then through a sequence of the second flavor horizontal
exchange moves we can make the horizontal position of $h_2$ be minimal.
We then perform a first flavor horizontal exchange move to make
it maximal.  We follow this by a sequence of second flavor horizontal exchange moves
to place the horizontal position of $h_2$ below and consecutive with the
horizontal position of $c^\prime$.  This may place the horizontal
position of $h_2$ within the vertical support of a rectangle of $\cR$.
After $h_2$ is positioned near $c^\prime$ we perform the inverse of all of the vertical exchange moves so as
to place all of the vertical arcs back in their original angular position.  (Thus,
the isotopy of $h_2$ through $R$ is achieved by sequence of horizontal and vertical exchange moves.)
The new $\notchdiscprime$ will be such that $|\notchdisc \cap C_1| = |\notchdiscprime \cap C_1|$.\qed

%


Building on $\cO_{1.1}$, the remaining $\cO_1$ operations will
allow us to simplify the middle boundary of $\partial \notchdisc$.
In particular, for the description that follows we assume that no component of
$\notchdisc \setminus (C_1 \cap \notchdisc)$ is a whole or half disc.
The middle boundary will be in a component
$D^m \subset \notchdisc \setminus \{ D_{+1} \cup \cR\}$.
To list the possibilities, $D^m$ is either a disc or an annulus.
Moreover, $\partial D^m$ intersects either one, two or three
$\partial$-vertical arcs of $\partial \notchdisc$.  If $\partial D^m$
intersects one $\partial$-vertical arc then $D^m$ is an annulus
and $\notchdisc$ is a destabilizing disc.  If $D^m$ intersects two
$\partial$-vertical arcs and is an annulus then $\notchdisc$ is an exchange disc.
However, for two $\partial$-vertical arcs if $D^m$ is a disc then
$\notchdisc$ can be either an exchange or flyping disc.  If
$D^m$ intersects three $\partial$-vertical arcs then $\notchdisc$ is a flyping disc
and $D^m$ can be either an annulus or a disc.  Using the $2$-tuple subscript
$\{ \_ , t \} \in \{ \{ A , D \} , \{ 1,2,3\}\}$
(where $A = {\rm annulus}$, $D= {\rm disc}$ and
$t$ is number of vertical arcs), in the above order just described
we have that $D^m_{\_ , t}$ can be $D^m_{A,1} \subset \notchdiscd$,
$D^m_{A,2} \subset \notchdisce$, $D^m_{D,2} \subset \notchdisce \  {\rm or} \  \notchdiscf$,
and $D^m_{A,3}$ or $D^m_{D,3} \subset \notchdiscf$. 

\noindent
\underline{$\cO_{1.2}$--Simplifying $D^m_{A,1} \subset \notchdiscd$.}
To start we observe that the \scc's of $\partial D^m_{A,1}$ are $c_{max}$ and
$$\partial \notchdiscd = h^\eta_1 \bu \a^\eta_{\theta_1^{\ }} \bu
h^\eta_{2} \bu v^\eta_2 \bu \cdots \bu h^\eta_{l} \bu v^\eta_l \bu.$$
(The reader may wish the refer back to Lemma \ref{lemma:first simplification of notch disc}.)
Our $\cO_{1}$ operation assumption on the occurrence of whole and half discs
gives us that $c_{max} = int(\notchdiscd) \cap C_1$.

There are in fact two cases to argue since the arc $\a^\eta_{\theta_1^{\ }}$
can be either a front edge or a back edge.  For the argument below we will assume
that we have a back edge.  We will alert the reader to where the argument
differs for a front edge but leave the details to the reader.

The steps for simplifying $\partial \notchdiscd$ are as follows.

\noindent
\underline{Step-1:}  Through a sequence of vertical and horizontal exchange moves,
we reposition horizontal arcs $h^\eta_{3}, \cdots , h^\eta_{l}$ such that for
each their horizontal position is arbitrarily close to and below that of $c_{max}$.

This step is readily achieved by applying operation $\cO_{1.1}$ to each horizontal
arc in the middle boundary that is not adjacent to the $\partial$-vertical arc
$\a^\eta_{\theta_1^{\ }}$.  We will then have that for each horizontal arc,
$h^\eta_{3}, \cdots , h^\eta_{l}$, over its angular support it is consecutive with
$c_{max}$, i.e. arbitrarily close to and below that of $c_{max}$.

\noindent
\underline{Step-2:}  Through a sequence of horizontal simplifications we reduce
$\partial \notchdiscd$ to $$ h^\eta_1 \bu \a^\eta_{\theta_1^{\ }} \bu
h^\eta_{2} \bu v^\eta_2 \bu h^\eta_{3} \bu v^\eta_3 \bu$$
where the horizontal position of $h^\eta_{3}$ is
below and consecutive with that of $c_{max}$.

Having the repositioned $\partial \notchdisc$ in hand from Step-1, we now
realize that over the angular support of the edgepath $h^\eta_{l-1} \bu v^\eta_{l-1}
\bu h^\eta_l$ we have that the horizontal positions of $h^\eta_{l-1}$ and
$h^\eta_l$ are consecutive.  We can thus perform a horizontal simplification
reducing the number of vertical arcs of $\sheararcpres$ by one and shortening
the middle boundary.  We iterate this procedure until we have
$\partial \notchdiscd$ to $ h^\eta_1 \bu \a^\eta_{\theta_1^{\ }} \bu
h^\eta_{2} \bu v^\eta_2 \bu h^\eta_{3} \bu v^\eta_3 \bu$.
By construction the horizontal position of $h^\eta_{3}$ is
below and consecutive with that of $c_{max}$.

\noindent
\underline{Step-3:}  Through a sequence of vertical exchange moves followed by a shear
vertical simplification we reduce $\partial \notchdiscd$ to
$ h^\eta_{3} \bu \a^\eta_{\theta_1^{\ }} \bu h^\eta_{2} \bu v^\eta_2 \bu$ where 
the horizontal position of $h^\eta_{3}$ is below and consecutive with that of $c_{max}$.

We use the $\notchdisc$ coming out of Step-2 and giving us
$\partial \notchdisc = h^\eta_1 \bu \a^\eta_{\theta_1^{\ }} \bu
h^\eta_{2} \bu v^\eta_2 \bu h^\eta_{3} \bu v^\eta_3 \bu$.
In Step-1 \& -2 we did not use our assumption that $\a^\eta_{\theta_1^{\ }}$
is a back edge.  Here is where it comes into play.  In particular, this assumption
(along with having no whole or half discs) implies that our
horizontal positions for the horizontal arcs of $\partial \notchdisc$ has ordering
$z^{h_2}_{{\ }_{\ }} < z^{h_1}_{{\ }_{\ }} < z^{h_3}_{{\ }_{\ }} < z_{max}$.
(A back edge assumption would result in a
$z^{h_1}_{{\ }_{\ }} < z^{h_2}_{{\ }_{\ }} < z^{h_3}_{{\ }_{\ }} < z_{max}$ ordering.)
From this ordering we can conclude that there exists a rectangular subdisc
$R \subset \notchdisc \cap \cT_{\infty}$ having the following features.
\bi
\item[1.] $\partial R = v_3 \bu h_1 \bu \g^\partial \bu \g^h \bu$.
\item[2.] $\g^\partial \subset H_{\theta_1^{\ }}$, the disc fiber of $\fib$
containing $\a^\eta_{\theta_1^{\ }}$ and has vertical support
equal to that of $v_3$.
\item[3.] $\g^h$ has horizontal position equal to that of $h^\eta_3$ and angular
support equal to that of $h^\eta_1$.
\item[4.] From features 2 \& 3 we observe that the vertical support of $R$
is the same as that of $v^\eta_3$.
\ei
We now use $R$ to push $v^\eta_3$ forward and into the shearing interval
thus performing a shear vertical simplification.  Similar to the argument
given for operation $\cO_{1.1}$ the subdisc $R$ gives us that
vertical arcs in the angular support of $R$ will be nested with $v^\eta_3$.
This isotopy which is just pushing $v^\eta_3$ through the leaves of the induced
foliation on $R$ will correspond to performing a sequence of vertical
exchange moves followed by a shear vertical simplification.

For a front edge assumption a similar $R$ can be seen however it involves
$h^\eta_2$ and $v^\eta_2$ with $v^\eta_2$ being pushed backwards.\qed

\noindent
\underline{$\cO_{1.3}$--Simplifying $D^m_{A,2} \subset \notchdisce$.}
Here, we observe that $\partial D^m_{A,2}$ is a $c_{max}$ union
$$ \partial \notchdisce = h^\eta_1 \bu \a^\eta_{\theta_1^{\ }} \bu h^\eta_{2}
\bu v^\eta_2 \bu \cdots \bu  h^\eta_{l} \bu \a^\eta_{\theta_2^{\ }} \bu .$$
(See Lemma \ref{lemma:first simplification of notch disc}.)
Again, our $\cO_{1}$ operation assumption gives us that
$c_{max} = int(\notchdiscd) \cap C_1$.

As in $\cO_{1.2}$ there two cases: $\a^\eta_{\theta_1^{\ }}$ is a back edge
and $\a^\eta_{\theta_2^{\ }}$ is a front edge; or the edge assignments
are reversed.  For the argument below we will assume
that we have the former edge assignment.  We will alert the reader to where the argument
differs for the latter edge assignment but leave the details to the reader.

The steps for simplifying $\partial \notchdisce$ are as follows.

\noindent
\underline{Step-1:}  Through a sequence of vertical and horizontal exchange moves,
we reposition horizontal arcs $h^\eta_{3}, \cdots , h^\eta_{l-1}$ such that for
each one we have that over its angular support its horizontal position is minimal.

The step is achieved in essentially the same manner as Step-1 of $\cO_{1.2}$
by utilizing $\cO_{1.1}$.  However, this will repositioned all our
horizontal arcs near $c_{max}$.  (If we were arguing the latter edge assignment,
we would be interested in achieving this repositioning and would be done.)
To achieve minimal horizontal position over the angular support of one of our
horizontal arcs we perform a first flavor horizontal exchange move.
After a sequence of such first flavor moves we achieve the goal of this step.

If it happens that our initial middle boundary has exactly two horizontal
arcs then we skip this step and Step-2, and proceed to Step-3.

\noindent
\underline{Step-2:}  Through a sequence of horizontal simplifications we reduce
$\partial \notchdisce$ to 
$$h^\eta_1 \bu \a^\eta_{\theta_1^{\ }} \bu h^\eta_{2}
\bu v^\eta_2 \bu  h^\eta_{3} \bu v^\eta_3 \bu h^\eta_4 \bu \a^\eta_{\theta_2^{\ }} \bu$$
where $h^\eta_{3}$ has minimal horizontal position.

This step is a repeat of the argument of Step-2 of $\cO_{1.2}$ except that all of
our horizontal arcs are at a minimal (not maximal) horizontal position in their angular support.

For the latter edge assignment case we would end up with $h^\eta_3$ being at
maximal horizontal position and below $c_{max}$.

\noindent
\underline{Step-3:}  Through a sequence of vertical exchange moves and shear
vertical simplifications we reduce $\partial \notchdisce$ to
$h^\eta_1 \bu \a^\eta_{\theta_1^{\ }} \bu h^\eta_{2}
\bu \a^\eta_{\theta_2^{\ }} \bu$.

We again adapt operation $\cO_{1.1}$ to the resulting Step-2 middle boundary.
Specifically, we consider the edgepath
$\a^\eta_{\theta_1^{\ }} \bu h^\eta_{2} \bu v^\eta_2$.  Observe that
the vertical support of $\a^\eta_{\theta_1^{\ }}$ is above the horizontal position
of $h^\eta_{2}$, whereas the vertical support of $v^\eta_2$ is below the position
of $h^\eta_{2}$.  So we can push any vertical arc having vertical support below
$h^\eta_{2}$ and angular position in the angular support of $h^\eta_{2}$
backwards and past $\a^\eta_{\theta_1^{\ }}$.  But, pushing them past
$\a^\eta_{\theta_1^{\ }}$ is pushing them into a shearing interval.  So through a
sequence of vertical exchange moves and shear vertical simplification we can
alter $\sheararcpres$ so as to assume that there are no vertical arcs below
the horizontal position of $h^\eta_2$ in the angular support of $h^\eta_2$.

Once this is achieve we consider the {\em virtual rectangular disc}
$R = [\theta_1^{h^\eta_2}, \theta_2^{h^\eta_2}] \times v^\eta_2 \subset \mathbb R^3$, where
the angular interval is the angular support of $h^\eta_2$.
Since we have no vertical arcs of $\sheararcpres$ below and in the angular support
of $h^\eta_2$ we have that $int(R) \cap \notchdisc = \emptyset$.  We can then
utilize our virtual $R$ in a similar manner to our $R$-subdisc in Step-3 of $\cO_{1.2}$
to push $v^\eta_2$ backwards into the shearing interval associated with
$\a^\eta_{\theta_1^{\ }}$.
Arguing in a similar fashion but with a forward push of $v^\eta_3$ past
$\a^\eta_{\theta_2^{\ }}$ we can achieve the $\partial \notchdisc$ described in
this step.

If our initial middle boundary had only two horizontal arcs then we can still
construct a virtual rectangular disc and we will push the single vertical
arc in the direction of the shortest (in the vertical support sense)
$\partial$-vertical arc.

For the latter edge assignment case we can use $R$-subdiscs similar to the
one in Step-3 of $\cO_{1.2}$.\qed

\noindent
\underline{$\cO_{1.4}$--Simplifying $D^m_{D,2} \subset \notchdisce \  {\rm or} \  \notchdiscf$.}
Here we have $\partial D^m_{D,2}$ being of the form
$$ \g_+ \bu \g^\partial_1 \bu h^\eta_2 \bu v^\eta_1 \bu h^\eta_{2}
\bu v^\eta_2 \bu \cdots \bu  h^\eta_{l} \bu \g^\partial_2 \bu $$
were $\g_+$ is a horizontal boundary of a rectangle of $\cR$,
and $\g^\partial_1$ \& $\g^\partial_2$ are subarcs in
distinct $\partial$-vertical arcs of $\notchdisc$.
Again, our $\cO_{1}$ operation assumption gives us that
$int(D^m_{D,2}) \cap C_1 = \emptyset$.
The steps for simplifying $\partial \notchdisc$ are as follows.

By our assumptions in our discussion after Lemma
\ref{lemma:first simplification of notch disc} on the start/end of the middle boundary we
know that for $\notchdisce$ we must have $\g_+$ being a type-2 \sa and
a type-3 \sa for $\notchdiscf$.  For the $\notchdiscf$ case
there are two cases since we have either a positive or negative flype.
We will assume a positive flype and leave the similar negative case to the
reader.

For the case of $D^m_{D,2} \subset \notchdisce$ we appeal to our $\cO_{1.3}$
argument as a model.  As just observed, we have that $\g^\partial_1$ is a subarc
on a back edge and $\g^\partial_2$ is on a front edge.  Then imitating
Step-1 of $\cO_{1.3}$ we can use operation $\cO_{1.1}$ to reposition
all of the horizontal edges of $D^m_{D,2} \cap \sheararcpres$ to be
near $\g_+$.  Continuing and imitating Step-2 of $\cO_{1.3}$ we can
perform horizontal simplifications to alter $\partial D^m_{D,2}$ to
be $\g_+ \bu \g^\partial_1 \bu h^\eta_2 \bu v^\eta_2 \bu h^\eta_{3}
\bu v^\eta_3 \bu  h^\eta_{4} \bu v^\partial_2 \bu$.
Finally, we can imitate Step-3 of $\cO_{1.3}$ altering
$\partial D^m_{D,2}$ to
be $\g_+ \bu \g^\partial_1 \bu h^\eta_2 \bu \g^\partial_2 \bu$
through a sequence of vertical exchange moves and vertical
simplifications.

We also appeal to $\cO_{1.3}$ in the $D^m_{D,2} \subset \notchdiscf$
positive flype case with $\g^\partial_1$ \& $\g^\partial_2$ both being subarcs
on differing front edges.  Again, imitating
Step-1 of $\cO_{1.3}$ we can use operation $\cO_{1.1}$ to reposition
all of the horizontal edges of $D^m_{D,2} \cap \sheararcpres$ to be
near $\g_+$. Continuing and imitating Step-2 of $\cO_{1.3}$ we can
perform horizontal simplifications to alter $\partial D^m_{D,2}$ to
be $\g_+ \bu \g^\partial_1 \bu h^\eta_2 \bu v^\eta_2 \bu h^\eta_{3}
\bu v^\eta_3 \bu  h^\eta_{4} \bu v^\partial_2 \bu$.
Finally, we can imitate Step-3 of $\cO_{1.3}$ altering
$\partial D^m_{D,2}$ to
be $\g_+ \bu \g^\partial_1 \bu h^\eta_2 \bu \g^\partial_2 \bu$
through a sequence of vertical exchange moves and vertical
simplifications.  However, this totally simplifies the middle
boundary to one horizontal arc.  Thus, the resulting notch disc
cannot imply the existence of a flype.  So we could not have had
$D_0$ be contained in a flyping disc.\qed

\noindent
\underline{$\cO_{1.5}$--Simplifying $D^m_{A,3} \subset \notchdiscf$.}
Briefly, by Lemma \ref{lemma:first simplification of notch disc} statement d.
we have $\partial D^m_{D,3}$ is the union of $c_{max}$ and
$$h^\eta_1 \bu \a^\eta_{\theta_1^{\ }} \bu h^\eta_{2}
\bu \a^\eta_{\theta_2^{\ }} \bu h^\eta_{3} \bu v^\eta_3 \bu
\cdots \bu h^\eta_{l} \bu \a^\eta_{\theta_3^{\ }} \bu$$
Again, the $\cO_{1}$ operation assumption gives us that $c_{max} = int(\notchdiscd) \cap C_1$.

Our edge assignment assumption for $\notchdiscf$ (see FA-a) gives us
two similar but different cases to argue.  Either $\a^\eta_{\theta_1^{\ }}$,
$\a^\eta_{\theta_2^{\ }}$ and $\a^\eta_{\theta_3^{\ }}$
are a front edge, back edge, front edge, respectively; or,
they are a back edge, front edge and back edge, respectively.
We will argument the former case (the positive flype case)
and, with hints, leave the latter case to the reader.

The steps for simplifying $\partial \notchdiscd$ imitates the steps
for $\cO_{1.2}$.  Specifically, starting with the middle boundary
$ h^\eta_{3} \bu v^\eta_3 \bu
\cdots \bu h^\eta_{l} $ by a sequence vertical and horizontal
exchange moves we can reposition horizontal arcs
$h^\eta_4$ through $h^\eta_{l-1}$ so that each has minimal horizontal
position over its angular support.  Next, by a sequence of horizontal
simplifications we can shorten the middle boundary until it
corresponds to $ h^\eta_{3} \bu v^\eta_3 \bu
h^\eta_4 \bu v^\eta_4 \bu h^\eta_{5}$ with the horizontal position of
$h^\eta_4$ being minimal.  Finally, using our back edge assumption
for $\a^\eta_{\theta_2^{\ }}$ and $\a^\eta_{\theta_3^{\ }}$
we realize that we can imitate the argument in Step-3 of $\cO_{1.3}$ to
construct a virtual rectangular disc for pushing $v^\eta_3$ forward
and past $\a^\eta_{\theta_1^{\ }}$.  As in Step-3 of $\cO_{1.3}$, this will be achieved through a
sequence of vertical exchange moves and shear vertical simplifications.
We then obtain $ h^\eta_{3} \bu v^\eta_3 \bu
h^\eta_4 $ as the middle boundary.

In the front-front edge assignment case we would construct a virtual rectangular disc for
pushing $v^\eta_4$ backward past $\a^\eta_{\theta_3^{\ }}$.\qed

\noindent
\underline{$\cO_{1.6}$--Simplifying $D^m_{D,3} \subset \notchdiscf$.}
This remaining simplification is handled exactly as the previous simplification
in $\cO_{1.5}$.  In middle boundary of $\partial \notchdiscf$ contained
in $\partial D^m_{D,3}$ becomes $ h^\eta_{3} \bu v^\eta_3 \bu h^\eta_4 $
through a sequence of our moves from \S\ref{subsection:rectangular diagrams}.
It is useful to observe that $\g_+$ must be either a type-2 or -4 \sa.\qed

\noindent
{\bf $\cO_2$--Elimination of whole discs.}
There are three differing configurations of the tiled foliation which can be used to
detect the occurrence of horizontal and vertical exchange moves for the
elimination of whole discs.

\begin{figure}[htb]
\centerline{\includegraphics[scale=.6, bb=0 0 630 564]{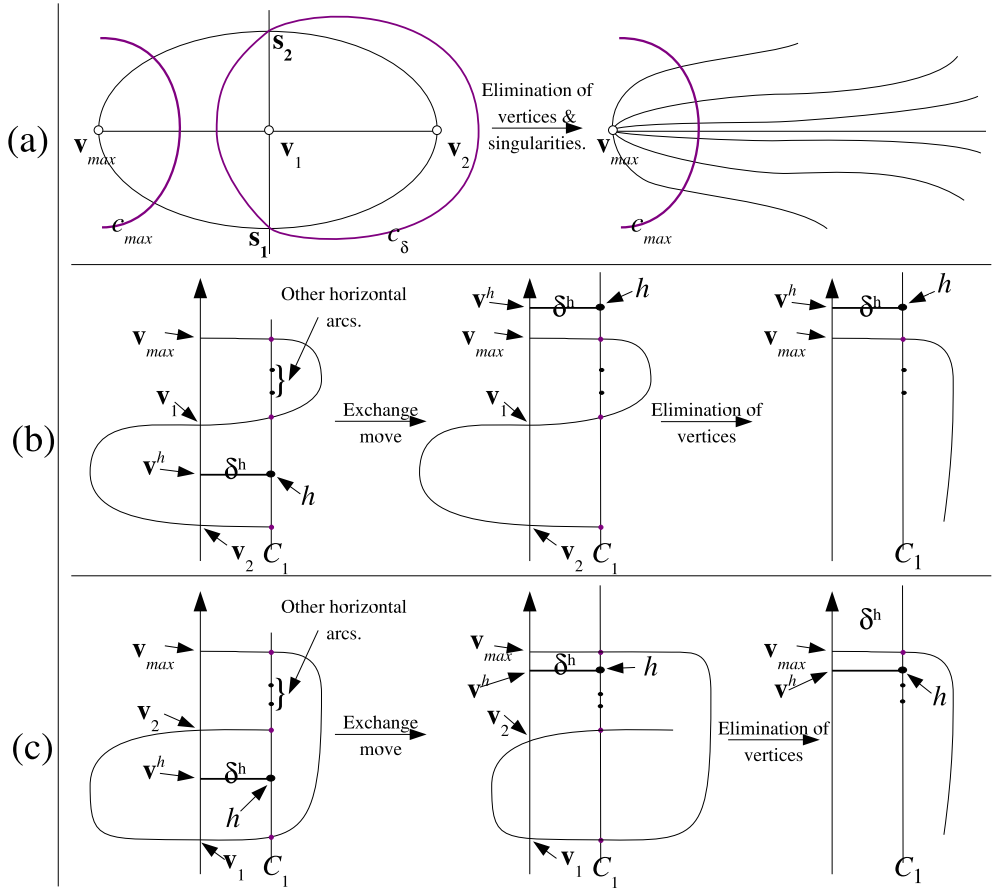}}
\caption{In (a) we depict the tiled foliation in a neighborhood of a valence
two vertex and show the change to the foliation when we eliminate
positive/negative pairs of vertices and singularities.  In (b) we illustrate
a cross-section of the neighborhood in (a) along with the cross-section of the
coning disc $\delta^h$ and the cylinder $C_1$ in the case where the
order of vertices on $\axis$ is $ \ver_2 < \ver^h < \ver_1 < \ver_{max}$.
In (c) we have a similar cross-section for the order
$\ver_1 < \ver^h < \ver_2 < \ver_{max}$.}
\label{figure:valence two}
\end{figure}

\noindent
\underline{$\cO_{2.1}$--Elimination of a whole disc near $c_{max}$.}  We consider
a coning disc $\delta^h$ having an intersection arc $\g \subset \delta^h \cap \notchdisc$
which is {\em innermost} on the component of $\notchdisc \setminus \cR$
containing $D_{+1}$.  That is, for every leaf $a_\theta$ in the radial
foliation of $\notchdisc$ intersecting $\g$ we have that when traversing
$a_\theta$ with $\ver_{max}$ as the starting point, the first time $a_\theta$ intersects
any coning disc is at the point $a_\theta \cap \g $.
For such a choice of coning disc we apply the transition
$ (\sheararcpres,\notchdisc) \Ttran (\sheararcpres,\tilenotchdisc)$
to produce our tiled notch disc.
Referring back to Figure \ref{figure:tiles}(a), if $\g$ is innermost then all of the
``top of the page'' leaves adjacent to the negative vertex are ``coned'' to $\ver_{max}$.
Moreover, the curve $c_\g$ will be ``nearest'' to $c_{max}$.
Figure \ref{figure:valence two}(a) illustrates these
features.  The illustrated foliated subdisc of $\tilenotchdisc$
contains $\ver_{max}$ with two additional vertices, $\ver_1$ and $\ver_2$, and the two
singularities, $\sin_1$ and $\sin_2$.
$\ver_{max}$ and $\ver_2$ are positive vertices, whereas $\ver_1$ has negative parity.
The \scc $c_\gamma \subset \tilenotchdisc \cap C_1$
passes through the two singularities and
encircles $\ver_1$ and $\ver_2$.  Between $\ver_1$ and $\ver_2$ on the axis $\axis$ there is a vertex, $\ver^h$, of coning disc
$\delta^h$ for our horizontal arc $h \subset \sheararcpres$.
And, the $b$-arcs adjacent to $\ver_{max}$ \&
$\ver_1$ intersect $\tilenotchdisc \cap C_1$ only at $c_{max}$ and $c_\g$.

We now perform the classical exchange move of \cite{[BM1], [BM2], [BM3], [BM4], [BF]} that
isotop $\delta^h$ so that the resulting $\ver^h$ and $\ver_{max}$
are at consecutive horizontal positions.  After the isotopy
$|\delta^h \cap int(\tilenotchdisc)|$ will decrease by at least one:  $\g$ will be eliminated; and,
all the arcs of
$\delta^h \cap int(\tilenotchdisc)$ in $\delta^h$ between $\ver^h$ and $\g$ will also be
eliminated.  The \scc's associated with these arcs will also be eliminated.
However, any arcs of $\delta^h \cap int(\tilenotchdisc)$ in $\delta^h$ between $\g$ and $h$
will remain along with their associated \scc. 
Thus, the intersection curve $c_\g$ is
eliminated to reduce $|int(\tilenotchdisc) \cap C_1|$, and
we replace $\tilenotchdisc$ with a radially foliated notch disc $\notchdiscprime$.

We have two cases depending on the initial order of our four vertices on $\axis$:
when the initial order is $\ver_2 < \ver^h < \ver_1 < \ver_{max}$ the classical exchange move
will isotop $\delta^h$ so the $\ver_{max} < \ver^h$
(Figure \ref{figure:valence two}(b)-middle); and,
when the initial order is $\ver_1 < \ver^h < \ver_2 < \ver_{max}$ the classical exchange move
will isotop $\delta^h$ so the $\ver^h < \ver_{max}$ (Figure \ref{figure:valence two}(c)-middle).
(Treating vertices as their horizontal-position's real values on $\axis$ is a
convenient abuse of notation.)
In the Figure \ref{figure:valence two}(b) case $\sin_2$ is necessarily a negative singularity
and $\sin_1$ is positive.  Whereas in the Figure \ref{figure:valence two}(c)
case $\sin_2$ is positive and $\sin_1$ is negative.  
Our difficulty is in seeing that these two isotopy cases in fact correspond
to a sequence of horizontal and vertical exchange moves
as defined in \S \ref{subsection:rectangular diagrams}.
In particular, for the second flavor of horizontal exchange moves the key features of consecutive and nested have to be verified.  This verification will be similar but different for
the Figure \ref{figure:valence two}(b) and (c) cases.

For convenience, let $\sin_1$ and $\sin_2$
have angular position $\theta_1$ and $\theta_2$ respectively.
Then the $b$-arcs adjacent
to $\ver_{max}$ and $\ver_1$ occur over the angular interval $[\theta_1 , \theta_2] \subset S^1$.
Moreover, the angular support of $h$, $[\theta^h_1 , \theta^h_2]$, is properly contained in but arbitrarily close to
$[\theta_1 , \theta_2]$.

We start with the case in Figure \ref{figure:valence two}(b).  If each
horizontal arc of $\sheararcpres$ having horizontal position in the
interval $(\ver_1 , \ver_{max}) \subset \axis$ also has angular support containing the interval
$[\theta^h_1 , \theta^h_2]$ then the nested feature is satisfied and
the isotopy in (b) is just a sequence of horizontal
exchange moves (second flavor).

Now, suppose that $h^\prime \subset \sheararcpres$ is a horizontal arc having horizontal
position in the interval $(\ver_1 , \ver_{max})$ and angular support
$[\theta_{1^{\ }}^{h^\prime} , \theta_{2^{\ }}^{h^\prime}]$ with
$\theta_{2^{\ }}^{h^\prime} \in (\theta^h_1 , \theta^h_2)$.  Without loss of generality
we can assume that the angle $\theta_{2^{\ }}^{h^\prime}$ is the nearest such angle
to $\theta^h_2$, and let $v^\prime \subset \sheararcpres \cap
H_{\theta_{2}^{h^\prime}}^{{\ }_{\ }}$ be the vertical arc
adjacent to $h^\prime$.  Since the $b$-arcs adjacent to $\ver_1$ and $\ver_{max}$ have
vertical support $[\ver_1 , \ver_{max}] \subset \axis$, the vertical
support of $v^\prime$ must be contained in $(\ver_1 , \ver_{max})$, i.e.
$v^\prime \cap int(\tilenotchdisc) = \emptyset$.  Moreover,
by our `$\theta_{2^{\ }}^{h^\prime}$ nearest $\theta^h_2$' assumption,
all other vertical arcs of $\sheararcpres$ having angular position inside the interval
$(\theta_{2^{\ }}^{h^\prime} , \theta^h_2)$ must have vertical support
strictly below $\ver_1$.  Thus, the key feature of nested for
vertical exchange moves is satisfied for $v^\prime$ and the next consecutive vertical
arc in the angular interval $(\theta_{2^{\ }}^{h^\prime} , \theta^h_2)$.
Moving $v^\prime$ past this arc and iterating this procedure, through a sequence of
vertical exchange moves we can push $v^\prime$ forward and outside the interval
$[\theta^h_1 , \theta^h_2]$.  Thus, through a sequence of vertical exchange moves
we can assume that the angular support of all horizontal arcs having horizontal
position in $(\ver_1 , \ver_{max})$ properly contains $[\theta^h_1 , \theta^h_2]$.
The isotopy of Figure \ref{figure:valence two}(a) can thus be achieved through a sequence of vertical
and horizontal exchange moves.  Moreover, once we have placed
$h$ above that of $\ver_{max}$ we perform the inverse of all of the vertical exchange moves
to place all of the vertical arcs back in their original angular position.

But, we are not quite done since the isotopy in Figure \ref{figure:valence two}(a)
leaves the horizontal position of $h$ above that of $\ver_{max}$.  To resolve this
remaining issue we note that the horizontal position of $h$ is now maximal and we
can perform a first flavor horizontal exchange move on $h$ making its
horizontal position minimal.  Notice this will not change the value of
$|\notchdisc \cap C_1|$.  Moreover, since a first flavor horizontal exchange move
on $h$ corresponds to re-choosing the point a infinity on $\axis$ (a choice between the
horizontal positions of $\ver_{max}$ and $h$), the foliation of $\notchdisc$ and its new intersection
with $C_1$ remains the same after this exchange move.

We next deal with the case in Figure \ref{figure:valence two}(c).
It is argued in a similar fashion to the previous case except we work with
arcs of $\sheararcpres$ whose horizontal position and vertical support are
below that of $\ver_1$.  Specifically, if each
horizontal arc of $\sheararcpres$ having horizontal position in the
interval $(-\infty , \ver_1) \subset \axis$ also has angular support containing the interval
$[\theta^h_1 , \theta^h_2]$ then the nested feature is satisfied and
the isotopy in (b) is just a sequence of second flavor horizontal
exchange moves followed by a single first flavor horizontal exchange move.

Similar to before,
suppose that $h^\prime \subset \sheararcpres$ is a horizontal arc having horizontal
position in the interval $(-\infty , \ver_1)$ and angular support
$[\theta_{1^{\ }}^{h^\prime} , \theta_{2^{\ }}^{h^\prime}]$ with
$\theta_{2^{\ }}^{h^\prime} \in (\theta^h_1 , \theta^h_2)$.  Again, without loss of generality
we can assume that the angle $\theta_{2^{\ }}^{h^\prime}$ is the nearest such angle
to $\theta^h_2$, and let $v^\prime \subset \sheararcpres \cap
H_{\theta_{2}^{h^\prime}}^{{\ }_{\ }}$ be the vertical arc
adjacent to $h^\prime$.  Since the $b$-arcs adjacent to $\ver_1$ and $\ver_{max}$ have
vertical support $[\ver_1 , \ver_{max}] \subset \axis$, the vertical
support of $v^\prime$ must be contained in $(-\infty , \ver_1)$, i.e.
$v^\prime \cap int(\tilenotchdisc) = \emptyset$.  As before,
all other vertical arcs of $\sheararcpres$ having angular position inside the interval
$(\theta_{2^{\ }}^{h^\prime} , \theta^h_2)$ must have vertical support
strictly above $\ver_1$.  Thus, the key feature of nested for
vertical exchange moves is satisfied for $v^\prime$ and the next consecutive vertical
arc in the angular interval $(\theta_{2^{\ }}^{h^\prime} , \theta^h_2)$.
Moving $v^\prime$ past this arc and iterating this procedure, through a sequence of
vertical exchange moves we can push $v^\prime$ forward and outside the interval
$[\theta^h_1 , \theta^h_2]$.  Thus, through a sequence of vertical exchange moves
we can assume that the angular support of all horizontal arcs having horizontal
position in $(-\infty , \ver_1)$ properly contains $[\theta^h_1 , \theta^h_2]$.
Once we have placed
$h$ below that of $\ver_{max}$, we perform the inverse of all of the vertical exchange moves
to place all of the vertical arcs back in their original angular position.
The isotopy of Figure \ref{figure:valence two}(a) can be achieved through a sequence of vertical
and horizontal exchange moves.\qed

\begin{figure}[htb]
\centerline{\includegraphics[scale=.7, bb=0 0 469 380]{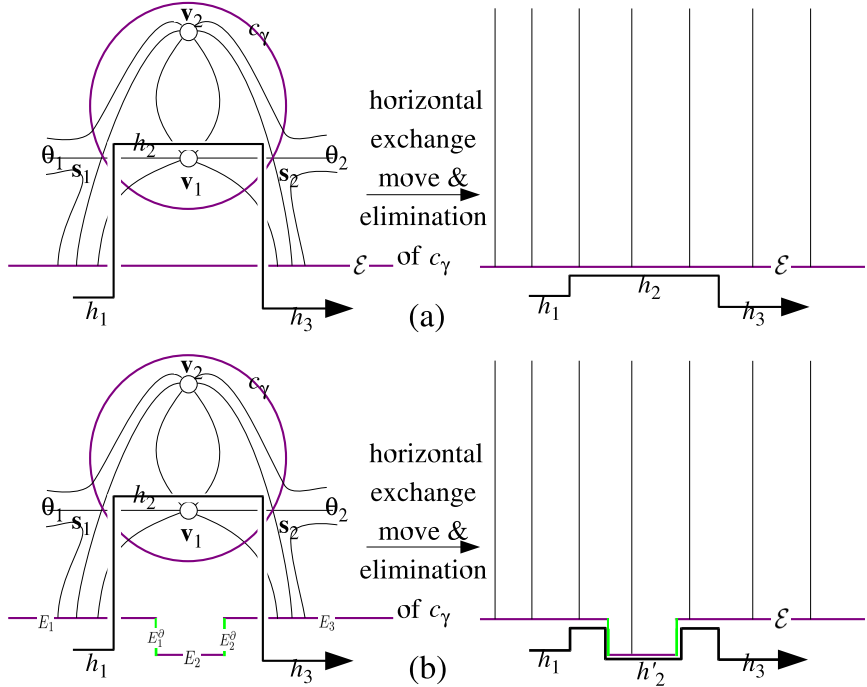}}
\caption{The top sequence illustrates the use of a horizontal exchange move in
the elimination of a valence two vertex and a $c_\gamma$ curve that is
nearest $\cE$.  The bottom sequence illustrates
shear horizontal exchange moves when $cE$ contains $\partial$-vertical arcs (possibility-B).}
\label{figure:boundary exchange moves large}
\end{figure}

\noindent
\underline{$\cO_{2.2}$--Elimination of a whole disc near an edgepath
$\cE$ that is away from $\partial$-vertical arcs.}
We consider
a coning disc $\delta^h$ having an intersection arc $\g \subset \delta^h \cap \notchdisc$
which is {\em nearest} on a $D_0$ or the $D_{+1}$ component, $D^\star \subset \notchdisc \setminus int(\cR)$,
to an edge-path $\cE \subset \partial D^\star$.
That is, for every leaf $a_\theta$ in the radial
foliation of $\notchdisc$ intersecting $\g$ we have that it also intersects $\cE$ and, when traversing
$a_\theta$ from $a_\theta \cap \cE$ to $a_\theta \cap \g$, the first time $a_\theta$ intersects
any coning disc is at $a_\theta \cap \g $.  (Picture $\cE$
having slightly larger angular support than that of $c_\g \subset \notchdisc \cap C_1$,
the \scc associated with $\g$.)
Finally, assume that no leaf $a_\theta$ in the
radial foliation of $\notchdisc$ intersecting $\g$ also contains a $\partial$-vertical arc
of $\notchdisc$.

Such an edgepath $\cE$ is a portion of either $\sheararcpres \cap D^\star$ or
a horizontal boundary arc of a component of $\cR$.  Taking the most straight forward situations
first, assume that $\cE$ is a portion of either a horizontal arc of $\sheararcpres$ or
a horizontal boundary arc of a component of $\cR$.

For such a choice of coning disc we apply the transition
$ (\sheararcpres,\notchdisc) \Ttran (\sheararcpres,\tilenotchdisc)$
to produce our tiled notch disc.
This choice of coning disc yields a configuration similar to that
in $\cO_{2.1}$.
Referring back to Figure \ref{figure:tiles}(a)-left, if $\g$ is nearest $\cE$
then the curve $c_\g$ will be ``nearest'' to $\cE \subset \partial D^\star$.
Figure \ref{figure:boundary exchange moves large}(a) illustrates these
features and we now appeal to this figure's labeling.
The illustrated foliated subdisc of $\tilenotchdisc$
has horizontal boundary arc $\cE$ along with two additional
vertices, $\ver_1$ and $\ver_2$, and the two singularities, $\sin_1$ and $\sin_2$
in $int(\tilenotchdisc)$.
$\ver_2$ and $\ver_1$ have opposite parity.
(When $\cE$ is a horizontal arc of $\sheararcpres$ then $\ver_1$ is positive and $\ver_2$ is negative.)
The \scc $c_\gamma \subset \tilenotchdisc \cap C_1$
passes through the two singularities and
encircles $\ver_1$ and $\ver_2$.  Between $\ver_1$ and $\ver_2$ on the axis $\axis$ there is a vertex,
$\ver^{h_2}_{\ }$, of coning disc
$\delta^{h_2}_{\ }$ for our horizontal arc $h_2 \subset \sheararcpres$.
And, the $a$-arcs in the induced foliation of $D^\star$ adjacent to
$\ver_1$ and having endpoints on $\cE$
intersect $\tilenotchdisc \cap C_1$ only at $\cE$ \& $c_\g$.

We perform the classical boundary exchange move of
\cite{[BM1], [BM2], [BM3], [BM4], [BF]} that is predicated on the existence
of a valence two vertex near the boundary of a surface.
Figure \ref{figure:boundary exchange moves}(a) illustrates the isotopy
and elimination of the \scc $c_\g$.
Specifically, we
isotop $\delta^{h_2}_{\ }$ so that the resulting $\ver^{h_2}_{\ }$ and $\cE$
are at consecutive horizontal positions---we allow for the resulting horizontal position
of $\ver^{h_2}_{\ }$ to possibly be in the vertical support of a rectangle of $\cR$
when $\cE$ is a horizontal boundary arc of a component of $\cR$.
After the isotopy
$|\delta^{h_2}_{\ } \cap (int(\tilenotchdisc) \setminus \cR)|$ has been reduced by
at least one:  the arc $\g$ has been eliminated; and, any arc of
$\delta^{h_2}_{\ } \cap (int(\tilenotchdisc) \setminus \cR)$ between $\ver^{h_2}_{\ }$ and
$\g$ has been eliminated.
Thus, the intersection curve $c_\g$ is
eliminated reducing $|int(\tilenotchdisc) \cap C_1|$, and
we replace $\tilenotchdisc$ with a radially foliated notch disc $\notchdiscprime$.
(Again, any \scc of $\tilenotchdisc \cap C_1$ 
associated with the previously mentioned arcs between $\ver^{h_2}_{\ }$ and
$\g$ is eliminated
by this classical exchange move.)

The details of this argument are exactly like those for $\cO_{2.1}$ when we substitute
the horizontal position of $\cE$ for the horizontal position of $\ver_{max}$.
In fact, the cross-sections illustration in Figure \ref{figure:valence two}(b) \& (c)
can be appealed to if we alter our understanding of the illustration by thinking
of the point $\ver_{max} \subset \axis$ as being the cone point for a coning disc
of $\cE$.  The isotopies of Figure \ref{figure:valence two}(b) \& (c)
may possibly result in the horizontal position of $h_2$ being in the vertical
support of a rectangle of $\cR$.
As such we will not be repetitive by imitating the argument of $\cO_{2.1}$, but conclude that
through a sequence of vertical and horizontal (both flavors) exchange moves
we can eliminate the $c_\g$ intersect curve and replace $\notchdisc$ with
a new $\notchdiscprime$ having $|\notchdisc \cap C_1| > |\notchdiscprime \cap C_1|$.
(We do not forget to perform the inverse of all of the vertical exchange moves so as
to place all of the vertical arcs back in their original angular position.)

The remaining case is when $\cE$ is in $\sheararcpres$ but is not a portion of a single
horizontal arc.  However, we can transform $\cE$ into a portion of such a horizontal arc
by utilizing $\cO_{1.1}$ and imitating the argument in Step-1 \& Step-2 of $\cO_{1.2}$.
Such a transformation of $\cE$ will necessarily employ the use of horizontal and vertical
exchange moves, and horizontal simplifications.
\qed

\noindent
\underline{$\cO_{2.3}$--Elimination of a whole disc near an edgepath
$\cE$ containing $\partial$-vertical arcs.}
Again, we consider
a coning disc $\delta^h$ having an intersection arc $\g \subset \delta^h \cap \notchdisc$
which is {\em nearest} on a $D_0$ or the $D_{+1}$ component, $D^\star \subset \notchdisc \setminus int(\cR)$,
to an edge-path $\cE \subset \partial D^\star$.
As before this means that for every leaf $a_\theta$ in the radial
foliation of $\notchdisc$ intersecting $\g$ we have that it also intersects $\cE$ and, when traversing
$a_\theta$ from $a_\theta \cap \cE$ to $a_\theta \cap \g$, the first time $a_\theta$ intersects
any coning disc is at $a_\theta \cap \g $.
However, we now assume that there are leaves $a_\theta$ in the
radial foliation of $\notchdisc$ intersecting $\g$ containing $\partial$-vertical arcs
of $\notchdisc$.

Over the angular support of $\g$ we list the initial possibilities: (i) $\notchdiscd$ 
with one such leaf containing a $\partial$-vertical arc; (ii) $\notchdisce$ with one or two such leaves
containing $\partial$-vertical arcs; and, (iii)
$\notchdiscf$ with one, two or three such leaves containing $\partial$-vertical arcs.
From this list of six initial possibilities it is useful to eliminate
outright three of them.  In particular, for possibilities (i), (ii) with two leaves, and
(iii) with three leaves, the \scc $c_\g$ will be {\em close} to $c_{max}$.  That is, when traversing
a leaf $a_\theta$ in the radial foliation of $\notchdisc$ from $a_\theta \cap c_{max}$ to $a_\theta \cap c_\g$,
$a_\theta$ will only intersect other \scc's of $\notchdisc \cap C_1$, i.e. it does not intersect any \sa.
Clearly \scc's that are close to $c_{max}$ can be eliminated through a sequence of $\cO_{2.1}$ starting
with the innermost first.  Thus, we have to eliminate only (ii) with one leaf (call this possibility-A) and
(iii) with two leaves (call this possibility-B) containing $\partial$-vertical arcs.

To continue, we consider the simplest situations first.  For possibility-A, let
$\cE = E_1 \bu E^\partial \bu E_2$ where $E^\partial$ is a portion of a $\partial$-vertical arc and
$E_1$ \& $E_2$ are portions of horizontal arcs of either $\sheararcpres$ or $\partial \cR$.
For possibility-B, let $\cE = E_1 \bu E^\partial_1 \bu E_2 \bu E^\partial_2 \bu E_3$ where
$E^\partial_1$ \& $E^\partial_2$ is a portions of a $\partial$-vertical arcs and
$E_1$, $E_2$ \& $E_3$ are portions of horizontal arcs of either $\sheararcpres$ or $\partial \cR$.
(Figure \ref{figure:boundary exchange moves large}(b)-left illustrates this simple situation for possibility-B.)

Now, we return to our choice of coning disc and we apply the transition
$ (\sheararcpres,\notchdisc) \Ttran (\sheararcpres,\tilenotchdisc)$
to produce our tiled notch disc.
Referring back to Figure \ref{figure:tiles}(b)-left, if $\g$ is nearest $\cE$
then the curve $c_\g$ will be ``nearest'' to $\cE \subset \partial D^\star$.
Figure \ref{figure:boundary exchange moves large}(b) illustrates these
features and we now appeal to this figure's labeling.
The illustrated foliated subdisc of $\tilenotchdisc$
has edgepath $\cE(= E_1 \bu E^\partial_1 \bu E_2 \bu E^\partial_2 \bu E_3)$ for the simple-possibility-B along with two additional
vertices, $\ver_1$ and $\ver_2$, and the two singularities, $\sin_1$ and $\sin_2$
in $int(\tilenotchdisc)$.
$\ver_2$ and $\ver_1$ have opposite parity.
(When $\cE$ contains horizontal arcs of $\sheararcpres$ then $\ver_1$ is positive and $\ver_2$ is negative.)
The \scc $c_\gamma \subset \tilenotchdisc \cap C_1$
passes through the two singularities and
encircles $\ver_1$ and $\ver_2$.  Between $\ver_1$ and $\ver_2$ on the axis $\axis$ there is a vertex,
$\ver^{h_2}_{\ }$, of coning disc
$\delta^{h_2}_{\ }$ for our horizontal arc $h_2 \subset \sheararcpres$.
And, the $a$-arcs in the induced foliation of $D^\star$ adjacent to
$\ver_1$ and having endpoints on $\cE$
intersect $\tilenotchdisc \cap C_1$ only at $\cE$ \& $c_\g$.

We cannot quite perform the classical boundary exchange that is predicated on the existence
of a valence two vertex near the boundary of a surface since the $\partial$-vertical arcs obstruct
such a move.  However, for each portion of a $\partial$-vertical arc in our edgepath $\cE$ we
can perform a shear horizontal exchange move.
Figure \ref{figure:boundary exchange moves}(b) illustrates the isotopy
and elimination of the \scc $c_\g$.  This move can be thought of as ``shearing'' the coning discs
into two (for possibility-A) or three (for possibility-B) coning discs and then 
isotoping them so that the resulting sheared-coning discs and horizontal sub-arcs of $\cE$
are at consecutive horizontal positions.  (See Figure \ref{figure:boundary exchange moves}(b)-right.)
After the isotopy
$c_\g$ has been eliminated and for the resulting $\notchdisc$ we have that
$|[ \notchdisc \setminus \cR ] \cap [ C_1 \setminus (C_1 \cap \cI) ]|$
has been reduced.

Once we have sheared $\delta^h$ the details establishing that this entire operation
eliminating $c_\g$ corresponds to a sequence of our elementary moves are
exactly like the details in $\cO_{2.2}$ expect that we are applying them to
the newly formed horizontal arcs coming from shearing $h$.
The remaining case is when the $E_i$ portions of $\cE$ are in $\sheararcpres$ but is not a portion of a single
horizontal arc.  However, we can transform each such $E_i$ into a portion of such a horizontal arc
by utilizing $\cO_{1.1}$ and imitating the argument in Step-1 \& Step-2 of $\cO_{1.2}$.
Such a transformation of $\cE$ will necessarily employ the use of horizontal and vertical
exchange moves, and horizontal simplifications.
\qed


\noindent
{\bf $\cO_3$--Eliminating half discs.}
Similar to the count of $\cO_2$ operations,
there are three differing configurations of the tiled foliation which can be used to
detect the occurrence of our elementary moves for the
elimination of half discs.
Recalling Figure \ref{figure:sa arc types} and the arguments of
Lemmas \ref{lemma:discs that sa and scc bound} and \ref{lemma:discs with label 0},
we note that we have been careful to deal with both case-1 and case-2 type-1 {{\bf sa}'s} individually.
However, our $\cO_3$ moves for eliminating half discs treat the \sa of these two
cases essentially the same.  So our discussion (and accompanying illustrations) will
describe the properties of these moves mostly in terms of case-1 type-1 \sa.
More important to this discussion is whether $x$ ``can see the half disc''
where $x \subset \notchdisc \cap C_1$.  (The possibilities for $x$ are
$c_{max}$, a horizontal arc of $\sheararcpres$, or a horizontal boundary arc of $\cR$.)
Specifically, let $c_\g \subset C_1 \cap \notchdisc$ be a type-1 \sa
with $h$ the associated horizontal arc of statement d. of Lemma \ref{lemma:discs that sa and scc bound}.
Moreover, let $v^\partial \subset H_{\theta^\partial}^{\ }$ be the $\partial$-vertical arc
that $c_\g$ is adjacent to along with the disc fiber in $\fib$.  Now we consider
the disc components of $(H_{\theta^\partial}^{\ } \cap \cT_{\infty}) \setminus \notchdisc$.
Suppose that $Z$ is a disc-component for which $h \cap H_{\theta^\partial}^{\ } \subset \partial \overline{Z}$.
(Notice that there is only one such $Z$.)
If $\partial \overline{Z}$ also contains $x \cap H_{\theta^\partial}^{\ }$ then $x$
{\em can see the half disc} associated with $h$.\qed

\begin{figure}[htb]
\centerline{\includegraphics[scale=.6, bb=0 0 630 564]{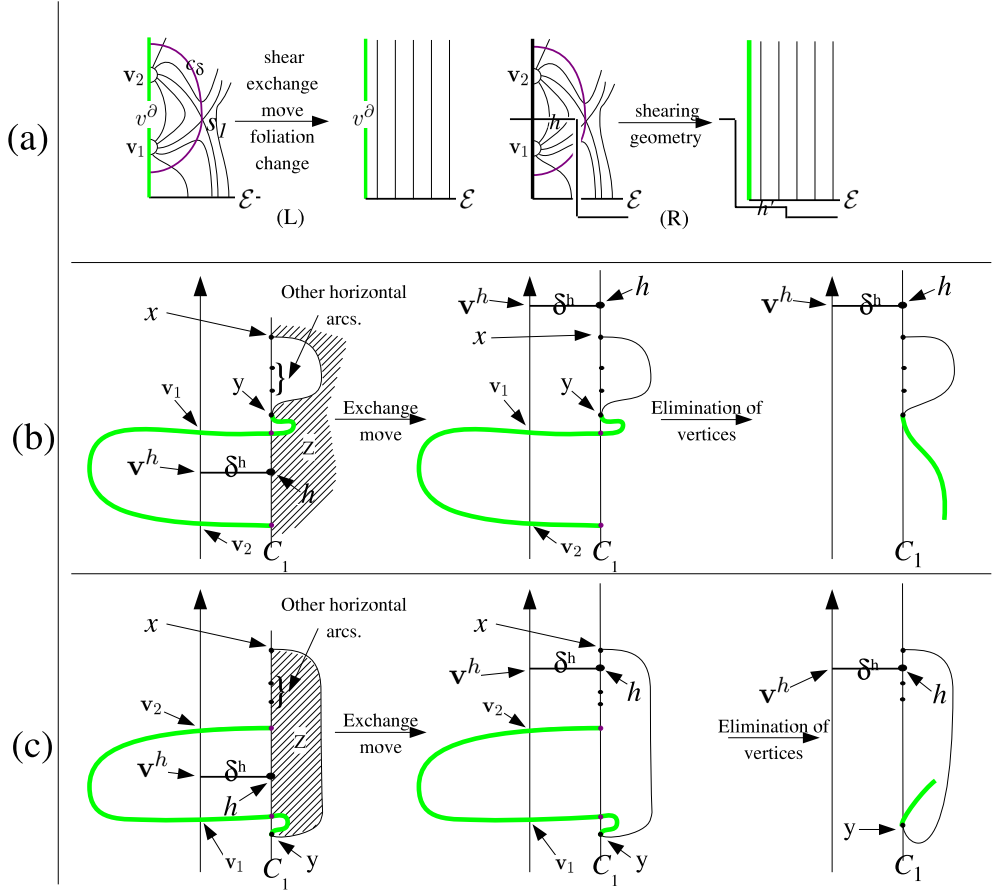}}
\caption{The illustrations should allow for the following interruption.  For operation $\cO_{3.1}$
we let $x = c_{max}$.  For $\cO_{3.2}$ we let $x$ be an edgepath in $\sheararcpres$ or a
horizontal arc in $\partial \cR$.  For $\cO_{3.3}$ we let $x$ be an edgepath (similar to
those in $\cO_{2.3}$) containing portions of $\partial$-vertical arcs.}
\label{figure:valence two half discs}
\end{figure}

\noindent
\underline{$\cO_{3.1}$--Eliminating half discs near $c_{max}$ which $c_{max}$ can see.}
Let $h \subset \sheararcpres$ be a portion of a horizontal arc associated
with a half disc as described in statement c. of Lemma \ref{lemma:discs that sa and scc bound}.
Assume $\delta^{h}$ has an intersection arc $\g \subset \delta^{h} \cap \notchdisc$
which is {\em innermost} on the component of $\notchdisc \setminus \cR$
containing $D_{+1}$.  That is, for every leaf $a_\theta$ in the radial
foliation of $\notchdisc$ intersecting $\g$ we have that when traversing
$a_\theta$ with $\ver_{max}$ as the starting point, the first time $a_\theta$ intersects
any coning disc is at the point $a_\theta \cap \g $.
For such a choice of coning disc we apply the transition
$ (\sheararcpres,\notchdisc) \Ttran (\sheararcpres,\tilenotchdisc)$
to produce our tiled notch disc.  This transition necessarily involves positioning
the $\partial$-vertical arc $v^\partial$ that is adjacent to $c_\g$
so that it intersects $\axis$ at (at least) two points.

Referring to Figure \ref{figure:valence two half discs},
this transition yields a tiled foliated subdisc of $\tilenotchdisc$
containing $\ver_{max}$ with two additional vertices, $\ver_1$ and $\ver_2$, and a
singularity $\sin_0$.  (In Figure \ref{figure:valence two half discs}(a) $\ver_{max}$ is at
the ``top of the page'' along with $c_{max}$ and not shown.  However, in the
Figure \ref{figure:valence two half discs}(b) \& (c) cross-sections the
labeled point $x (\subset C_1)$ is to be
interrupted as $c_{max} \cap H_{\theta^\partial}^{\ }$.)
$\ver_{max}$ and $\ver_2$ are positive vertices, whereas $\ver_1$ has negative parity.
The \sa $c_\gamma \subset \tilenotchdisc \cap C_1$
passes through this singularity and
splits off $\ver_1$ and $\ver_2$.  Between $\ver_1$ and $\ver_2$ on the axis $\axis$ there is the vertex, $\ver^{h}$, of coning disc $\delta^{h}$.
And, the $b$-arcs adjacent to $\ver_{max}$ \&
$\ver_2$ intersect $\tilenotchdisc \cap C_1$ only at $c_{max}$ \& $c_\g$.

More importantly, we assume that $c_{max}$ can see $h$ (along with $c_\g$).
We illustrate this assumption in
the cross-sections of Figure \ref{figure:valence two half discs}(b) \& (c).  The
disc region $Z \subset [(H_{\theta^\partial}^{\ } \cap \cT_{\infty}) \setminus \notchdisc] \subset H_{\theta^\partial}^{\ }$ (shaded regions in figure) has that both $c_{max} \cap \partial \overline{Z}$
and $h \cap \partial \overline{Z}$ are nonempty.  (Our cross-sections
are in fact the union of the two disc fibers of $\fib$ that contain
the $v^\partial \subset \partial \tilenotchdisc$.)
Having $c_{max}$ seeing $h$ allows us to eliminate $c_\g$ by a shear horizontal
exchange move of $h$ as described next. 

We can perform a shear exchange move which will (similar to $\cO_{2.3}$) ``shear'' our
coning disc $\delta^h$ into two coning disc $\delta^{h\prime}$ and $\delta^{h\prime \prime}$.
Let $\delta^{h\prime}$ be the coning disc associated with our \sa $c_\g$.
Now isotop $\delta^{h\prime}$ so that the resulting $\ver^{h\prime}$ and $\ver_{max}$
are at consecutive horizontal positions.  After the isotopy
$\g$ and any arcs of $\delta^{h\prime} \cap C_1$ between
$\ver^{h\prime}$ and $\g$ on $\delta^{h\prime}$
will be eliminated.  Thus, the intersection curve $c_\g$ can be
eliminated to reduce $|[ \notchdisc \setminus \cR ] \cap [ C_1 \setminus (C_1 \cap \cI) ]|$, and
we replace $\tilenotchdisc$ with a radially foliated notch disc $\notchdiscprime$.
(There may be additional \sa's of $\tilenotchdisc \cap C_1$ that are eliminated
by this shear exchange move since $\g$ was innermost on $\notchdisc$ with respect
to $\ver_{max}$ but may not be innermost on $\delta^{h}$
with respect to $\ver^{h}$.)

It is our ``seeing'' assumption that allows that us to perform
our shear horizontal exchange move on $h^\prime$ enabling the elimination of
$c_\g$ as illustrated in cross-section sequences.  More precisely, throughout all disc fibers
in the angular support of $c_\g$ we have that $c_{max}$ can see $h$ and $c_{max}$.  Thus,
we have the $3$-ball that is the shaded-region-$Z$ cross the angular-support-of-$c_\g$
in which the shear horizontal exchange move on $h^\prime$ can be performed. 

Similar to the argument for $\cO_{2.1}$ there are
two cases depending on the initial ordering of vertices
on the axis:  $\ver_2 < \ver^{h} < \ver_1 < \ver_{max}$
(Figure \ref{figure:valence two half discs}(b))
and $\ver_1 < \ver^{h} < \ver_2 < \ver_{max}$
(Figure \ref{figure:valence two half discs}(c)).  For Figure \ref{figure:valence two half discs}(b) a shear exchange move will position $h^\prime$ so that its horizontal position is
just above $\ver_{max}$.  Whereas,for Figure \ref{figure:valence two half discs}(c) a shear exchange move will position $h^\prime$ so that its horizontal position is
just below $\ver_{max}$.

As in $\cO_{2.1}$, our difficulty is in seeing that these two cases in fact correspond
to a sequence of vertical exchange moves followed by shear horizontal exchange moves 
as defined in \S\ref{subsection:rectangular diagrams}.
In particular, we need to argue that through the uses of
vertical exchange moves the nested feature for the shear horizontal exchange move
is satisfied.  But, the argument for using vertical exchange moves to move
vertical arcs past
the angular position of $\sin_0$ (so as to ensure
nested) is exactly the same as the argument in $\cO_{2.1}$, and we will
not repeat it here.

Again, we are not quite done since the isotopy in Figure \ref{figure:valence two}(b)
leaves the horizontal position of $h^\prime$ above that of $\ver_{max}$.  To resolve this
remaining issue we note that the horizontal position of $h^\prime$ is now maximal and we
can perform a first flavor horizontal exchange move on $h^\prime$ making its
horizontal position minimal.  Notice this will not change the value of
$|[ \notchdisc \setminus \cR ] \cap [ C_1 \setminus (C_1 \cap \cI) ]|$.
Again, we do not forget to place all of the previously
moved vertical arcs back into their original position.

We conclude that through a sequence of vertical exchange moves and shear
horizontal exchange moves we can eliminate half discs that
are associated with innermost intersections of cone discs
on the component of $\notchdisc \setminus \cR$ containing $D_{+1}$
when $c_{max}$ can see $h$.
The resulting notch disc will have fewer intersection components with
$C_1 \setminus \cI$.\qed

\noindent
\underline{$\cO_{3.2}$--Elimination of half disc near edgepath
$\cE$ that $\cE$ can see, away from $\partial$-vertical arcs.}
Again, let $h \subset \sheararcpres$ be a portion of a horizontal arc associated
with a half disc, $c_\g$, as described in statement c. of Lemma \ref{lemma:discs that sa and scc bound}.
Assume
the coning disc $\delta^h$ has an intersection arc $\g \subset \delta^h \cap \notchdisc$
which is {\em nearest} on a $D_0$ or the $D_{+1}$ component, $D^\star \subset \notchdisc \setminus int(\cR)$,
to an edge-path $\cE \subset \partial D^\star$.
That is, for every leaf $a_\theta$ in the radial
foliation of $\notchdisc$ intersecting $\g$ we have that it also intersects $\cE$ and, when traversing
$a_\theta$ from $a_\theta \cap \cE$ to $a_\theta \cap \g$, the first time $a_\theta$ intersects
any coning disc is at $a_\theta \cap \g $.  (Again, picture $\cE$
having slightly larger angular support than that of $c_\g$.)
Additionally, assume that the only $\partial$-vertical arc $\cE$ is adjacent to
is $v^\partial$, the $\partial$-vertical arc that $c_\g$ is adjacent to.
(Thus, $int(\cE)$ does not intersect any $\partial$-vertical arc.)
Finally, assume that $\cE$ can see $h$.

As in $\cO_{2.2}$, $\cE$ is a portion of either $\sheararcpres \cap D^\star$ or
a horizontal boundary arc of a component of $\cR$.  And, the most straight forward situations
are when that $\cE$ is a portion of either a horizontal arc of $\sheararcpres$ or
a horizontal boundary arc of a component of $\cR$.

We apply the transition
$ (\sheararcpres,\notchdisc) \Ttran (\sheararcpres,\tilenotchdisc)$
for this choice of coning disc.
Referring to Figure \ref{figure:valence two half discs},
this transition yields a tiled foliated subdisc of $\tilenotchdisc$
containing vertices, $\ver_1$ and $\ver_2$, and a
singularity $\sin_0$.  (In Figure \ref{figure:valence two half discs}(a)
the edgepath $\cE$ in our straight-forward-situation is the labeled horizontal arc.
In Figure \ref{figure:valence two half discs}(b) \& (c) cross-sections the
labeled point $x (\subset C_1)$ is to be
interrupted as $\cE \cap H_{\theta^\partial}^{\ }$.)
$\ver_2$ and $\ver_1$ are of opposite parity.
Again, the \sa $c_\gamma \subset \tilenotchdisc \cap C_1$
passes through this singularity and
splits off $\ver_1$ and $\ver_2$.  Between $\ver_1$ and $\ver_2$ on the axis $\axis$ there is the vertex, $\ver^{h}$, of coning disc $\delta^{h}$.
And, the $a$-arcs adjacent to $\cE$ \&
$\ver_1$ intersect $\tilenotchdisc \cap C_1$ only at $\cE$ \& $c_\g$.

Again, the important assumption is that $\cE$ can see $h$ (along with $c_\g$).
We see this in
the cross-sections of Figure \ref{figure:valence two half discs}(b) \& (c)
when interrupting the labeled point $x$ as $(\cE \cap C_1) \cap H_{\theta^\partial}^{\ }$.  The
disc region $Z \subset [(H_{\theta^\partial}^{\ } \cap \cT_{\infty}) \setminus \notchdisc] \subset H_{\theta^\partial}^{\ }$ has that both $\cE \cap \partial \overline{Z}$
and $h \cap \partial \overline{Z}$ are nonempty.
Inside the $3$-ball that is the shaded-region-$Z$ cross the angular-support-of-$c_\g$
we can perform the shear horizontal exchange move on $h^\prime$. 

We perform a shear exchange move by ``shearing''
$\delta^h$ into two coning disc $\delta^{h\prime}$ and $\delta^{h\prime \prime}$.
$\delta^{h\prime}$ will be the coning disc associated with our \sa $c_\g$.
Now isotop $\delta^{h\prime}$ so that the resulting $\ver^{h\prime}$ and $\cE$
are at consecutive horizontal positions.  After the isotopy
$\g$ and any arcs of $\delta^{h\prime} \cap C_1$ between
$\ver^{h\prime}$ and $\g$ on $\delta^{h\prime}$
will be eliminated.  Thus, the intersection curve $c_\g$ can be
eliminated to reduce $|[ \notchdisc \setminus \cR ] \cap [ C_1 \setminus (C_1 \cap \cI) ]|$, and
we replace $\tilenotchdisc$ with a radially foliated notch disc $\notchdiscprime$.
(Again, there may be additional \sa's of $\tilenotchdisc \cap C_1$ that are eliminated
by this shear exchange move since $\g$ was innermost on $\notchdisc$ with respect
to $\ver_{max}$ but may not be innermost on $\delta^{h}$
with respect to $\ver^{h}$.)

As in $\cO_{2.2}$ and $\cO_{3.1}$, our difficulty is in seeing that this isotopy of $h^\prime$
corresponds to a sequence of vertical exchange moves
followed by shear horizontal exchange moves 
as defined in \S\ref{subsection:rectangular diagrams}.
The argument for using vertical exchange moves to move vertical arcs past
the angular position of $\sin_0$ (so as to ensure
nested) is exactly the same as the argument in $\cO_{2.2}$.
The value of $|[ \notchdisc \setminus \cR ] \cap [ C_1 \setminus (C_1 \cap \cI) ]|$
is reduced with the elimination of $c_\g$.
Again, we do not forget to place all of the previously
moved vertical arcs back into their original position.

The remaining case is when $\cE$ is in $\sheararcpres$ but is not a portion of a single
horizontal arc.  As in $\cO_{2.2}$,
we can transform $\cE$ into a portion of such a horizontal arc
by utilizing $\cO_{1.1}$ and imitating the argument in Step-1 \& Step-2 of $\cO_{1.2}$.
Such a transformation of $\cE$ will necessarily employ the use of horizontal and vertical
exchange moves, and horizontal simplifications.

We conclude that through a sequence of vertical exchange moves and shear
horizontal exchange moves we can eliminate half discs in this situation.
\qed

\noindent
\underline{$\cO_{3.3}$--Elimination of half disc near edgepath
$\cE$ that $\cE$ can see, adjacent to $\partial$-vertical arcs.}
Our assumptions are exactly as those in $\cO_{3.2}$ except we allow for
$\cE$ to be adjacent to $\partial$-vertical arcs that $c_\g$ is not adjacent to.
We list the initial possibilities: (i) $\notchdisce$ and $int(\cE)$ intersects
one containing $\partial$-vertical arcs; and, (ii)
$\notchdiscf$ and $int(\cE)$ intersects with one or two $\partial$-vertical arcs.

The argument is a combining of the arguments of $\cO_{2.3}$ and $\cO_{3.2}$.
The details at this points should be readily accessible to the reader.
We concluded that through a sequence of elementary moves half discs in
this situation can be eliminated.  Moreover, the count
$|[ \notchdisc \setminus \cR ] \cap [ C_1 \setminus (C_1 \cap \cI) ]|$
is reduced.\qed

\noindent
{\bf $\cO_{4}$--Using shear vertical simplification on product of subdisc fibers.}
We generalize the idea of using virtual rectangular discs
to perform shear vertical exchange moves.  We assume that
the foliation of $\notchdisc \setminus (\notchdisc \cap C_1)$ contains no
whole or half discs.  Let
$[\theta_1 , \theta_2] \subset S^1$ be an angular interval
such that $H_{\theta_1^{\ }}^{\ }$
and $H_{\theta_2^{\ }}^{\ }$ are disc fibers of $\fib$ containing
$\partial$-vertical arcs of $\partial \notchdisc$ having differing
edge assignments.  To list the obvious cases: (A) there is a back
edge, $v^\partial_1$, having position $\theta_1$ and a front edge,
$v^\partial_2$ having position $\theta_2$;
or, (B) there is a front
edge, $v^\partial_1$, having position $\theta_1$ and a back edge,
$v^\partial_2$, having position $\theta_2$.

We assume that over the angular interval $[\theta_1 , \theta_2]$ our
$\partial \notchdisc$ contains no vertical arcs.  If we have case (A)
(respectively, case (B)) then any \sa having endpoints on
$v^\partial _1$ and $v^\partial_2$ with angular support $[\theta_1, \theta_2]$
will be a type-2 (respectively, type-4).

We consider the components of 
$ 
[\cup_{\theta \in (\theta_1 , \theta_2)} H_\theta \cap \cT_{\infty} ] \setminus \notchdisc.
$
The reader should realize that these components are topological
open $3$-ball.  We wish to distinguish one of these $3$-ball components.
In particular, $\hat{B}$ will be the $3$-ball component whose closure contains the
horizontal subarcs of $ \sheararcpres \setminus \partial \notchdisc$
whose endpoints correspond to the endpoints of $v^\partial_1$ and $v^\partial_2$.

To expand on this terse description, we have case (B) (respectively, case (A)) as follows.
Referring back to Figure \ref{figure:boundary of disc} we observe that the
left $\partial$-vertical arc is a front edge with its ``bottom of the page'' endpoint
(respectively, ``top of the page'' endpoint)
having the braid `entering'
(respectively, `exiting')
$\partial \notchdisc$.  The right
$\partial$-vertical arc is a back edge with the braid `exiting'
(respectively, `entering')
$\partial \notchdisc$ at the
``top of the page''
(respectively, ``bottom of the page'')
endpoint.  $\hat{B}$ is the component
whose closure contains the portion of the braid entering and exiting
(respectively, exiting and entering)
$\partial \notchdisc$.

For both case (A) and (B) it should be clear that each $3$-ball component has a
natural product structure, $ D \times (\theta_1 , \theta_2)$ where $D$
is an open disc component of
$ (H_{\theta} \cap \cT_{\infty}) \setminus  \notchdisc \subset
H_\theta$ for $\theta \in (\theta_1 , \theta_2)$.
In particular, for any $3$-ball component $ D \times (\theta_1 , \theta_2)$
which is not $\hat{B}$, if $v \subset \overline{D \times (\theta_1 , \theta_2)}$
is a vertical arc of $\sheararcpres$
having angular position nearest to, say, $\theta_1$, then
we have a virtual rectangular disc $ R = v \times [\theta_1 , \theta_2]
\subset \overline{D \times (\theta_1 , \theta_2)}$.  Moreover, $R \cap \notchdisc = \emptyset$
and we can use $R$ to push $v$ out of $\overline{D \times (\theta_1 , \theta_2)}$.
As with our previous use of virtual rectangular discs, this push will
correspond to a sequence of vertical exchange moves followed by a single
shear vertical simplification.\qed
%

\section{Proof of Theorem \ref{theorem:main result}.}
\label{section:proof of theorem 1}

With the operations of
\S\ref{subsection:the tiling machinery on notchdisc} in place we are now in a position to
prove Proposition \ref{proposition:obvious disc} and, thus, Theorem \ref{theorem:main result}.
We start by enhancing our complexity measure on $\sheararcpres$ so as to
incorporate information from our notch disc.  Specifically,
we define the lexicographically
ordered complexity $\cC(\sheararcpres , \notchdisc) =
(n_1,n_2,n_3)$ where: $n_1 = \cC(\sheararcpres)$;
$n_2$ is the number of components of
$[ \notchdisc \setminus \cR ] \cap [ C_1 \setminus (C_1 \cap \cI) ]$; and,
$n_3 $ is the number of vertical arcs in $\partial \notchdisc$.
Then we will have an
obvious destabilizing disc when $\cC(\sheararcpres, \notchdiscd) = (\_ ,0,1)$;
an obvious exchange disc when $\cC(\sheararcpres, \notchdisce) = (\_ , 0,0)$;
and, an obvious flyping disc when $\cC(\sheararcpres, \notchdiscf) = (\_ , 0,1)$.
Ideally, for an obvious disc all of the vertically arcs not in $\partial \notchdisc$
have been pushed into $\cI$ by shear-vertical simplification.  Thus, although not
necessary for recognizing an obvious disc, the absolute minimal
complexity has: $\cC(\sheararcpres, \notchdiscd) = (2,0,1)$;
$\cC(\sheararcpres, \notchdisce) = (0, 0,0)$;
and, $\cC(\sheararcpres, \notchdiscf) = (1, 0,1)$.

Now, let $X$ be a closed braid which admits one of our three braid isotopies and
let $X \ntran \arcpres$ be a transition to an arc presentation.  Since associated
to $X$ there was a disc $\Delta_\v$ depicting a braid isotopy there
exists a $\notchdisc$ associated to $\arcpres$ also depicting the same braid isotopy.
We can apply Lemma \ref{lemma:first simplification of notch disc} without
changing $\cC(\arcpres)$.  The position of the $\partial$-vertical arcs
of $\notchdisc$ gives us the existence and position of our shearing intervals $\cI$.
We can then apply Lemmas
\ref{lemma:control of arc intersections}, \ref{lemma:discs that sa and scc bound} and
\ref{lemma:discs with label 0} so as to place
$\notchdisc$ into a position for calculating the values of our $3$-tuple
complexity measure $\cC(\sheararcpres, \notchdisc)$ pair.

The argument for achieving minimal complexity should be clear.  Applying a
sequence of $\cO_2$ and $\cO_3$ operations we can eliminate all whole and half
discs while reducing $n_2$ and, thus, $\cC(\sheararcpres , \notchdisc)$.
Specifically, if a half disc is near an edgepath that is either
a portion of $\sheararcpres$ or a horizontal boundary arc of a component of $\cR$
then the edgepath can see the horizontal arc associated with the associated \sa.
So starting with \scc that are close to $c_{max}$ or near an edgepath; or \sa that
in near and can be seen by an edgepath or $c_{max}$, 
we can iterate these two
operations until every \scc and \sa is eliminated.

Next applying the appropriate $\cO_1$ sub-operation we can simplify the
middle boundary of $\notchdisc$ until $n_3$ is either $0$ (for $\notchdisce$) or
$1$ (for $\notchdiscd$ or $\notchdiscf$).  Along the way we may have applied
some number of horizontal, vertical and shear vertical simplifications to
reduce $n_1 = \cC(\sheararcpres)$.  Since all of our $\cO$-operations involve
the use of elementary moves the four statements of Proposition
\ref{proposition:obvious disc} are satisfied.

We pause briefly to discuss the nature of the finiteness of our
complexity measure.  In view of Remarks \ref{remark:elementary exchange moves and flypes},
\ref{remark: equivalency class for EM and F} and \ref{remark: double destabilization}
plus the ``connect the dots'' discussion just prior to the statement of
Theorem \ref{theorem:main result}, for
a specified $\cI$ and middle boundary edge path there may be infinitely many
possible $\partial$-vertical arcs yielding infinitely many possible
$\notchdisc$.  Thus, an initial value of $n_2$ can be arbitrarily large.
However, initial values of $n_1$ and $n_3$ are bounded.  Moreover,
for any values of $n_1$ and $n_3$ there are only finitely many possible
cyclic orderings of horizontal and vertical arcs.  Thus, within this
finite collection of orderings there will exists some $\notchdisc$ having
its associated $n_2$ being zero.  With an $n_2 = 0$ we can reduce
$n_3$ to its associated minimal value using $\cO_1$ operations.

Having established Proposition \ref{proposition:obvious disc} we have
also established statements 1 and 2 of Theorem \ref{theorem:main result}.
In order to establish statement 3 of Theorem \ref{theorem:main result} we
need to think about how the rectangles of $\cR$ alter the
$\sheararcpresprime$ of statement 3 and, thus, the standard projection of
a $X^\prime$ coming from $\sheararcpresprime \btran X^\prime$. 

It is simplest to deal with each braid isotopy individually in turn.

\noindent
\underline{A destabilizing disc}--We briefly summarize the salient features
of the pair $(\sheararcpresprime, \notchdiscprimed)$
(with $\cC(\sheararcpresprime, \notchdiscprimed) = (\_ , 0 , 1)$)
that established the validity
of Proposition \ref{proposition:obvious disc} for a destabilization.
Since $\partial \notchdiscprimed$ has a single $\partial$-vertical arc we know that
there are no type-2, -3, or -4 \sa in $\notchdiscprimed \cap C_1$.
In fact, as a result of the simplification in $\cO_{1.2}$ we know that:
\bi
\item[1.] $int(\notchdiscprimed) \cap C_1 = c_{max}$.
\item[2.] $\partial \notchdiscprimed = h_1 \bu v^\partial \bu h_2 \bu v_2 \bu$.
\item[3.] If $\notchdiscprimed$ corresponds to a positive (respectively, negative)
destabilization then $v^\partial$ is a front (respectively, back) edge, and
$h_2$ (respectively, $h_1$) has maximal horizontal position over all horizontal
arcs of $\sheararcpresprime$.
\item[4.] Recalling the decomposition $\notchdiscprimed = D_{+1} \cup_{c_{max}} N$
where $\partial N = c_{max} \cup \partial \notchdiscprimed$ we have that
$N = R_1 \cup R_2$ where:
\bi
\item[a.] $R_1$ has the same angular support as $h_1$
($[\theta_2 , \theta^\partial]$) and vertical support
$[z^{h_1}_{{\ }^{\ }}, z_{max}]$.
\item[b.] $R_2$ has the same angular support as $h_2$
($[\theta^\partial , \theta_2]$)
and the vertical
support $[z^{h_2}_{{\ }^{\ }}, z_{max}]$.
\ei
\ei
There are the positive and negative destabilizing cases to argue.  The arguments
are similar, so we will present the positive destabilization case and leave the
details of the negative case to the reader.

The existence of $N$ divides $\cT_{\infty}$ up into three potential
regions where braiding of $\sheararcpresprime \setminus \partial \notchdiscprimed$
can occur: $\cH_1$ is over the interval $[\theta^\partial , \theta_2]$ and below
$z^{h_2}_{{\ }^{\ }}$; $\cH_2$ is over the interval $[\theta_2 , \theta^\partial]$ and above
$z^{h_1}_{{\ }^{\ }}$; and, $\cH_3$ is over the interval $[\theta_2 , \theta^\partial]$ and below
$z^{h_1}_{{\ }^{\ }}$.  We will achieve a standard projection easily seen as
equivalent (through type-III Reidemeister moves)
to $W \sigma_{n-1}$ when: $\cH_3$ contains no vertical arcs, thus no braiding;
and $\cH_2$ contains a single vertical arc which will be adjacent to $h_2$.

Having no vertical arcs in $\cH_3$ is easily achieved by vertical exchange moves.
Specifically, every vertical arc in $\cH_3$ has vertical support either
below or above $z^{h_1}_{{\ }^{\ }}$.  Thus, using the idea in Step-3 of
$\cO_{1.3}$ for constructing virtual rectangular discs,
through a sequence of vertical exchange moves
we can push all vertical arcs out of $\cH_3$ and into $\cH_1$.

Finally, we consider the vertical arcs in $\cH_2$.  Since $z^{h_2}_{{\ }^{\ }}$
is maximal, there must be a vertical arc $v^\prime \subset \sheararcpresprime$
having angular position in
$(\theta_2 , \theta^\partial)$ and adjacent to $h_2$.  Let $\theta^\prime$ be
it angular position.  Then repeating the vertical exchange moves arguments
of Step-3 of $\cO_{1.2}$ we can construct virtual rectangular discs that are in the
complement of $\notchdiscprimed$ for pushing all the vertical arcs in $\cH_2$
having angular position in $(\theta^\prime , \theta^\partial)$
forward past $v^\partial$, and in $\cH_2$ having
angular position in $(\theta_2 , \theta^\prime)$ backward past $v_2$.
All such pushes move vertical arcs into $\cH_1$.

All of these vertical exchange moves do no change the complexity measure of our
pair $(\sheararcpresprime, \notchdiscprimed)$.  Thus, we have established
statement 3 of Theorem \ref{theorem:main result} for a destabilization.

\noindent
\underline{An exchange disc}--We consider
a pair $(\sheararcpresprime, \notchdiscprimee)$
with $\cC(\sheararcpresprime, \notchdiscprimee) = (\_ , 0 , 0)$
which establishes the validity
of Proposition \ref{proposition:obvious disc} for an exchange move.
Since $\partial \notchdiscprimee$ has a two $\partial$-vertical arcs of
differing edge assignment we know that
there are only type-2 and -4 \sa in $\notchdiscprimee \cap C_1$.
From our $\cO_{1.3}$ and $\cO_{1.4}$ operations we know that it is useful
to breakup our discussion into two cases: $\cR = \emptyset$ and $\cR \not= \emptyset$.
When $\cR = \emptyset$ we know from $\cO_{1.3}$ that:
\bi
\item[1.] $\notchdiscprimee \cap C_1 = c_{max} \cup \partial \notchdiscprimee$.
\item[2.] $\partial \notchdiscprimee = h_1 \bu v^\partial_1 \bu h_2 \bu v^\partial_2 \bu$.
For convenience, we say $v^\partial_1$ is a back edge and $v^\partial_2$ is
a front edge with $z^{h_2}_{{\ }^{\ }} < z^{h_1}_{{\ }^{\ }}$.
\item[3.] Recalling the decomposition $\notchdiscprimee = D_{+1} \cup_{c_{max}} N$
where $\partial N = c_{max} \cup \partial \notchdiscprimee$ we have that
$N = R_1 \cup R_2$ where:
\bi
\item[a.] $R_1$ has the same angular support as $h_1$
($[\theta^\partial_2 , \theta^\partial_1]$) and vertical support
$[z^{h_1}_{{\ }^{\ }}, z_{max}]$.
\item[b.] $R_2$ has the same angular support as $h_2$
($[\theta^\partial_1 , \theta^\partial_2]$)
and the vertical
support $[z^{h_2}_{{\ }^{\ }}, z_{max}]$.
\ei
\ei
The decomposing disc $R_1 ,R_2 \subset N$ divides $\cT_{\infty}$ up into four potential
regions where braiding of $\sheararcpresprime \setminus \partial \notchdiscprimed$
can occur: $\cH_1^+$ is over the interval $[\theta^\partial_2 , \theta^\partial_1]$ and above
$z^{h_1}_{{\ }^{\ }}$ (and under $R_1$ in $\cT_{\infty}$);
$\cH_1^-$ is over the interval $[\theta^\partial_2 , \theta^\partial_1]$ and below
$z^{h_1}_{{\ }^{\ }}$ (and below $R_1$);
$\cH_2^+$ is over the interval $[\theta^\partial_1 , \theta^\partial_2]$ and above
$z^{h_2}_{{\ }^{\ }}$ (and under $R_2$ in $\cT_{\infty}$); and, $\cH_2^-$ is over the interval
$[\theta^\partial_1 , \theta^\partial_2]$ and below
$z^{h_2}_{{\ }^{\ }}$ (and below $R_2$).  We will achieve a standard projection for $X^\prime$ easily seen as
having the needed braid word structure $W U$
($W$ a word in $\cW^t$ and $U$ a word in $\cU^s$ with  $s \leq t$)
when $\cH_1^+$ and $\cH_2^-$ contains no vertical arcs, thus no braiding.
The number of horizontal arcs of $\sheararcpresprime$ at the
angle, say $\theta^\partial_1$ having horizontal position
in the interval $(z^{h_2}_{{\ }^{\ }} , z^{h_1}_{{\ }^{\ }})$ corresponds to
$t-s+2$.

For two angular intervals, $[\theta^\partial_2 , \theta^\partial_1]$ and
$[\theta^\partial_1 , \theta^\partial_2]$, we notice that we are in the
situation described in $\cO_{4}$.  Moreover, for
$[\theta^\partial_2 , \theta^\partial_1]$ (respectively,
$[\theta^\partial_1 , \theta^\partial_2]$) we have that
$\cH_1^-$ (respectively, $\cH_2^+$) corresponds to the distinguish $\hat{B}$ component.
Thus, by $\cO_{4}$ we can push all vertical arcs in $\cH_1^+$
(respectively, $\cH_2^-$) into $\cH_2^+$ (respectively, $\cH_1^-$) using a
sequence of vertical exchange moves.

Proceeding to the case where $\cR \not= \emptyset$ we acknowledge that
using only the collection of elementary moves on arc presentations with
shearing intervals, it may not be possible to place $\sheararcpresprime$ in a position
such that for $\sheararcpres \btran X^\prime$ we have $\b(X^\prime) = WU$
(where $W \in \cW^t$ and $U \in \cU^s$,  $s \leq t$).  To be clear,
it is readily seen that by allowing additional moves on the our shear arc presentations
that correspond to the inverse of horizontal simplification we can produce
presentations $\sheararcpresprimeprime$ such that $\cR = \emptyset$ and, by
the above argument,
for $\sheararcpresprimeprime \btran X^{\prime \prime}$ we have
the desired braid word form $\b(X^{\prime \prime}) = WU$.  Unfortunately, we lose
the feature that our complexity measure is non-increasing.

However, referring back to Remarks \ref{remark:elementary exchange moves and flypes}
and \ref{remark: equivalency class for EM and F} we note: first,
since $\sheararcpresprimeprime$ admits an exchange move, it admits a thin exchange move; and,
second, there are infinitely many exchange move discs associated with a
thin exchange move of $\sheararcpresprimeprime$.
In particular, from Remark \ref{remark: equivalency class for EM and F} we have the following.

Let $X^{\prime \prime}$ be an $n$-braid such that $\b(X^{\prime \prime}) = W U$, where
$W \in \cW^t$ and $U \in \cU^s$ with $t - s$ minimal, i.e. a thin exchange move.
Then for every integer pair $(p,q)$ the $n$-braid associated with the word
$$W^\prime U^\prime = [\tau_{[1, s-1]}^p W \tau_{[1,s-1]}^{-p}][\tau_{[t+2,n]}^q U \tau_{[t+2, n]}^{-q}]$$
admits the same thin exchange move.  Referring the reader back to
\S\ref{subsection:d,e,f discs} and our description of
exchange-move-disc/braid-presentation, we realize that for the closed $n$-braid
$X$ there are infinitely many exchange discs $\Delta_e (p,q)$ satisfying
conditions E-a through E-e with
$\partial \Delta_e (p,q) = \a^{\ }_{h_1} \bu [\a^{\partial}_{v_1} (p,q) ] \bu
\a^{\ }_{h_2} \bu [\a^{\partial}_{v_2}(p,q)] \bu$, i.e. all such discs have the
same two horizontal arcs, but differing $\partial$-vertical arcs indexed by
$(p,q)$.  Shifting this picture over to exchange-move-disc/arc-presentation,
we realize that there are infinitely many notch exchange discs, $\notchdisce(p,q)$,
satisfying conditions EA-a through EA-e.  Moreover,
$\partial \notchdisce(p,q) = h_1 \bu [v_1^\partial (p,q)] \bu h_2 \bu [v_2^\partial (p,q)]
\bu $.  That is, the boundary all such discs have the
same two horizontal arcs, but differing $\partial$-vertical arcs indexed by
$(p,q)$.

Using the radial foliations of $\Delta_e (p,q)$ we take note of the following
{\em winding behavior}.
As in condition E-a we have $\a_{v_1}^\partial (p,q) \subset H^{\ }_{\theta^\partial_1}$ and
$\a_{v_2}^\partial (p,q) \subset H^{\ }_{\theta^\partial_2}$.  So in the interval
$(\theta^\partial_1, \theta^\partial_2)$ there are, say, $s-1$ trivial strands
of $X^{\prime \prime}$; and,
in $(\theta^\partial_2 , \theta^\partial_1)$ there are $n - (t+1)$
trivial strands of $X^{\prime \prime}$.
Moreover, each leaf $a = \Delta_e(p,q) \cap H^{\ }_{\theta}
, \ \theta \in (\theta^\partial_1 , \theta^\partial_2)$
winds $p$ times around these $s-2$ trivial strands as seen in $H_{\theta}$---one
strand corresponds to the boundary of $\Delta_e(p,q)$.  (This is winding in the
sense that any ray in $H_{\theta}$ starting at the point corresponding to where the
stand intersects $H_{\theta}$ must algebraically intersect the leaf $a$ a minimal of $p$ times.)
Similarly, each leaf $a = \Delta_e(p,q) \cap H^{\ }_{\theta} , \ \theta \in (\theta^\partial_2 , \theta^\partial_1)$
winds $q$ times around these $n-(t+2)$ trivial strands as seen in $H_{\theta}$.

This winding behavior extends to all braid representatives in $\cB_n(X^{\prime \prime})$.
Specifically, for any fixed representative $X \in B_n (X^{\prime \prime} )$, since
$X$ is braid isotopic to $X^{\prime \prime}$ there is an ambient isotopy of
$S^3 \setminus \axis$ taking the family $\Delta_e (p,q)$ to exchange discs for
$X$.  Moreover,  for $X$ there must exists a pair $(P_X , Q_X)$ such that
for $P_X < |p|$ and $ Q_X < |q| $ we have
two similar statements of {\em winding behavior}:
\bi
\item[$(\theta^\partial_1 , \theta^\partial_2)$] As $|p|$ increases, each leaf $a = \Delta_e(p,q) \cap H^{\ }_{\theta} , \ \theta \in (\theta^\partial_1 , \theta^\partial_2)$
increasingly winds (geometrically) around these $s-1$ strands as seen in $H_{\theta}$.
\item[$(\theta^\partial_2 , \theta^\partial_1)$] As $|q|$ increases, each leaf $a = \Delta_e(p,q) \cap H^{\ }_{\theta} , \ \theta \in (\theta^\partial_2 , \theta^\partial_1)$
increasingly winds (geometrically) around these $n-(t+1)$ strands as seen in $H_{\theta}$.
\ei
Now taking a transition $X \ntran \arcpres$ we see that we have these two winding statements
also for $\notchdisce (p,q)$ with $P_X < |p|$ and $ Q_X < |q| $.

Next, consider an infinite collection of pairs
$(\sheararcpres, \notchdisce(p,q))$ which are associate as above to a thin exchange move.
In particular, for any two $2$-tuples, $(p,q)$ and $(p^\prime , q^\prime)$,
we have $\partial \notchdisc(p,q) \cap \sheararcpres =
\partial \notchdisc(p^\prime , q^\prime) \cap \sheararcpres $.
For any fixed $(p,q)$ we can apply our $\cO$-operations to reduce
the complexity measure for the shear-arc-presentation/notch-disc pair
to $(\_ , 0 , 0)$.
Now recall the observation just prior to Theorem \ref{theorem:main result}
that there are only finitely many shear arc presentations resulting from applying
a sequence of elementary moves to $\sheararcpres$.
This observation implies that there exists a resulting $\sheararcpresprime$
and an infinite set of $2$-tuples,
$\cS \subset \mathbb{Z} \times \mathbb{Z}$ such that:
\bi
\item[1.] the values $|p|$ \& $|q|$ are unbounded over all $(p,q) \in \cS$;
\item[2.] $h_1 \cup h_2 = \sheararcpresprime \cup \partial \notchdiscprimee (p,q)$
for all $(p,q) \in \cS$;
\item[3.] $\cC(\sheararcpresprime, \notchdiscprimee (p,q)) = (\_ , 0 ,0)$ for
all $(p,q) \in \cS$.
\ei
The reader should realize that increasing the values of $|p|$ \& $|q|$ causes
an increase in the number of type -2 \& -4 \sa's in $\notchdiscprimee (p,q) \cap C_1$ and, thus,
the value of $|\cR|$.
Our immediate goal is to show that $\sheararcpresprime$ is ``close'' to
being a shear arc presentation that transitions to a closed braid
which illustrates an obvious exchange move.  By ``close'' we will mean
that we may have to apply some number of horizontal/shear horizontal exchange moves.

To achieve this stated goal we proceed by understanding the implications of
the winding behavior in the radial foliation of our exchange discs.
We know that for large enough values $|p_0|$ and $|q_0|$, $(p_0,q_0) \in \cS$,
the radial leaves of $\notchdiscprimee(p_0,q_0)$ will wind around a maximal number
(due to having a thin exchange move) of
strands of $\sheararcpresprime$ within the specified angular intervals.  This winding implies the
existence of a subdisc $D^w_\theta \subset H_\theta$ for
$\theta \in (\theta^\partial_1 , \theta^\partial_2)$ (respectively, $\theta \in (\theta^\partial_2, \theta^\partial_1)$)
which satisfies the following.
\bi
\item[a.] $\partial D^w_\theta = e^1_\theta \bu e^2_\theta \bu$ where
$e^1_\theta \subset \notchdiscprimee (p_0,q_0) \cap H_\theta$
and $e^2_\theta \subset C_1 \cap H_\theta$ with 
$\theta \in (\theta^\partial_1, \theta^\partial_2)$ (respectively, $\theta \in (\theta^\partial_2, \theta^\partial_1)$);
\item[b.] $e^1_\theta \cap e^2_\theta = \partial e^1_\theta = \partial e^2_\theta$;
\item[c.] $D^w_\theta \cap \sheararcpresprime = int(D^w_\theta) \cap \sheararcpresprime$;
\item[d.] $|D^w_\theta \cap \sheararcpresprime| = s-1$ (respectively,
$|D^w_\theta \cap \sheararcpresprime| = n-(t+1)$).  That is, $D^w_\theta$ intersects
a maximal number of strands of $\sheararcpresprime$ since $\notchdiscprimee (p_0 ,q_0)$
is associated with a thin exchange move.
\ei

Since $\cC(\sheararcpresprime, \notchdiscprimee (p_0,q_0)) = (\_ , 0 ,0)$,
we have the product structure discussed in $\cO_4$ inside the $3$-ball
$B^{1,2} = \cup_{\theta \in (\theta^\partial_1 , \theta^\partial_2)}^{\ } D^w_\theta$ (or, $ B^{2, 1} = \cup_{\theta \in (\theta^\partial_2 , \theta^\partial_1)}^{\ } D^w_\theta$)
that allow us to assume after performing some number of shear vertical
simplifications that this $3$-ball contains no vertical arcs of
$\sheararcpresprime$.  Similarly, the bounded $3$-ball components of
$\cup_{\theta \in (\theta^\partial_1, \theta^\partial_2)}^{\ } (H_\theta \cap \cT_{\infty}) \setminus B^{1,2}$ have
product structures that allow us to apply $\cO_4$ removing vertical arcs.
There is still the possibility that there are horizontal arcs intersecting
the boundary of the closure of these $3$-balls components.  We will shortly
deal with this possibility.
(We have a similar statement for $B^{2,1}$ and interval $(\theta^\partial_2 , \theta^\partial_1)$.
However, if $B^{1,2}$/$B^{2,1}$ are simple enough there may be no such 
unbounded complementary $3$-balls components.)

Now consider the rectangular shaped disc components of
$$ [ \cup_{\theta \in (\theta^\partial_1 , \theta^\partial_2)} (C_1 \cap H_\theta) ] \setminus B^{1,2}. $$
We observe that no component contains horizontal subarc of $\sheararcpresprime$
having angular support containing $(\theta^\partial_1 , \theta^\partial_2)$.  Since, if such
a horizontal subarc occurred it would contradict our assumption that $D^w_\theta$
intersected a maximal number of trivial strands for $\theta \in (\theta^\partial_1 , \theta^\partial_2)$,
i.e. a contradiction of our thin exchange move assumption.
We can make a similar statement for components of
$$ [ \cup_{\theta \in (\theta^\partial_2 , \theta^\partial_1)} (C_1 \cap H_\theta) ] \setminus B^{2,1}. $$

We observe that, after performing some number of horizontal exchange moves
to ensure that the horizontal position of all horizontal arcs intersecting (say)
$B^{2,1}$ are above the horizontal position of all the horizontal arcs intersecting
$B^{1,2}$ we can assume that the vertical support of $B^{1,2}$ and $B^{2,1}$ are
disjoint.  (Since we now have $B^{1,2}$ and
$B^{2,1}$ containing no vertical arcs our necessarily nested condition 
is easily satisfied for horizontal exchange moves.)
Moreover, over the vertical support of $B^{1,2}$ (respectively, $B^{2,1}$)
within the angular interval $(\theta^\partial_1 , \theta^\partial_2)$ (respectively, $(\theta^\partial_2 , \theta^\partial_1)$)
we have that the only horizontal arcs are those contained in $B^{1,2}$
(respectively, $B^{2,1}$).

Finally, for the angular interval $(\theta^\partial_1 , \theta^\partial_2)$ (respectively,
$(\theta^\partial_2 , \theta^\partial_1)$) if there are any horizontal arcs below (respectively, above)
the vertical support of $B^{1,2}$ (respectively, $B^{2,1}$) then, through a
sequence of horizontal exchange moves and shear horizontal exchange moves,
we can move them above (respectively, below) $B^{1,2}$ (respectively, $B^{2,1}$).
(Again, the nested condition necessary for the performance of these exchange move
is readily seen to be satisfied.)  For the resulting $\sheararcpresprime$ it should
be apparent that for $\sheararcpresprime \btran X^\prime$ we have
$\b (X^\prime)$ admitting an exchange move.

\noindent
\underline{A flyping disc}--As with our previous two isotopies we consider
a pair $(\sheararcpresprime, \notchdiscprimef)$
with $\cC(\sheararcpresprime, \notchdiscprimef) = (\_ , 0 , 1)$
for an elementary flype move.
$\partial \notchdiscprimef$ has a three $\partial$-vertical arcs with
edge assignments front-back-front (respectively, back-front-back) for
a positive (respectively, negative) flype.  We will argue the positive
flype and leave the details of the negative flype to the reader.

From our $\cO_{1.4}$, $\cO_{1.5}$  and $\cO_{1.6}$ operations we know that it is useful
to breakup our discussion into two cases: $\cR = \emptyset$ and $\cR \not= \emptyset$.
When $\cR = \emptyset$ we know from $\cO_{1.5}$ that:
\bi
\item[1.] $\notchdiscprimef \cap C_1 = c_{max} \cup \partial \notchdiscprimef$.
\item[2.] $\partial \notchdiscprimef = h_1 \bu v^\partial_1 \bu h_2 \bu v^\partial_2 
\bu h_3 v^\partial_3 \bu h_4 \bu v_4 \bu$.
For convenience, we say $v^\partial_1$ \& $v^\partial_3$ are front edges and $v^\partial_2$ is
a back edge with $z^{h_1}_{{\ }^{\ }}$ minimal and $z^{h_3}_{{\ }^{\ }} < z^{h_2}_{{\ }^{\ }}$.
\item[3.] Recalling the decomposition $\notchdiscprimef = D_{+1} \cup_{c_{max}} N$
where $\partial N = c_{max} \cup \partial \notchdiscprimef$ we have that
$N = R_1 \cup R_2 \cup R_3$ where:
\bi
\item[a.] $R_1$ has the same angular support as $h_1$
($[\theta^{v_4}_{\ } , \theta^\partial_1]$) and vertical support
$[z^{h_1}_{{\ }^{\ }}, z_{max}]$.
\item[b.] $R_2$ has the same angular support as $h_2$
($[\theta^\partial_1 , \theta^\partial_2]$)
and the vertical
support $[z^{h_2}_{{\ }^{\ }}, z_{max}]$.
\item[c.] $R_3$ has the same angular support as $h_3$
($[\theta^\partial_2 , \theta^\partial_3]$)
and the vertical
support $[z^{h_3}_{{\ }^{\ }}, z_{max}]$.
\item[d.] $R_4$ has the same angular support as $h_4$
($[\theta^\partial_3 , \theta^{v_4}_{\ }]$)
and the vertical
support $[z^{h_4}_{{\ }^{\ }}, z_{max}]$.
\ei
\ei

We emulate the $\cH^\pm_i$-analysis in the exchange move situation.  In particular, we
have $N$ dividing $\cT_{\infty}$ up into eight potential
regions where braiding of $\sheararcpresprime \setminus \partial \notchdiscprimed$
can occur: $\cH_i^+$ is over the interval corresponding the the angular support
of $h_i$ and above $z^{h_i}_{{\ }^{\ }}$; and, $\cH_i^-$ is over the interval
corresponding the the angular support of $h_i$ and below
$z^{h_i}_{{\ }^{\ }}$, for $i \in \{ 1,2,3,4\}$.  $\cH_1^+ , \cH_2^+ , \cH_3^+ , \cH_4^+$
will be in $\cT_{\infty}$ under the discs $R_1 , R_2 , R_3 , R_4$, respectively.

Through the application of vertical exchange moves and shear horizontal exchange
moves in $\cH_1^-$ we obtain our assumption of statement 2 that $z^{h_1}_{{\ }^{\ }}$ is minimal.
(The assumption that $z^{h_3}_{{\ }^{\ }} < z^{h_2}_{{\ }^{\ }}$ will be dealt with
shortly.)
Similarly, after some number of shear vertical exchange moves to move any
vertical arcs in $\cH_4^-$ past $v^\partial_3$ we can perform shear horizontal exchange
moves so that, in the angular interval $(\theta^\partial_3 , \theta^{v_4}_{\ } )$,
$z^{h_4}_{\ }$ is two horizontal positions above $z^{h_1}_{\ }$.

Now, to determine whether the resulting $\sheararcpresprime$ has a transition
$\sheararcpresprime \btran X^\prime$ with $\b(X^\prime)$
corresponding to a word admitting an elementary flype, we need to inspect
the horizontal position of $h_2$ to see if it is two above minimal over its angular support,
and if the horizontal position of $h_3$ is minimal over its angular support.
(We will deal with the unsatisfactory nature of having to inspect momentarily.)

If the horizontal positions of the boundary are correct
then we will have braid blocks occurring as follows (see Figure \ref{figure:types of flyping disc} $FA^+$):
\bi
\item[$W_1$--] A braid block in $\cH_1^+$ that is under the $R_1 (\subset N)$ disc.
Since there are no available strands in $\cH_1^-$ there is no braiding.
\item[$U$--] A braid block in $\cH_2^-$.  The horizontal position of $h_2$ makes this
a two strand braid block.  Using $\cO_{4}$ we can push all vertical arc in
$\cH_2^+$ that are under $R_2$ past either $v^\partial_1$ or $v^\partial_2$.  
\item[$W_2$--] A braid block in $\cH_3^+$ that is under the $R_3$ disc.
Since there are no available strands in $\cH_3^-$ there is no braiding.
\item[$\emptyset$--]  There is no braiding in $\cH_4^+$ under $R_4$ since any
vertical arcs can be pushed past $v_4$ into $\cH_1^+$. And, there is no
braiding in $\cH_4^-$ since there is but one strand.
\ei

We finally address the issue of having to inspect the
horizontal positioning of $h_2$ and $h_3$. 
As in the exchange move situation, Remark \ref{remark: equivalency class for EM and F} gives us  
that $\sheararcpresprime$ has infinitely many flyping discs.
That is, if $\Delta_f$ illustrates
that a closed braid $X$ admits an elementary
flype---$\b(X) = W_1 U W_2 \sigma^{\pm 1}_{n-1}$---then there are infinitely
flyping discs illustrating the same flype.  Specifically,
for $W_1 , W_2 \in \cW^{n-2}$, $U \in \cU^{n-1}$ and
any word $\a \in \cW^{n-3}$ we have an associated braid/flyping-disc pair
$(X^\prime , \Delta_f (\a))$ with $\b(X^\prime) = [W_1 \a]  U [\a^{-1} W_2] \sigma^{\pm 1}_{n-1}$.
So for $X^\prime \ntran \sheararcpresprime$ we have the family of
resulting pairs $(\sheararcpresprime , \notchdiscprimef (\a))$.

Thus, we have an argument similar to the one for thin exchange moves.
For a positive (respectively, negative) elementary flype we
take a sufficiently complex $\a$ (say powers of full twists on $n-2$ strands)
to `grab' all the stands that are common to the $W_1$ and $W_2$ braid blocks
in the angular support of $h_2$
(respectively, $h_3$).  The winding behavior of leaves in the radial
foliation of $\notchdiscprimef (\a)$ can be used to create a product structure $3$-ball
$B^{1,2}$ (respectively, $B^{2,3}$)
that intersects only these $n-2$ strands and has $h_2 {\rm 's}$ 
(respectively, $h_3 {\rm 's}$) angular support.
Next, using $\cO_4$ we can push out all vertical arcs in $B^{1,2}$/$B^{2,3}$.
Finally, using $B^{1,2}$/$B^{2,3}$ we can reposition that the remaining braiding
occurring in the angular support $h_2$/$h_3$ below all of the
$n-2$ strands intersecting $B^{1,2}$/$B^{2,3}$.
This yields a $\sheararcpresprime$ such that for $\sheararcpres \btran X^\prime$
we have $\b(X^\prime) = W_1^\prime U W_2^\prime \sigma_{n-1}^{\pm 1}$.

\section{Algorithm.}
\label{section:algorithm}
We now describe how to utilize Theorem \ref{theorem:main result} to construct an
algorithm for determining whether a give closed braid admits one of our specified
braid isotopies, i.e. the algorithms of Corollary \ref{corrollary:main result}.
In each case we start with a closed $n$-braid $X$ and readily produce
an arc presentation $\arcpres$.  For each isotopy type the algorithm proceeds by
first specifying a choice of shearing intervals $\cI$. 

\subsection{Algorithm for destabilization isotopy.}
\label{subsection:destabilization isotopy}
For the destabilization isotopy we need only specify the location of one shearing
interval.  But, by the argument of Lemma
\ref{lemma:first simplification of notch disc}, we can chose
any angle $\theta$.  Let
$[\theta - \e,\theta + \e] = \cI$.

Now, we list the possibilities.  For $X$ be a link of $k$-components there are $k$
possible destabilizations, where $k \leq n-1$.
So for a choice of component which has winding around
$\axis$ greater than one we pick an edge-path
$\cE = h_1  \bu v_1 \bu h_2 \bu \cdots \bu h_l$ such that:
1) $h_1 \cap \cI \not= \emptyset \not= h^\eta_l \cap \cI$; and, 2)
$\cE \cap \cI = (h^\eta_1 \cup h^\eta_l) \cap \cI$.

We now apply our elementary moves to $\sheararcpres$ with the proviso that no shearing
occur on the horizontal arcs of $\cE$.

By Proposition \ref{proposition:obvious disc}
after a finite number of elementary moves either there are no more complexity reducing
elementary moves to apply, or we have a resulting edge-path
$\cE^\prime = h^{\prime}_1 \bu v^{\prime}_1 \bu h^{\prime}_2$ such that:
1) $h^{\prime}_1 \cap \cI \not= \emptyset \not= h^{\prime}_2 \cap \cI$; and, 2)
the horizontal positions of $h^{\prime}_1$ and $h^{\prime}_2$ are consecutive.
If it is the latter then we are done.  If it is the former then we start again with a
new edge-path on a different component of $X$.  After a finite number of iterations
this process will stop by either finding a destabilization or determining that
$X$ does not admit a destabilization.

\subsection{Algorithm for thin exchange move.}
\label{subsection:exchange move isotopy}
For the exchange move isotopy we need only specify the location of two shearing
interval which will necessarily have differing edge assignments.
Again, by the argument of Lemma
\ref{lemma:first simplification of notch disc}, we see that we in fact can chose
$\cI = \{ [\theta_1 - \e_1, \theta_1 + \e_1] , [ \theta_2 - \e_2 , \theta_2 + \e_2] \}$ where
over the interval $[\theta_2 , \theta_1]$ there are no vertical arcs.
Specify the interval $[\theta_1 - \e_1, \theta_1  + \e_1]$ as a back edge and
$[ \theta_2 - \e_2 , \theta_2 + \e_2]$ as a front edge.

We consider the possibilities.  They include choices an edge-path $\cE \subset \arcpres$
and an additional horizontal arc $\hat{h}$ such that:
1) $\cE = h_1 \bu v_1 \bu h_2 \bu \cdots \bu h_l $;
2) $h_1$ has angular support containing $[\theta_2 , \theta_1]$;
3) $\cE$ has angular support containing $[\theta_1 , \theta_2]$.
Obviously the number of possible choices for $\hat{h}$ is bounded by the
$n$.  Then our choice of $\cE$ is bounded $n-1$ and must satisfy the condition that
$\hat{h} \cap \cE = \emptyset$.
With a choice of $\cE$ and $\hat{h}$ in hand
we proceed.

We now apply our elementary moves to $\sheararcpres$ with the proviso that no shear
horizontal exchange move
occur on the horizontal portion of $\cE$ and on $\hat{h}$.

By Proposition \ref{proposition:obvious disc}
after a finite number of elementary moves either there are no more complexity reducing
elementary moves to apply, or we have a resulting edge-path
$\cE^\prime = h^{\prime}_1$ such that the angular support of $h^{\prime}_1$
contains the intervals $[\theta_1 , \theta_2]$.
Since the horizontal arc $\hat{h}$ remains unaltered, the angular
support of $\hat{h}$ still contains the interval $[\theta_2 , \theta_1]$.
If it is the latter then we are done.
(We then can use
the horizontal arcs $h^{\prime}_1$ and $\hat{h}$ as the horizontal
boundary of an obvious exchange disc.)
If it is the former then we start again with a
new choice of $\cE$ \& $\hat{h}$ pair.  After a finite number of iterations
this process will stop by either finding an exchange move or determining that
$X$ does not admit an exchange move.

To find thin exchange moves we need to repeat this procedure over all
possible choices of shearing intervals, initial horizontal arcs and edgepaths.  Finally,
over all exchange moves that are recognized by an easy strand count we can pick the ones that
are thin exchange moves.

\subsection{Algorithm for an elementary flype.}
\label{subsection:flype isotopy}
For the flype isotopy we need only specify the location of three shearing
interval.  Again, by Lemma
\ref{lemma:first simplification of notch disc}, we can choose
$\cI $ to be $ \{ [\theta_1 - \e_1, \theta_1 + \e_1] ,
[ \theta_2 - \e_2, \theta_2 + \e_2] , [\theta_3 - \e_3, \theta_3 + \e_3] \}$
where over $[\theta_1 , \theta_2] \bu [\theta_2 , \theta_3]$ there
are no vertical arcs of $\arcpres$.

We concern the final possibilities.  They include choices an edge-path $\cE \subset \arcpres$
and two additional horizontal arcs $\hat{h}_1 , \hat{h}_2$ such that:
1) $\cE = h_1  \bu v_1 \bu h_{2} \bu \cdots \bu h_{l}$;
2) $\cE$ has angular support containing $[\theta_3 \theta_1]$;
3) $\hat{h}_1$ has angular support containing $[\theta_1, \theta_2]$; and,
4) $\hat{h}_2$ has angular support containing $[\theta_2, \theta_3]$.
Obviously the number of possible choices for $\hat{h}_1$ \& $\hat{h}_2$
is bounded by the $\frac{n(n-1)}{2}$---$n$ chosen $2$ ways.
Our choice of $\cE$ is bounded
by $n-2$.  With a choice of $\cE$ and
$\hat{h}_1 , \hat{h}_2$ in hand we continue.

We now apply our elementary moves to $\sheararcpres$ with the proviso that no shear
horizontal move
occur on the horizontal portion of $\cE , \ \hat{h}^\eta_1 , \ {\rm or} \ \hat{h}^\eta_2$.
Moreover, we require that side on which shear horizontal exchange move or shear
vertical simplification be consistent with the edge assignment of our shearing
intervals as specified in \S \ref{subsection:destabilization isotopy}
and \ref{subsection:exchange move isotopy}.

By Proposition \ref{proposition:obvious disc}
after a finite number of elementary moves either there are no more complexity reducing
elementary moves to apply, or we have a resulting edge-path being
$\cE^\prime = h^{\prime}_1 \bu v^{\prime} \bu h^{\prime\prime}_1$
such that the angular support of $\cE^\prime$
contains the interval $[\theta_3 , \theta_1]$.
Since the horizontal arc $\hat{h}$ remains unaltered, the angular
support of $\hat{h}_1$ \& $\hat{h}_2$ still contain the appropriate intervals of $\cI$.
(We then can use the horizontal arcs $\cE^\prime , \hat{h}_1 , \hat{h}_2$
as the horizontal boundary of an obvious flyping disc.)
If it is the latter then we are done.  If it is the former then we start again with a
new choice of triple $\cE ,\hat{h}_1 , \hat{h}_2$.
After a finite number of iterations
this process will stop by either finding an elementary flype or determining that
$X$ does not admit a flype.
Reiterating this procedure over all possible choices (shearing intervals, initial
horizontal arcs and edgepaths) we can recognize all possible braid preserving flypes.

\subsection{Proof of Theorem \ref{Theorem:destably equivalent}}
\label{subsection:destably equivalent}
Since we now have an algorithm for when a closed braid admits a
thin exchange move or elementary flype we take note
of the finite possibilities for a closed braid $X$.
Referring back to the ``connect the dots''discussion prior
to the statement of Theorem \ref{theorem:main result},
the finiteness of $\cC(\sheararcpres)$ tells us that there are only
finitely many possible thin exchange moves or elementary flypes to perform.
Performing a thin exchange move or elementary flype will result in a
$X^\prime$ braid.
We then use the classical solutions to the conjugacy problem
to see whether $X^\prime$ and a given closed braid $Y$ are braid isotopic.
For thin exchange move we are only checking to see if $Y$ is conjugate
to the appropriate closed braid associated with $X^\prime = W \tau_{[s , t+1]}^{\ } U \tau_{[s , t+1]}^{-1}$,
i.e. are $X$ and $Y$ related by a {\em single} thin exchange move.
As an aside, the case where $Y = X$ is of particular interest.
In a fashion similar to the  parallel strands in a double destabilization,
A. V. Malyutin's work points to the possibility that $WU$ may be conjugate to
$W \tau_{[s , t+1]}^{k} U \tau_{[s , t+1]}^{-k}$ for any integer value of $k$.

To determine whether $X$ admits a double destabilization we first apply our
exchange move algorithm to
determine whether $X$ admits a thin exchange move with $t-s = 2$ and $U \in \cU^{n-3}$.
If $X$ admits such a thin exchange move then over all such thin exchange moves,
we apply our destabilization algorithm to determine if the $U$ block of $X = WU$
admits a destabilization.  If $U$ admits a destabilization, we do the
destabilization and applying our destabilization once more to the resulting
$U$-block.  If the resulting $U$-block also destabilizes then by inspection
we determine whether $U$ was conjugate to $\UU$. 

\begin{figure}[p]
\centerline{\includegraphics[scale=.75, bb=0 0 363 143]{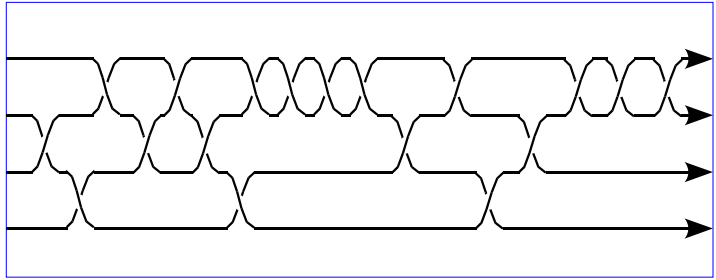}}
\caption{The rectangular disc is a portion of $C_1$ where we identify the left and the
rights edges.  With this representation of $C_1$ we then have a standard projection of
a $4$-braid.}We
\label{figure:example-01}
\end{figure}

\subsection{Destabilization example.}
\label{subsection:destabilization examples}
\begin{figure}[p]
\centerline{\includegraphics[scale=.75, bb=0 0 363 146]{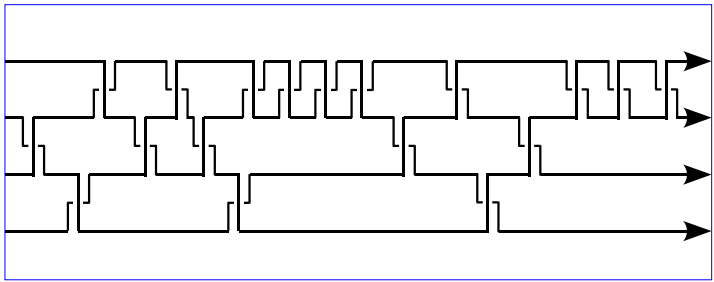}}
\caption{The arc presentation is directly obtained from the standard projection
in Figure \ref{figure:example-01}.}
\label{figure:example-02}
\end{figure}

We now give an extended example illustrating how the algorithm is implemented on a
closed braid that admits a destabilization.  Figure \ref{figure:example-01} is a standard projection of a $4$-braid, $X$, drawn on a portion of $C_1$.
($C_1$ is represented as a rectangular disc with its left and right edges identified.)
We will step through the implementation of the algorithm by first performing the transition
$X \ntran \arcpres$ which gives us the arc presentation illustrated in Figure
\ref{figure:example-02}.  The complexity measure of the arc presentation in
Figure \ref{figure:example-02} is just the number of vertical arcs, i.e. $54$.

\begin{figure}[p]
\centerline{\includegraphics[scale=.75, bb=0 0 363 143]{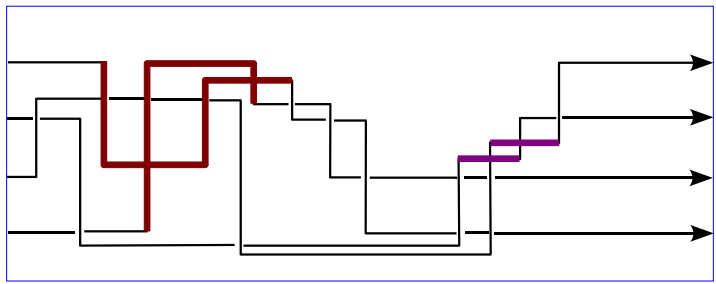}}
\caption{We indicate a vertical-horizontal-vertical (red) edge-path upon which we will
perform a horizontal exchange move; a vertical-horizontal-vertical-horizontal
(red) edge-path upon which we will perform two horizontal simplifications; and two
(magenta) horizontal arcs upon which we will perform two more horizontal exchange moves.}
\label{figure:example-03}
\end{figure}

The arc presentation in Figure \ref{figure:example-02} clearly has numerous places where
the number of arcs can be reduced by the application of either vertical or horizontal
simplifications.  Without altering the number of crossings in the arc presentation
we apply a number of such simplifications until we obtain the arc presentation
illustrated in Figure \ref{figure:example-03}.  The complexity measure of $\arcpres$ is
now $13$.

\begin{figure}[h]
\centerline{\includegraphics[scale=.75, bb=0 0 436 146]{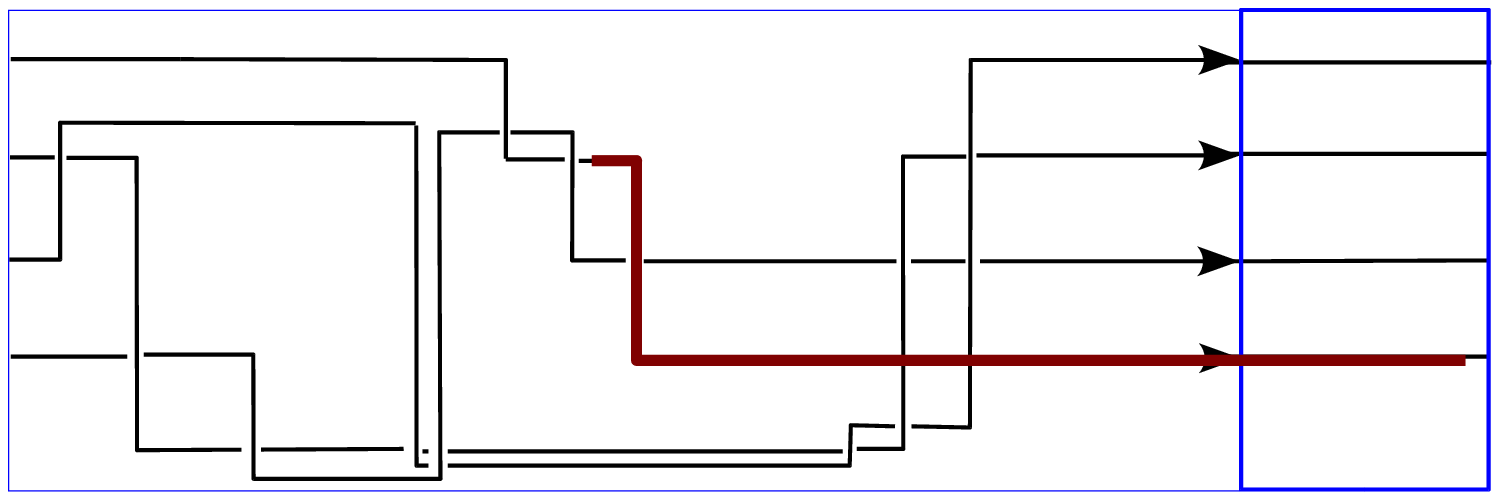}}
\caption{The right extended portion of the rectangular disc is our shearing interval.
The red edge-path will be used in a shear horizontal exchange move.}
\label{figure:example-04}
\end{figure}

Next we perform, first a horizontal exchange move.  The portion of the arc presentation
where this horizontal exchange move occurs is indicated by a transparent heavy red edge-path
that is a vertical-horizontal-vertical path.  After the performance of this exchange move
it is possible to perform two horizontal simplifications.  The portion of the arc
presentation involved in these simplification is indicated by a transparent heavy red
edge-path that is a vertical-horizontal-vertical-horizontal path.  These two
horizontal simplifications reduced the complexity of $\arcpres$ to $11$.
Finally, we perform two more horizontal exchange moves that can be sequenced
independent of our previous operations.  They are the two horizontal arc that are
indicated by the heavy/transparent magenta horizontal arcs.  These operations
result in the arc presentation illustrated in Figure \ref{figure:example-04}.

\begin{figure}[htpb]
\centerline{\includegraphics[scale=.75, bb=0 0 433 146]{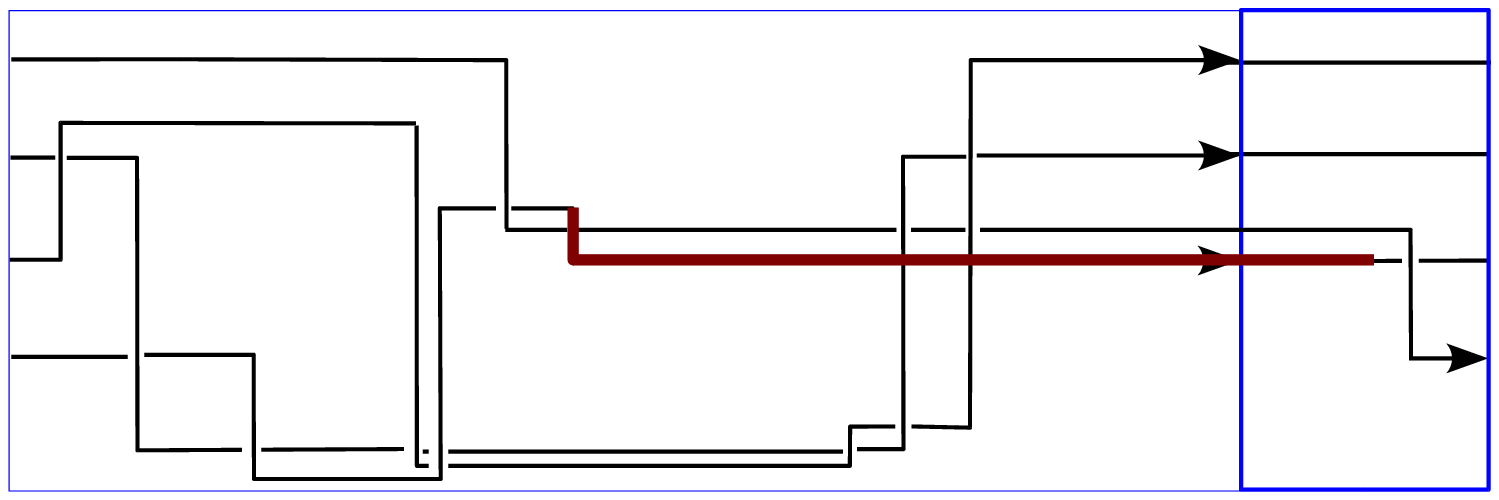}}
\caption{The red edge-path will be used in a shear vertical simplification.}
\label{figure:example-05}
\end{figure}

We now ``extend'' the right side of the rectangular disc so as to include a front edge
shearing interval and we consider $\sheararcpres$.  The complexity of $\sheararcpres$
is initially the same as that of $\arcpres$, namely $11$.
We now perform a shear horizontal exchange move and then a horizontal
simplification using the edge-path in Figure \ref{figure:example-04} that is
indicated by the heavy/transparent red path.  This will yield the arc presentation
illustrated in Figure \ref{figure:example-05}.
(This could also be seen as the performance of a shear vertical
simplification.)  The complexity of the resulting $\sheararcpres$ is reduced by one
to $10$.

\begin{figure}[htpb]
\centerline{\includegraphics[scale=.75, bb=0 0 433 146]{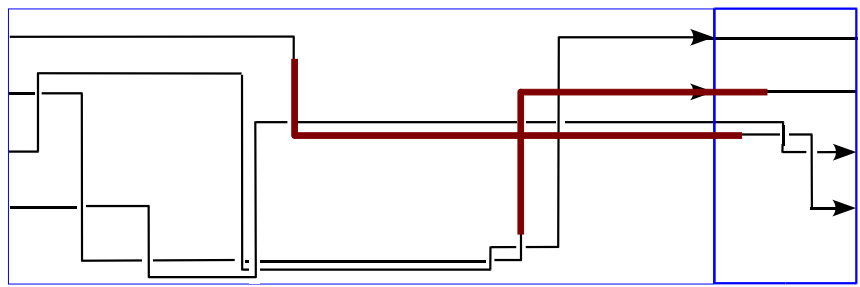}}
\caption{We indicate again by red edge-paths where we will perform two shear
horizontal exchange moves.}
\label{figure:example-06}
\end{figure}

The alteration from Figure \ref{figure:example-05} to \ref{figure:example-06} is
achieved by a shear vertical simplification applied to the vertical-horizontal
red edge-path.  The complexity measure of $\sheararcpres$ is reduced by a count of
one to $9$.

On Figure \ref{figure:example-06} we perform two move shear horizontal exchange moves
which will leave the complexity measure of $\sheararcpres$ unchanged.  The edge-paths
where the shear horizontal exchange moves occur are indicated in heavy/transparent
red again.  The operation leaves unchanged the complexity of $\sheararcpres$.

\begin{figure}[htpb]
\centerline{\includegraphics[scale=.75, bb=0 0 433 146]{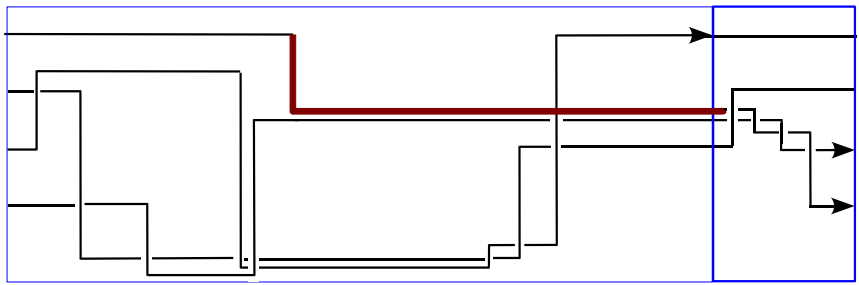}}
\caption{The red edge-path indicates where we will perform one last shear
horizontal exchange move to achieve the obvious destabilizing arc presentation
diagram in Figure \ref{figure:example-08}.}
\label{figure:example-07}
\end{figure}

\begin{figure}[htb]
\centerline{\includegraphics[scale=.75, bb=0 0 433 146]{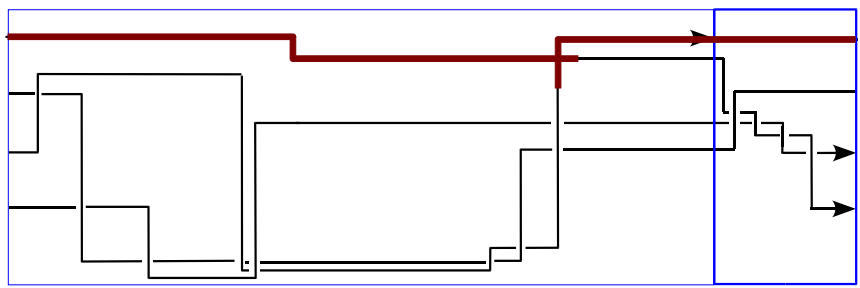}}
\caption{The red heavy/transparent edge-path picks out the configuration for
an obvious destabilization.}
\label{figure:example-08}
\end{figure}

Finally, using the edge-path indicated by the red heavy/transparent path we perform
one last shear horizontal exchange move to obtain the obvious destabilizing
configuration illustrated in Figure \ref{figure:example-08}.  This last
operation leaves the complexity of $\sheararcpres$ unchanged.


\subsection{Final remarks.}
\label{subsection:final remarks}
It is now useful to engage is some informed speculation on reasonable directions to
push the machinery used to establish Theorem \ref{theorem:main result}.  
The most reasonable direction is
in adapting our machinery to recognize other isotopies of the MTWS calculus.  Briefly, the
statement of the MTWS says that of a given pair of positive integers $(m,n)$ with
$m\geq n$ there are only finitely many {\em templates}---isotopic pairs of block-strand
diagram---which need be used to carry all closed $m$-braids to corresponding closed
$n$-braids.  The block-strand diagrams in Figure \ref{figure:moves} illustrate the three
simplest templates.  In general, in order to recognize when a braid is carried by a
specified block-strand diagram we need to account for both the blocks and the braiding
of the strands.  Adapting our $\notchdisc$-disc machinery we can clearly add more
shearing intervals to our set $\cI$.  Since the angular support of any block
can be made arbitrarily small, for a sufficient number of shearing intervals we
can always place a $\notchdisc$ in the complement of the braid the ``jumps'' from
strand to strand as its boundary winds around the axis once.  Moreover,  we can allow
the usage of multiple disjoint $\notchdisc$-discs so that there boundaries correspond to
an unlink.  Our machinery does not require that we have a $\notchdisc$-disc of a single
component.

\begin{figure}[htpb]
\centerline{\includegraphics[scale=.75, bb=0 0 306 239]{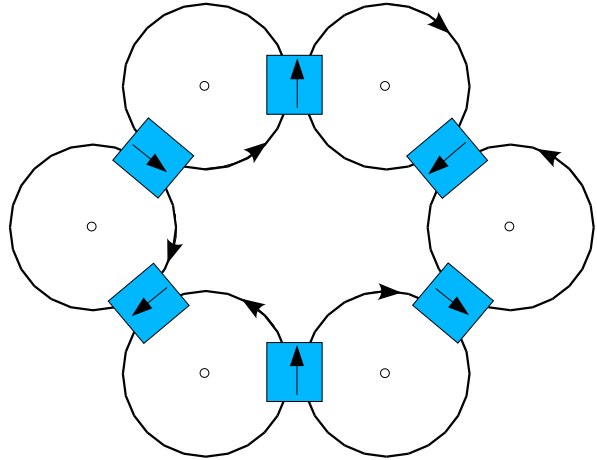}}
\caption{The $6$ block can contain any braiding that is consistent with the direction
of the indicate orientations.  This block-strand diagram can be re-embedded in
a braid fibration such that it is knotted.  There will then exist a non-peripheral
essential torus that will be a type $2k$ torus as described in \cite{[BM3]}.  We can
use this torus to alter the braid presentation so as to recognize the admission
of a cyclic move.}
\label{figure:cyclic move}
\end{figure}

Finally, we do not have to restrict ourselves to just surfaces that are discs.  Depending
on the configuration of the block-strand diagrams used in the isotopy template it may
be more appropriate to use a closed surface in the complement of closed braid.  For example,
in \cite{[BM4]} {\em cyclic moves} were analyzed.  In Figure \ref{figure:cyclic move}
we illustrate a block-strand diagram that has $6$ blocks attached to $6$ circles.  This
block-strand diagram can be embedded in braid fibration $\fib$ such that it is in fact
knotted.  Thus, we can think of this block-strand diagram contained in a solid torus
and we knot this diagram by knotting the solid torus.  The embedding of the boundary
of this solid torus will be a standard type $2k$ embedding as describe in
\cite{[BM3]}.  We can then use the arc presentation machinery along with
the coning discs to alter the arc presentation through horizontal and vertical
exchange moves until our type $2k$ torus has a standard tiling.

Based upon this discussion there should be optimism about the possibility of
recognizing any specified isotopy of the MTWS calculus.  But, there is a cautionary point.
Once it has been determined that a given closed braid admits a specified isotopy there may
still be an ambiguity on how to apply the isotopy.  For a closed braid that admits a
destabilization, it is clear that one need only destabilize.  But, for a closed braid
that admits an exchange move, it is not clear which way to apply the exchange move to
reduce.  (See Figure 6 of \cite{[BM4]} for an example of this pathology.)  On this point more investigation is needed.

\end{document}